\documentclass{amsart}
\usepackage{graphicx} 
\usepackage{amsmath}
\usepackage{a4wide}
\usepackage{amsthm}
\usepackage{amsfonts}
\usepackage{amssymb}
\usepackage{mathdots}
\usepackage{tikz}
\usepackage{url}
\usepackage{quiver}

\newcommand{\att}[2]{\raisebox{-.5\height}{ 
\includegraphics[scale = #2]{diagrams/#1.pdf}}}

\usepackage[boxsize = 1em]{ytableau}

\usepackage{pdfpages}
\usepackage{fancyhdr}
\usepackage{parskip}

\newtheorem{definition}{Definition}
\newtheorem{bigDefinition}{Definition}
\newtheorem{theorem}{Theorem}
\newtheorem{proposition}{Proposition}
\newtheorem{lemma}{Lemma}

\newtheorem{corollary}{Corollary}
\newtheorem{remark}{Remark}
\newtheorem{example}{Example}

\newtheorem{bigTheorem}{Theorem}

\newcommand{\directSumOver}[1]{\bigoplus_{\begin{substack}{#1}\end{substack}}}

\newcommand{\inner}[2]{\left\langle #1 , #2 \right\rangle}

\newcommand{\I}{\cap}
\newcommand{\U}{\cup}
\newcommand{\comp}{\setminus}

\newcommand{\fanB}{\mathcal{B}}
\newcommand{\fanC}{\mathcal{C}}
\newcommand{\fanD}{\mathcal{D}}
\newcommand{\C}{\mathbb{C}}

\newcommand{\N}{\mathbb{N}}
\newcommand{\Z}{\mathbb{Z}}

\newcommand{\R}{\mathbb{R}}

\newcommand{\inv}[1]{{#1}^{-1}}
\newcommand{\g}{\mathfrak{g}}
\newcommand{\h}{\mathfrak{h}}
\newcommand{\blank}{{-}}

\newcommand{\aff}[1]{\mathfrak{\tilde{#1}}}

\newcommand{\tensor}{\otimes}

\newcommand{\sms}[1]{\overline{#1}}

\newcommand{\V}{\mathcal{V}}
\newcommand{\so}{\mathfrak{so}}
\newcommand{\fsp}{\mathfrak{sp}}

\newcommand{\sbs}{\subset}
\newcommand{\set}[1]{\left\{#1\right\}}
\newcommand{\defset}[2]{\left\{#1 : #2\right\}}

\newcommand{\gen}[1]{\left\langle#1\right\rangle}

\title{Clifford BMW categories And Conformal Embeddings of Types B C and D}
\author{August Bergquist}

\begin{document}

\begin{abstract}

In this paper, we give diagrammatic descriptions of the quantum subgroups coming from conformal embeddings of the form $\V(\g, \check{h}) \sbs \V(\so_{\dim \g}, 1)$ for $\g \in \set{\so_N, \fsp_N}$, as well as the conformal embeddings of the form $\V(\fsp_{2r + 2}, 2r) \sbs \V(\so_{\dim \fsp_{2r}}, 1)$, extending the techniques of Edie-Michell and Snyder which applied to conformal embeddings of the form $\V(\mathfrak{sl}_r, r) \sbs \V(\so_{\dim \mathfrak{sl}_r} , 1)$. Our main conceptual innovation is the use of fermion condensation to construct a continuous family of planar superalgebras, which we call the Clifford BMW supercategories. The tower of endomorphism algebras, $CW_n(\sigma, \epsilon)$ is a combination of the Birman Wenzl Murakami algebras and the Clifford algebras, much in the spirit of the Hecke Clifford algebras. This tower of algebras appears to be new, and already yields some exciting combinatorics. 
\end{abstract}

\maketitle
\section{Introduction}

This paper establishes the existence of a new interpolation category, which we call the Clifford BMW supercategory. This category extends the Dubrovnik\footnote{Distinct from the Kauffman skein theory in that $\att{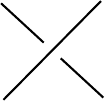}{0.25} - \att{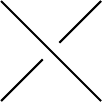}{0.25} = z\left(\att{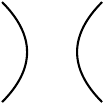}{0.25} - \att{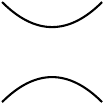}{0.25}\right)$} skein theory by a single odd dot which acts as a spinor. The endomorphism algebras, which we call the Clifford BMW superalgebras, are interesting in their own right. We believe these algebras will have consequences beyond the study of tensor categories. The representation theory of these algebras appears to be tractable directly from the presentation, without reference to the larger story in which they first appeared. For the reader who is more interested in towers of algebras and algebraic combinatorics, this paper serves to introduce these algebras and show that their relations are consistent, all while revealing much of their combinatorial structure. 

The Clifford BMW supercategory interpolates over certain infinite series of quantum subgroups. 

Quantum subgroups form a rich variety of examples of fusion categories, many of which come from conformal embeddings. One takes an Etale algebra $A$ in $C = \mathcal{C}(\mathfrak{g}, k)$ coming from a conformal embedding $\V(\mathfrak{g}, k) \sbs \V(\mathfrak{h}, 1)$, giving the category of modules $C_A$. Quantum subgroups coming from conformal embeddings have been fully classified in \cite{Schellekens}, but until recently their internal structure has been little understood. Throughout the last decade, diagrammatic techniques have made significant headway in this area.

The first diagrammatic description of a quantum subgroup coming from a conformal embedding was given in by Xu in \cite{E8} for the conformal embedding 
\[
\V(\mathfrak{sl}_2, 28) \sbs \V(\mathfrak{g}_2, 1). 
\]
His diagrammatic presentation gave the planar algebra E8, which involves a six-valent skein theory.

The first complete diagrammatic description of an infinite family of quantum subgroups coming from conformal embeddings was given by Liu in \cite{ZhengweiLiu}. His description applied to the family of conformal embeddings of the form 
\[
\V(\mathfrak{sl}_{N}, N \pm 2) \sbs \V(\mathfrak{sl}_{N(N \pm 1)/2}, 1),
\]
as shown in \cite{CainAndJacques}, Theorem 2.10.
This description also allowed for a Deligne interpolation result, revealing a continuous family of planar algebras which ``interpolated'' over these conformal embeddings.

In addition to the diagrammatic descriptions of these infinite families, in \cite{Caleb} Hill has applied graph planar algebra techniques to give a diagrammatic description of the conformal embeddings 
\[
\V(\g_2, 4)\sbs \V(\mathfrak{so}_{\dim\g_2}, 1)
\]
as well as 
\[
\V(\g_2, 3) \sbs \V(\mathfrak{e}_6, 1).
\] 
The techniques derived in \cite{Caleb} likely extend to the other finite families of conformal embeddings which give connected Etale algebras.


Among the conformal embeddings classified in \cite{Schellekens} and listed in \cite{WittGroup}, a large class of these fit into the pattern
\[
\mathcal{V}(\g, \check{h}) \sbs \V(\so_{\dim \g}, 1).
\]
Such conformal embeddings are called adjoint conformal embeddings. Related to these, when $\g = \so_{N}$ or $\g = \fsp_{2r}$, one also has a corresponding infinite family derived from these by level rank duality.

In \cite{CainAndSnyder} Edie-Michell and Snyder gave a diagrammatic description of the family of conformal embeddings of the form
\[
\mathcal{V}(\mathfrak{sl}_N, N +1) \sbs \mathcal{V}(\mathfrak{so}_{\dim \mathfrak{sl}_N}, 1).
\]
Their description yielded another Deligne interpolation result, establishing the consistency of a continuous family of Hecke-Clifford categories, as well as certain extensions of these categories at specified roots of unity. Subsequently in \cite{CainAndHans}, Edie-Michell and Wenzl leveraged this description to give a uniform parametrization of the simple objects, fusion rules, and the quantum dimensions in each of these categories. Edie-Michell and Snyder conjectured the existence of interpolation categories for the remaining infinite families: a result which this paper establishes.

In this paper, we provide diagrammatic descriptions for the remaining infinite families of conformal embeddings of the form $\mathcal{V}(\g, \check{h}) \sbs \V(\so_{\dim \g}, 1)$, excluding modules involving spinor representations in the type B and D cases. By a simple adjustment, we also address the conformal embeddings of the form $\V(\fsp_{2r + 2}, 2r) \sbs \V(\so_{\dim \fsp_{2r}}, 1)$. Our techniques almost certainly apply directly to the infinite families of conformal embeddings in types B and D derived from adjoint conformal embeddings by level rank duality, although we have omitted these cases in this paper.

The main difference between our techniques and those of \cite{CainAndSnyder} is the use of fermion condensation, without which the non-oriented nature of the Dubrovnik skein theory becomes unwieldy. The same techniques in this paper can also be used to significantly simplify their presentation, although their evaluation algorithm for the Hecke Clifford category is fairly straightforward. This paper could be interpreted as showing that the most natural language for diagrammatically treating quantum subgroups coming from adjoint conformal embeddings is actually super-pivotal categories in the sense of \cite{SuperPivotalCats}. This eliminates the need for a two colored skein theory, reduces the number of relations, and shows how many of the relations in the even subcategories are simple consequences of the axioms of a super-pivotal category.

 For each of the conformal embeddings of the form $\V(\g, k) \sbs \V(\so_{\dim \g}, 1)$ for $\g \in \set{\fsp_{2r}, \so_{2r}, \so_{2r + 1}}$, the distinguished simple $V_{(1)}\in \mathcal{C}(\g, k)$ is sent to a simple object $X = F_A(V_{(1)})\in \mathcal{C}(\g, k)_A$.
 Moreover, in each case there is a local module $g\in \mathcal{C}(\g, k)_A$, with the property that 
 \[
g \tensor_A g \cong A.\footnote{$A$ is the tensor unit in $\mathcal{C}(\g, k)_A$}
 \]

 This local module $g$ lifts to a fermion object in the Drinfeld center with a half braid $\gamma$, allowing us to form the fermion condensation $C_A/\gamma$. The fermion $g$ has the property that $g \tensor X \cong X$. As a result, the object $X$ is a type $Q$ simple object in $C_A/\gamma$, with an additional odd automorphism $v : X \to X$. 
 
 We shall denote the identity map on $X$ by a black strand, both in $C_A$ and $C_A/\gamma$. Diagrammatically, we denote this odd automorphism $v: X \to X$ by a dot:
 \[
 v = \att{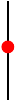}{0.5}
 \]

The composition of two odd morphisms is even, and the space of even endomorphisms of $X$ is $\C$. Therefore, with an appropriate normalization,
\[
\att{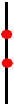}{0.5} = \att{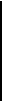}{0.5}.
\]

Because the odd part of $Hom_{C_A}(X\to X)$ is one dimensional and $X$ is self dual, we must have $v^* = \lambda v$ for some scalar $\lambda$. Moreover, since $v$ is odd and $C_A/\gamma$ is super pivotal, we must have $v^{**} = - v = \lambda^2 v$. As a result, $\lambda = \pm i$. The particular sign choice is entirely determined by the twist parameter $\theta$ so that $\att{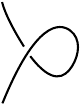}{0.25} = \theta \att{stick}{0.25}$ in $C$.

Because the free module functor $F_A: C \to C_A$ is central, the usual braiding on $V_{(1)}$ in $C$ becomes a half braid on the object $X$ in $C_A$. Diagrammatically, this means that the first and second Reidemeister moves apply, as well as an additional relation 
\[
\att{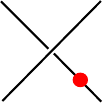}{0.5} = \att{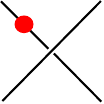}{0.5}
\]

Geometrically, these relations can be understood by introducing a back wall with a spin structure, where tangles interact with the back wall at finitely many points, representing instances of the morphism $v$. For example:
 \[
    \att{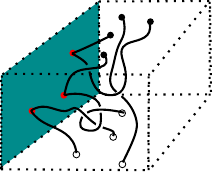}{1}
    \]
The tangles in the three dimensional half space are subject to the Dubrovnik relations, while the spinor dots are subject to the relations described above as well as the super interchange rule. This interpretation can be made even more vivid by imagining that the spinor dots pick up a phase shift as isotopies of the three dimensional half space deform the back wall. 

Certain calculations reveal that the above relations can only be consistent when $\att{twist}{0.5} = \sigma i q^{-\epsilon}$ for $\sigma, \epsilon \in \set{\pm 1}$, where the first sign choice is dictated by the spin structure on the back wall. This explains why the Clifford BMW supercategory has only a single continuous parameter instead of two.

The relations just described suggest one ought to define the following category for generic $q$, over some localization $R$ of $\C[q]$.

\begin{bigDefinition}
    Let $R = \C[q, \inv{q}, \frac{1}{q - \inv{q}}, \frac{1}{(1 + i)q + (i - 1)\inv{q}}, \frac{1}{(1 + i)\inv{q} + (i - 1)q}]$. The Clifford $BMW$ supercategory, denoted $CW(\sigma, \epsilon)$ is the quotient of the free pivotal $R$-linear supercategory generated by the $BMW$ generators subject to the $BMW(q, \sigma i q^{-\epsilon})$ relations, all of which are declared to be even, along with a single odd morphism $\att{oddDot}{0.5}: 1 \to 1$, subject to three additional relations

$$
\begin{array}{ccc}
   \text{(Rotate) } \left(\att{oddDot}{0.5}\right)^* = -\sigma i \att{oddDot}{0.5} & \text{(Half Braid) } \att{oddHalfBraid1}{0.5} = \att{oddHalfBraid2}{0.5} & \text{(Cancel) } \att{doubleDot}{0.5} = \att{stick}{0.5}
\end{array}
$$
\end{bigDefinition}

The reason for our particular choice of the ring $R$ is as follows. Since we will be uniting all of the quantum subgroups of types B, C, and D by interpolating over infinitely many specializations of $q$, the ring $R$ over which we define $CW(\sigma, \epsilon)$ cannot invert every polynomial in $q$. On the other hand, in order define the category, we must invert $q$. In order to make many of our results work, we also need to invert a number of other rational functions. By inverting little bit more polynomials than necessary, the $R$-linear planar superalgebras $CW(\sigma, \epsilon)$ are conjugate under symmetries of $R$ for any choices of $\sigma$ and $\epsilon$. This reveals deep regularities between all of the quantum subgroups treated in this paper.

In a sense, we regard the Clifford BMW supercategories as more natural algebraic objects than their even subcategories $CW(\sigma, \epsilon)[0]$. Before realizing that fermion condensation would bring massive simplifications to this project, we had originally only worked with $CW(\sigma, \epsilon)[0]$, to which we had given the same name. Due to the non-orientedness of the strands, the evaluation algorithm and basis theorem were extraordinarily difficult to prove, and the resulting two-colored skein theory was significantly less geometrically motivated. Worse still, the even endomorphism algebras appeared to be practically impossible to present.

The usual spanning set for the BMW algebras can be modified to accommodate the Clifford BMW skein theory.

 \begin{bigTheorem}
     For all $\sigma, \epsilon\in \set{\pm 1}$ and for all $n\in \N$, $End_{CW(\sigma, \epsilon)}(n)$ is an $R$-module of rank at most $2^{n}(2n - 1)!!$. A spanning set is indexed by pairings of the boundary points, and a choice of whether to place a spinor decoration on each of the strands. 
 \end{bigTheorem}

 In particular, the skein relations for the Clifford BMW category are evaluable, which is essential for proving the following equivalences of categories. 
 
\begin{bigTheorem}
    Let $r > 1$. Let $A\in \mathcal{C}(\fsp_{2r}, r + 1)$ be the Etale algebra coming from the conformal embedding $\V(\fsp_{2r}, r + 1) \sbs \V(\so_{\dim \fsp_{2r}}, 1)$. With $q = e^{\frac{2\pi i}{4(2r + 2)}}$ there exists a full and dominant even super-pivotal functor $\Phi: CW_q(-1, 1)\to  \mathcal{C}(\fsp_{2r}, r + 1)_A)/\gamma$.
    
    The even part
    $\Phi[0]: CW_q(-1, 1)[0]\to  \mathcal{C}(\fsp_{2r}, r + 1)_A$ gives a pivotal equivalence
    \[
    Ab(\sms{\Phi[0]}):Ab(\sms{CW_q[0](-1, 1)} \simeq \mathcal{C}(\fsp_{2r}, r + 1)_A.
    \]
\end{bigTheorem}

Level rank duality gives a mirror result for the Etale algebra $A\in \mathcal{C}(\fsp_{2r + 2}, 2r)$ coming from the conformal embedding $\V(\fsp_{2r + 2}, 2r) \sbs \V(\so_{\dim \fsp_{2r}}, 1)$. 

By interpolating over all of these conformal embeddings, we are able to establish the following basis theorem, showing that our spanning set is linearly independent.

\begin{bigTheorem}
    For all $\sigma, \epsilon\in \set{\pm 1}$ and for all $n\in \N$, $End_{CW(\sigma, \epsilon)}(n)$ is a free $R$-module of rank $2^{n}(2n - 1)!!$, and for every specialization of $CW(\sigma, \epsilon)$ to any $q\in \C$ such that $q^{8} \ne 1, 0$, $\dim CW_q( \sigma, \epsilon) = 2^{n} (2n - 1)!!$.
\end{bigTheorem}

We also obtain faithful functors from $\sms{CW_q(\sigma, \epsilon)}$ into the fermion condensations of quantum subgroups of types B and D. The parameters for which this happens are given in the following table:
 
\begin{center}
\begin{tabular}{c|c|c|c}
    Conformal Embedding & $q$ & $\epsilon$ & $\sigma$ \\\hline
    $\V(\fsp_{2r}, r + 1) \sbs \V(\so_{\dim \fsp_{2r}}, 1)$ & $e^{\frac{2\pi i}{4(2r + 2)}}$ & $1$ & $-1$  \\\hline
    $\V(\so_{2r + 1}, 2r - 1) \sbs \V(\so_{\dim \so_{2r + 1}}, 1)$ & $e^\frac{2\pi i}{4(4r - 2)}$ & $-1$ & $1$ \\\hline
    $\V(\so_{2r}, 2r - 2) \sbs \V(\so_{\dim \so_{2r}}, 1)$ & $e^{\frac{2\pi i}{4(4r - 4)}}$& $-1$& $1$   \\\hline
    $\V(\fsp_{2r + 2}, 2r) \sbs \V(\so_{\dim \fsp_{2r}}, 1)$ &$e^{\frac{2\pi i}{4(2r + 2)}}$ & $-1$ & $-1$ 
\end{tabular}
\end{center}

These functors, while faithful, are only full for the conformal embeddings of type C. In order to get equivalences of categories for the adjoint conformal embeddings of types B and D, we extend the specializations of the Clifford BMW categories at the appropriate values of $q$, $\epsilon$, and $\sigma$ by an extra vertex.

\begin{bigDefinition}
    Let $q^{8} \ne 1, 0$ and $N > 2 \in \N$ and $\theta\in \C \comp \mathcal{S}$, and $\sigma, \epsilon\in \set{\pm 1}$, such that the idempotents $P_\lambda$ in $BMW(q, \theta)$ are defined when $\lambda^t_1 +\lambda^t_2 \le N$. Define the \textbf{Special Clifford BMW category}, denoted $SCW_N(q)$, to be $\C$-linear pivotal category with symmetrically self-dual objects $\N$, which extends $CW_q(1, -1)$ by the generator 
    \[
    \att{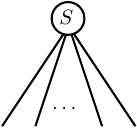}{0.5}\in Hom_{SCW_N(q)}(N \to 0)
    \]
    subject to the relations
    \[
    \text{(Braid)} \att{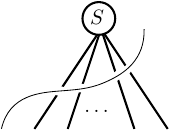}{0.5} = \att{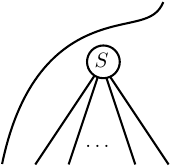}{0.5} \hspace{0.5cm} \text{(Split)} \att{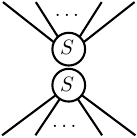}{0.5} = P_{(1^N)}
    \]
    where $P_{(1^N)}$ is the idempotent in $BMW(q, \theta)$ as in section~\ref{BMWIntro}, as well as the (rather awful) additional relation when $N = 4k$ for some integer $k$,
    \[
    \att{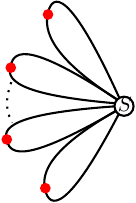}{1} = \lambda
    \]
    where the scalar $\lambda$ (which is rather ugly) is entirely determined by the other relations up to a sign choice.
\end{bigDefinition}

The additional relation \[
    \att{dotsOnPedels}{1} = \lambda,
    \]
while horrible, luckily appears for only half of the type $D$ conformal embeddings (those of the form $\mathcal{V}(\so_{4k}, 2k - 2)\sbs \mathcal{V'}$). 

This extension gives a presentation for the quantum subgroup of types B and D coming from adjoint conformal embeddings. 
\begin{bigTheorem}
    For every $N\in \N$ with $N> 2$ there is a pivotal equivalence
    \[
    Ab(\sms{SCW_N(q)}) \simeq \mathcal{C}(\so_N, \check{h})^{vec}_A
    \]
    where $\check{h} = \begin{cases}
        2r - 1 & N = 2r + 1\\
        2r - 2 & N = 2r
    \end{cases}$ and $q= \begin{cases}
        e^\frac{2\pi i}{4(4r - 2)} & N = 2r +1\\
        e^\frac{2\pi i}{4(4r - 4)} & N = 2r.
    \end{cases}$
    where $A$ is the Etale algebra coming from the conformal embedding $\V(\so_N, \check{h}) \sbs \V(\so_{\dim \so_N}, 1)$.
\end{bigTheorem}

The endomorphism algebras for $CW(\sigma, \epsilon)$ are easily presented, in contrast with their even subalgebras which we have been unable to find a presentation for. We doubt it can be done in a reasonable way.

\begin{bigTheorem}
    The algebra $CW_n(\sigma, \epsilon)$ has presentation 
\[
CW_n(\sigma, \epsilon) \cong \gen{v_k, t_\ell, e_\ell:\ell \in \set{1, \dots, n - 1}, k\in \set{1, \dots, n}}/R
\]
where $R$ is the ideal generated by the following relations
\begin{itemize}
    \item[BMW] For $\ell \in \set{1, \dots, n - 1}$ the generators $t_\ell$ and $e_\ell$ satisfy the relations of the Birman Wenzl Murakami algebra $BMW(q, \sigma i q^{-\epsilon})$
    \item[C] The $v_\ell$'s for $\ell\in \set{1, \dots, n}$ satisfy the relations of the Clifford algebra $Cl_n$, namely
    \begin{itemize}
        \item[C1] $v_\ell^2 = 1$ for all $\ell$, 
        \item[C2]  $v_\ell v_k = -v_k v_\ell$ for $\ell \ne k$
    \end{itemize}
    \item[M] As well as the mixed relations
    \begin{itemize}
        \item[M1] $v_{\ell + 1} t_\ell = t_\ell v_\ell$
        \item[M2] $v_\ell e_\ell = \sigma i v_{\ell + 1}e_\ell $ and $e_\ell v_\ell = -\sigma i e_\ell v_{\ell + 1}$
        \item[M3] $v_\ell e_k = e_\ell v_k$ when $|\ell - k| > 1$
        \item[M3] $v_\ell t_k = t_k v_\ell$ when $|\ell - k| > 1$
    \end{itemize}
\end{itemize}
\end{bigTheorem}

The Clifford BMW superalgebras appear to be new, and already appear to yield rich combinatorial structure. Morally, they give a ``wreath product'' of the Brauer algebra with the Clifford algebra, although there are good reasons why this cannot be done in a literal sense.

From this presentation, and by establishing the non-degeneracy of the trace by again interpolating over conformal embeddings of type C, we are able to establish the following decomposition for each $n$, for all but finitely many $q\in \C$.\footnote{$G_n(q)$ denotes the Hecke Clifford algebra, whereas $G_n(q)[0]$ denotes its even part}

\begin{bigTheorem}
    For all but finitely many $q\in \C$ and for $n > 2$, there are isomorphisms
    \[
    CW_n(\sigma, \epsilon) \cong (e_{n - 1}) \oplus G_n(q)
    \]
    \[
    CW_n(\sigma, \epsilon)[0] \cong (e_{n - 1})[0] \oplus G_n(q)[0].
    \]
    The cap-cup ideal $(e_{n - 1})$ is the Jones basic construction for the extension $CW_{n-2}(\sigma, \epsilon) \sbs CW_{n-1}(\sigma, \epsilon)$ and similarly for $CW_n(\sigma, \epsilon)[0]$.
\end{bigTheorem}

Using this decomposition, we can exploit the known representation theory for the Hecke Clifford algebras to inductively construct the Bratelli diagram for both the Clifford BMW superalgebras and their even subalgebras. 
This in turn gives an indexing for the simple objects in the quantum subgroups of types B C and D coming from conformal embeddings which stabilizes as the ranks in each series approach infinity, as well as in their fermion condensations by the fermion $g$. More than that, it also allows us to inductively construct stabilizing branching rules. We state the limiting rule for type C adjoint conformal embeddings.

\begin{bigTheorem}
    The simple objects in $\fanC(\fsp_{2r}, r + 1)_A$ appearing in all tensor powers $X^{\tensor n}$ for all $n < r/2- 1$ are indexed by strict partitions $\mu$ with $|\mu|\le n$, with a sign choice for $\ell(\mu)$ even. The branching rules are given by 
    \[
    V_{\mu} \tensor V_{(1)} \cong \bigoplus_{2\not|\ell(\lambda) , \mu \leftrightarrow \lambda} 2 V_\lambda \oplus \bigoplus_{2|\ell(\lambda), \mu \leftrightarrow \lambda } V_{\lambda^+} \oplus V_{\lambda^-}
    \]
    for $\ell(\mu)$ odd,
    \[
    V_{\mu^{\pm}} \tensor V_{(1)} \cong \bigoplus_{2 | \ell(\lambda)} V_{\lambda^+} \oplus V_{\lambda^-} \bigoplus_{2 \not | \ell(\lambda)} V_\lambda
    \]
     for $\ell(\mu)$ even, where both sums are over all $\lambda$ with $\ell(\lambda)\le r$ such that $\lambda$ is obtained from $\mu$ by adding or removing a box.

     The simple objects in the fermion condensation $\fanC(\fsp_{2r}, r + 1)_A/\gamma$ appearing in all tensor powers $X^{\tensor n}$ for $n < r/2 - 1$ are indexed by strict partitions $\mu$ with $|\mu| \le n$, with $V_\mu$ type $Q$ if $\ell(\mu)$ is odd, and $V_\mu$ type $M$ if $\ell(\mu)$ is even, and with fusion rules given by
     \[
    Hom_{CW}(V_\mu \tensor V_{(1)} \to V_\lambda) \cong \begin{cases}
        \C^{2|2} & \mu \leftrightarrow \lambda, \ell(\mu)\ell(\lambda) \mod 2 = 1\\
        \C^{1 | 1} & \mu \leftrightarrow \lambda, \ell(\mu)\ell(\lambda) \mod 2 = 0\\
        0 & \mu \not\leftrightarrow \lambda
    \end{cases}
     \]
\end{bigTheorem}

As a result, the limiting fusion graph for the fermion condensations $(\fanC(\fsp_{2r}, r + 1)_A)/\gamma$ up to level a distance of $6$ from the unit object is as follows for $r$ sufficiently large:
\[
\att{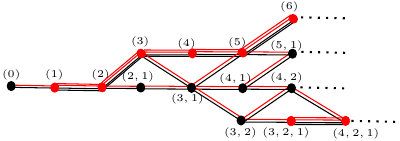}{1.5}
\]
Red edges represent odd fusion channels while black edges represent even fusion channels. Red dots represent type Q simples, while black dots represent type M simples. 

Similarly, the limiting fusion graph for $\fanC(\fsp_{2r}, r + 1)_A$ is as $ r\to \infty$ is
\[
\att{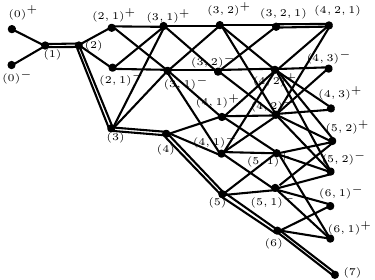}{1.5}.
\]

Note in particular that by specializing to $q = e^{\frac{2\pi i}{24}}$ one recovers the fusion graph for $\mathcal{C}(\fsp_{4}, 3)_A \cong \mathcal{C}(\so_{3}, 2)_A$,
\[
\att{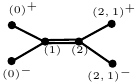}{2},
\]
which was already known from \cite{KnownCasesGB}. 

At all but countably many $q\in \C$, these are literally the fusion graphs of their respective categories. Further knowledge of the representation theory of the Clifford BMW superalgebras will be required to fully understand specializations at roots of unity. We have been able to work out some of the smaller specializations explicitly, as well as closed formulas for many of the quantum dimensions, but we have so far been unable to handle every case.

The story of the Clifford BMW superalgebras and the Hecke Clifford algebras gives a very strong analogy to the story of the BMW algebras and the Hecke algebras, as well as the Brauer algebras and the symmetric group algebras. We suspect the representation theory for the Clifford BMW superalgebras to be just as rich as each of the above examples, only connecting also to the combinatorics of the Sergeev algebras and the twisted symmetric group algebras!

\section*{Acknowledgments} I would like to thank my advisor, Edie-Michell, for introducing me to the incredible world of skein theory and diagram categories. It was also him who suggested I generalizing his, Wenzl, and Snyder's work in type $A$ to the remaining infinite families of conformal embeddings. The LaTeX command for inserting PDF diagrams into mathmode was given to me by Edie-Michell. Without this command, the paper would be far uglier. For the commutative diagrams, I used the online tool Quiver.\footnote{\url{https://q.uiver.app/}} I would also like to think the University of New Hampshire Department of Mathematics and Statistics for providing funding for my research.

\section*{AI Disclosure} In producing this manuscript, I have occasionally used ChatGPT to search the literature and do keyword searches. While ChatGPT was used to initially explore some of the mathematical literature, I have read the papers cited myself (with varying degrees of care). All of the editorial decisions and writing are my own work. No AI generated text is present in this paper. The mathematics in this paper is either from cited sources or my own work. 
\section{Preliminaries}

\subsection{Categorical Background}

The following proposition undergirds all of the diagrammatic techniques used in this paper, as well as those used in \cite{Caleb}, \cite{CainAndSnyder} and \cite{ZhengweiLiu}. A proof of this proposition was provided in \cite{CainAndSnyder}

\begin{proposition}[\cite{CainAndSnyder}]\label{semisimplification}
    Let $C$ be a spherical category dagger category. Let $D$ be a unitary tensor category. Let $F: C \to D$ be a pivotal functor which is also a dagger functor. Then there exists a unique faithful functor $\sms{F}: \sms{C} \to D$ which descends through the semisimplification functor. 
\[\begin{tikzcd}
	C & D \\
	{\sms{C}}
	\arrow["F", from=1-1, to=1-2]
	\arrow["S"', from=1-1, to=2-1]
	\arrow["{\sms{F}}"', from=2-1, to=1-2]
\end{tikzcd}\]
In this case, the category $\sms{C}$ is unitary.
\end{proposition}

Let $C$ be a pivotal, semisimple braided $\C$-linear tensor category with unit object $\mathbb{I}$. 

\begin{definition}
    An algebra object $A$ in $C$ is called an Etale algebra object if $A$ is commutative, and if the multiplication map $\mu: A\tensor A \to A$ admits a section $s: A \to A\tensor A$ which is also an $A$-bimodule homomorphism. 
\end{definition}

\begin{definition}
    Let $(A, \mu, u)$ be an algebra object in $C$, such that $Hom_\mathcal{C}(A \to \mathbb{I}) \cong \C$. Then $A$ is \textbf{connected}.
\end{definition}

Recall that a pivotal structure, when paired with a braiding on $C$, automatically turns $C$ into a ribbon category. Let $\theta_X : X \to X$ denote the twist on any object $X$ in $C$. 

Let $C_A$ denote the category of left $A$-modules in $C$. As described in \cite{CainAndSnyder} and \cite{qMacky}, and elsewhere, the free module functor $F_A: C \to C_A$ admits an adjoint $U: C_A \to C$ which maps $(M, \rho) \mapsto M$. For convenience, we recall the following well-known result:

\begin{proposition}\label{freeModuleFunctorPivotal}
    If $A$ is a connected Etale algebra, then $C_A$ has the structure of a semisimple, rigid monoidal category, such that the functor $F_A : C \to C_A$ is monoidal with a suitable tensorator. 
    If in addition $\theta_A = I_A$, then $C_A$ admits a pivotal structure such that $F_A: C \to C_A$ is pivotal. 
\end{proposition}

This was proven in \cite{qMacky} and elsewhere. The main application of this fact, for us, will be that it allows us to use Proposition \ref{semisimplification}. Edie-Michell and Snyder used this fact in \cite{CainAndSnyder} in order to obtain their results for type A, which we generalize to types B, C, and D. For types B and D, most of this can be thought of as book-keeping. Nonetheless, for type C, we will have to adjust the pivotal structure to fit our needs, which requires us to pay careful attention to what happens to $\theta_A$.

The following result describes a Galois-like correspondence between ``intermediate functors'' factoring the free module functor $F_A: C \to C_A$. This correspondence is well known, and was in fact used explicitly in \cite{CainAndSnyder}. Because we will be working with multiple examples, it will be convenient to lay it out explicitly for later reference.

\begin{proposition}[\cite{Exact}]\label{Galois}
    Let $C$ be a braided fusion category, and let $A$ be a connected Etale algebra. Then $F_A: C \to C_A$ is dominant. Moreover, if $F_1: C \to D$ and $ F_2: D \to C_A$ are dominant functors with $D$ finitely semisimple and with simple unit such that $F_A \cong F_2 \circ F_1$, then there exists a subalgebra $A'$ of $ A$ such that the diagram commutes up to natural isomorphisms
\[\begin{tikzcd}
	C & D & {C_A} \\
	& {C_{A'}}
	\arrow["{F_1}", from=1-1, to=1-2]
	\arrow["{F_A}"{description}, curve={height=-18pt}, from=1-1, to=1-3]
	\arrow["{F_{A'}}"{description}, from=1-1, to=2-2]
	\arrow["{F_2}", from=1-2, to=1-3]
	\arrow["K", from=1-2, to=2-2]
    \arrow["{F'}", from=2-2, to=1-3]
\end{tikzcd}\]    
\end{proposition}
\begin{proof}
    Since $C$ and $D$ are semisimple and $F_1$ is dominant, $F_1$ admits a right adjoint $F_1^{\vee}: D \to C$, which is also a left adjoint by semisimplicity of $C$ and $D$. Since two-sided adjoints preserve limits and colimits, $F_1$ is exact. Also, $F_1^{\vee}$ is lax monoidal. Similarly, $F_2$ admits a lax monoidal, two-sided adjoint $F_2^{\vee}$, so that $F_2 \circ F_1 \dashv F_1^{\vee} \circ F_2^{\vee}$. Since $F_1^{\vee}$ is lax monoidal, $F_1^{\vee}(\mathbb{I})$ is endowed with an algebra structure, given by the tensorator and the unit map. Letting $F_A^{\vee}: C_A \to C$ be the forgetful functor, by the uniqueness of adjoints we have $F_A^{\vee} \cong F_1^{\vee}\circ F_2^{\vee}$. 
    
    By~\cite{Exact} Proposition 5.1, there exists a half braiding $\sigma$ on $A' = F_1^{\vee}(\mathbb{I})$ and a monoidal equivalence $K : D \to C_{A'}$ such that the above diagram commutes up to natural isomorphisms ($F'$ is defined to make this commute). All that remains is to see that $A' \sbs A$ (that is, that there exists a monomorphism $A' \to A$ which is an algebra homomorphism). Letting $F'^{\vee}\cong \inv{K} \circ F_2^\vee$ be the adjoint to $F$, note that $F'^{\vee}$ is equipped with a (necessarily non-zero) map $\epsilon: \mathbb{I} \to F'^{\vee}(A)$\footnote{note that $A$ is the unit in $C_A$}. Since $D$ has simple unit, $\epsilon$ is a mono, and also algebra homomorphism. Since $F_{A'}^{\vee}$ is a double-adjoint, it preserves kernels (which are the same as monos in our context). By uniqueness of adjoints, $F_{A'}^{\vee}\circ F'^{\vee}\cong F_A$, and hence $ F_{A'}^{\vee}\circ F'^{\vee}(A) \cong F_A^\vee(A) = A$, and we have a mono\footnote{noting that the unit in $C_{A'}$ is $A'$}
    \[
    F_{A'}^\vee(\epsilon): A' \to  F_{A'}^{\vee}\circ F'^{\vee}(A) \cong A
    \]
    which is also an algebra homomorphism. 
\end{proof}

We will also need the following fact about the free module functor, which was proven in \cite{CainAndSnyder}.
\begin{proposition}[\cite{CainAndSnyder}]
    For any Etale algebra $A\in C$ with $C$ a braided tensor category, the free module functor $F_A: C \to C_A$ is central. 
\end{proposition}


\subsection{Fermion Condensation and Supercategories}

A \textbf{supercategory} $\mathcal{C}$ is a linear category over some field (we will require the characteristic not to be $2$), enriched over $\underline{sVec}$, the category of super spaces and even maps. Morphisms $f: x \to y$ in $\mathcal{C}$ decompose into an odd and even part, thus turning $Hom_\mathcal{C}( x\to y)$ into a superspace. Composition is defined to be even.

By a monoidal supercategory we mean a super category $C$ equipped with a super monoidal product satisfying the super interchange rule, as defined in \cite{MonoidalSupercats}. In diagrammatic terms, this means that one picks up a sign whenever two odd morphisms slide past each other. So, for $f$ and $h$ homogeneous morphisms, 
\[
\att{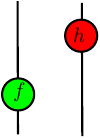}{1} = (-1)^{|f||h|}\att{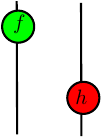}{1}
\]
This can be made more formal by asking that the monoidal product induce an even functor from a twisted Deligne product $C\boxtimes_S C$ of $C$ with itself.

By a pivotal super category, we mean a monoidal super category $C$ equipped with a supernatural transformation $\phi:(\blank)^{**} \Rightarrow I$ from the double dual functor to the identity functor. In terms of diagrammatic calculus applied to super pivotal categories, this means that the 360 degree rotation of an odd morphinism $f$ is $-f$. So, for $f$ homogeneous:
\[
\att{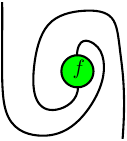}{1} = (-1)^{|f|}\att{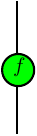}{1}.
\]
A pivotal supercategory $C$ is strictly pivotal if $f^{**} = -f$ and the pivotal structure is the supernatural transformation is the twisted identity functor which maps $f \mapsto (-1)^{|f|}f$ on morphisms. Super pivotal categories appear to have first been defined in \cite{SuperPivotalCats}.

We refer the reader to \cite{SuperPivotalCats} and \cite{MonoidalSupercats} for a more thorough treatment of supercategories of various kinds.

\subsubsection{Reconstructing supercategories}

A \textbf{superalgebra} $A$ is an algebra with an even and odd part $A = A_0 \oplus A_1$ such that $A_i A_j\sbs A_{i + j}$. If $A = A_0 \oplus A_1$ is a superalgebra, a \textbf{supermodule} is a superspace $V = V_0 \oplus V_1$ along with an action of $A$ on $M$ such that $A_i.M_j \sbs M_{i + j}$. If $A$ is semisimple, all of its modules are semisimple, including the graded ones. But it need not be the case that this decomposition respects the grading. This gives rise to two notions of irreducible supermodules: type M and type Q. A module $M$ is type $M$ if $End_{A}(V) \cong \C$, and type $Q$ if $End_A(V) \cong Cl_1$. 

Similarly, we have a notion of type $Q$ and type $M$ simples in any supercategory. An object $x$ in a supercategory $\mathcal{C}$ is a type $M$ simple if $End(x) \cong \C$, and a type $Q$ simple if $End(x) \cong Cl_1$.

We now introduce a slightly generalized notion of the Karoubi completion of a monoidal supercategory, which is probably already well known. Note that one definition was given in \cite{MonoidalSupercats}, but our definition seems to be slightly different, and we make no use of the $\Pi$-completion, treating this as a separate process. We will not place much emphasis on the $\Pi$-completion due to the fact that all fermion condensations are $\Pi$-complete.
\begin{definition}
    Let $\mathcal{C}$ be a rigid monoidal (super) category. Let $S\sbs \mathcal{C}$ be a set of idempotents in $\mathcal{C}$ which satisfy the following properties.

    \begin{enumerate}
        \item Every idempotent $p \in S$ is even. 
        \item For every two idempotents $p, q\in \mathcal{S}$, $p \tensor q \in S$. 
        \item For every $p \in S$, $p^*\in S$
    \end{enumerate}

    Then $Kar_S(\mathcal{C})$ is the partial Karoubi completion of $\mathcal{C}$ by $S$, which is the full subcategory of $Kar(\mathcal{C})$ spanned by $S$. Explicitly, for idempotents $p: x \to x, q: y \to y\in S$, 
    \[
    Hom_{Kar_S(\mathcal{C})}(p \to q) = q \circ Hom_{\mathcal{C}}(x \to y) \circ p.
    \]
    Since $p$ and $q$ are even, $q \circ Hom_{\mathcal{C}}(x \to y) \circ p$ is guaranteed to be a super-subspace of $ Hom_{\mathcal{C}}(x \to y)$.

    If $\mathcal{C}$ is a monoidal supercategory, then $SKar(\mathcal{C})$ shall be understood to be the partial Karoubi completion of $\mathcal{C}$ by the set of even idempotents.
\end{definition}

In general, $Kar(\mathcal{C})$ need not be a monoidal supercategory, or even a supercategory at all! It is entirely possible for an idempotent to be inhomogeneous. This is why we require that every idempotent in $S$ be even. 

An instructive example of a problematic idempotent would be the idempotents $\frac{1}{2}(1 + v)$ and $\frac{1}{2}(1 - v)$ in the algebra $Cl_1 = \gen{v|v^2 = 1}$, which is graded by $|1| = 0, |v| = 1$. One easily check that these are idempotents, and that they are minimal. But they are not homogeneous.

\begin{proposition}
    For every set $S$ of idempotents in a rigid monoidal supercategory $\mathcal{C}$, $Kar_S(\mathcal{C})$ is a monoidal supercategory. If $\mathcal{C}$ is a super-pivotal category, so is $Kar_S(\mathcal{C})$. 
\end{proposition}

We leave the proof of this to the reader.


 If a category $\mathcal{C}$ is semisimple and karoubi complete, then simple objects appearing as subobjects of $x\in \mathcal{C}$ are in correspondence with minimal idempotents in $\mathcal{C}$, which in turn are in correspondence with the category of algebra representations $Rep(End_{\mathcal{C}}(x))$. This carries over perfectly (super well in fact) into the super world, with a slight adjustment. 

\begin{proposition}
    Let $\mathcal{C}$ be any supercategory. Let $p\in End_\mathcal{C}(x)$ be an even idempotent with $x$ an object in $\mathcal{C}$. The following are equivalent. 
    \begin{enumerate}
        \item $p$ is a type $Q$ simple in $SKar(\mathcal{C})$
        \item $p = a + b$ for inhomogeneous minimal idempotents.
        \item $End_\mathcal{C}(x)p$ is a type $Q$ irreducible representation of $End_\mathcal{C}(x)$.
        \item $p End_{\mathcal{C}}(x) p \cong Cl_1$ as superalgebras. 
    \end{enumerate}
\end{proposition}

The proof of this proposition is routine. 

Recall that a \textbf{dagger structure} on a $\C$-linear category $\mathcal{D}$ is just a sesquilinear functor $\dag: \fanD^{op} \to \fanD$. A monoidal dagger structure on a monoidal $\C$-linear category $\fanD$ is simply a dagger structure which is also a monoidal functor. There are also notions of braid dagger structures, and so on, each corresponding compatibility with various categorical structures. 

We propose the following definition of a super-dagger structure:
\begin{definition}
    A \textbf{super-dagger} structure $\dag$ on a super-category $\fanD$ is a dagger structure on $\fanD$ which is also an even functor. A \textbf{monoidal super-dagger} structure is simply a super-dagger structure which happens to also be a super-monoidal functor.
\end{definition}

If a monoidal supercategory $\fanD$ is pivotal, one can define a categorical trace. In the case where $\fanD$ is a super-pivotal category, one can also define a categorical trace.
\begin{definition}
    Suppose that $\mathcal{D}$ has simple unit. Define a Hermetian form on $Hom_{\fanD}(x \to y)$ for all objects $x, y\in \fanD$ by
    \[
    \inner{f}{g} = tr(f \circ g^\dag)
    \]
    A super-dagger structure on a supercategory with simple unit, $\fanD$, is \textbf{unitary} if $\inner{f}{f} = |f|^2 > 0$ unless $f = 0$. 
\end{definition}

\begin{definition}
    A superpivotal structure is \textbf{spherical} if the left trace is the right trace. Let $\mathcal{D}$ be a monoidal supercategory. A \textbf{superideal} $\mathcal{J}$ of $\mathcal{D}$ consists of the following:
    \begin{enumerate}
        \item For all $x, y$ objects in $\fanD$, one has a linear subspace
        \[
        \mathcal{J}(x\to y) \sbs Hom_{\mathcal{D}}(x \to y)
        \]
        such that if $f\in \mathcal{J}(x\to y)$ then $f_0, f_1 \in \mathcal{J}(x\to y)$.
        \item Whenever $f\in \mathcal{J}( x\to y)$, then for all $g\in Hom_{\fanD}( x' \to y')$ we have $f \tensor g \in \mathcal{J} ( x \tensor x' \to y \tensor y')$ and $g \tensor f \in \mathcal{J} ( x \tensor x' \to y \tensor y')$ 
        \item Whenever $f\in \mathcal{J} ( x \to y)$ and $g \in Hom_\fanD( y \to z)$ and $h \in Hom_{\fanD} (w \to z$ we have $h \circ g \circ f $.
    \end{enumerate}
\end{definition}

Note that, just like the usual notion of a monoidal ideal in a monoidal category, superideals are \textit{not} subcategories: they lack identity morphisms. The fate of a superideal is to be quotiented out.

\begin{definition}
    If $\mathcal{J} \sbs \fanD$ is a superideal, define the monoidal supercategory 
    \[
    \fanD/\mathcal{J}
    \]
    with hom-spaces 
    \[
    Hom_{\fanD/\mathcal{J}} ( x \to y) = Hom_{\fanD}( x\to y) /\mathcal{J} ( x \to y)
    \]
\end{definition}

\begin{proposition}
    Whenever $\mathcal{J}$ is a superideal, $\mathcal{D}$ remains a monoidal supercategory. If $\mathcal{J}$ is pivotal, so is $\fanD/\mathcal{J}$.
\end{proposition}
\begin{proof}
    Clear.
\end{proof}

\begin{definition}
    Let $\fanD$ be a pivotal supercategory. For all objects $x, y\in \fanD$ define 
    \[
    \mathcal{N}(x \to y) = \defset{f: x \to y}{ \forall g: y \to x, tr(f \circ g) = 0}.
    \]
\end{definition}

\begin{proposition}
    For every strict pivotal supercategory $\fanD$, $\mathcal{N}\sbs \mathcal{D}$ is a superideal so long as $\fanD$ is spherical.
\end{proposition}
\begin{proof}
    The proof that $\mathcal{N}$ is closed under tensoring by other morphisms and composition is identical to the proof that the usual negligible ideal is an ideal. 

    All that remains to show is that the even and odd part of every morphism in $\mathcal{N}$ is also in $\mathcal{N}$. The key is to note that the categorical trace is even for any super pivotal category (i.e. $tr(f)= tr(f_0)$).

    Suppose that $f\in \mathcal{N} (x \to y)$. We must show that $f_0$ and $f_1 \in \mathcal{N} (x \to y)$. Suppose that $h \in Hom_{\fanD}( y \to x)$. Then $tr(f_0 h) = tr(fh) = 0$ since $(fh)_0 = f_0 h$. On the other hand, if $h \in h \in Hom_{\fanD}( y \to x)_1$ then $|f_0h| = 1$ and hence $tr(f_0h) = 0$. As a result, for any morphism $h: y \to x$ we have $tr(f_0 h) = 0$. Hence $f_0 \in \mathcal{N}(x \to y)$. Hence $f_1 = f - f_0 \in \mathcal{N}(x \to y)$.
\end{proof}

\begin{definition}
    The semisimplification $\sms{\fanD}$ of a spherical supercategory $\mathcal{D}$ is the spherical supercategory 
    \[
    \sms{\fanD} = \fanD/\mathcal{N}.
    \]
    The even pivotal functor $P:\fanD \to \sms{\fanD}$ which maps $f: x \to y$ to $f + \mathcal{N} (x \to y)$ is called the \textbf{semisimplification} functor.
\end{definition}

All of this was defined in order to prove the following super-analog of the descend lemma. 
\begin{proposition}[Super Descent Lemma]\label{superDescent}
    Let $\fanC$ be a strictly pivotal, spherical super category with simple unit. Let $\fanD$ be a super-pivotal super-unitary category. 
    
    Let $F: \fanC \to \fanD$ be an even pivotal functor. Then there exists a unique even faithful superpivotal functor $\sms{F} : \sms{\fanC} \to \fanD$ making the diagram commute:
\begin{tikzcd}
	{\mathcal{C}} & {\mathcal{D}} \\
	{\overline{\mathcal{C}}}
	\arrow["F", from=1-1, to=1-2]
	\arrow["P"', from=1-1, to=2-1]
	\arrow["{\overline{F}}"', dotted, from=2-1, to=1-2]
\end{tikzcd}

This descent respects even parts, as does the semisimplification. One has
\[
\sms{C[0]} = \sms{C}[0]
\]
and 
\[
\sms{F}[0] = \sms{F[0]}
\]
\end{proposition}

The proof is essentially identical to that of the usual Descent Lemma.




One can also define the additive completion of a supercategory, $Add(\mathcal{C})$, with no adjustment at all. The resulting category is super monoidal if $\mathcal{C}$ was, and also superpivotal if $\mathcal{C}$ was. 

\begin{definition}
    The Abelianization of a supercategory $\mathcal{C}$ is the category $Ab(\mathcal{C})= Add(SKar(\mathcal{C}))$. 
\end{definition}

If $\mathcal{C}$ is superpivotal or supermonoidal, so is $\mathcal{C}$. 

One also has the identity, 
\[
Ab(\mathcal{C}) = Ab(\mathcal{C}[0]).
\]
This is an equality, not just an equivalence. Objects in $Ab(\mathcal{C})$ are formal direct sums of even idempotents and morphisms are matrices morphisms in $\mathcal{C}$ which absorb those formal direct sums. Objects in $Ab(\mathcal{C}[0])$ are formal direct sums of even idempotents as well, while morphisms are even morphisms in $\mathcal{C}$ which absorb those formal direct sums of idempotents.

\subsubsection{Fermion Condensation}

We will now discuss a common way of constructing superpivotal categories: fermion condensation. Fermion condensation is far from new. Our formal definition of fermion condensation appears to be essentially the same as Lacabanne's in \cite{AbelSuperModular}, only in his setting he worked with a full braid. Our conventions differ slightly from Lacabanne's. The results of this section are not new, but the author is not sure how far back the formal definitions go, or who came up with them first.

Throughout the rest of this section, let $P$ be a strict monoidal, strictly pivotal category. Let $g\in P$ be any object such that $g\tensor g \cong \mathbb{I}$. We say that $g$ is a \textbf{fermion} with half braid $\gamma$ if there exists a half braid $\gamma$ on $g$, with the property that $\gamma_g = -1_{g\tensor g}$. Throughout the rest of this subsection, let $g$ be fixed, with a fixed half braid $\gamma$ on $g$ turning $g$ into a fermion. 

\begin{definition}
    The fermion condensation of $P$ by $\gamma$, denoted $P/\gamma$, is the category with morphisms and objects defined as follows. 
    \begin{itemize}
        \item $Ob(P/\gamma) = Ob(P)$.
        \item Given $x, y\in Ob(P/\gamma)$, $Hom_{P/\gamma}(x \to y) = Hom_P(x \to y) \oplus Hom_{P}(x \to y \tensor g)$. 
\end{itemize}
A morphism $f: x\to y$ is even if it is a morphism $f: x\to y$ in $P$, and odd of it is a morphism $f: x \to y \tensor g$. We write $Hom_{P/\gamma}(x \to y)_0 =Hom_P(x \to y)\oplus 0$ and $Hom_{P/\gamma}(x \to y)_1 = 0 \oplus Hom_{P^g}(x\to y \tensor g)$. This turns $Hom_{P/\gamma}(x \to y)$ into a superspace. 

Composition of morphisms is defined piecewise depending on the degree of composing homogeneous morphisms, and extended bilinearly. To distinguish composition in $P/\gamma$ from composition in $P$, we let $\circ$ denote composition in $P$, while $\circ_g$ shall denote composition in $P/\gamma$. 

\begin{itemize}
    \item If $h\in Hom_{P/\gamma} (x \to y), k\in Hom_{P/\gamma}(y \to z)$ are both even, then $k \circ_g h = k \circ h$.
    \item If $h$ is odd and $k$ is even $h\in Hom(x\to y \tensor g)$ and $k \in Hom(y\to z)$ we define
    \[
    k\circ_g h = (k \tensor 1_g)\circ h
    \]
    \item If $h$ is even and $k$ is odd, we define
    \[
    k \circ_g h = k \circ h
    \]
    \item If both $h$ and $k$ are odd, we define
    \[
    k \circ_g h = \epsilon \circ (k \tensor 1_g) \circ h
    \]
    where $\epsilon: g\tensor g \to \mathbb{I}$ is a fixed isomorphism. 
\end{itemize}
\end{definition}

Note that composition of morphisms has been defined to be an even map. We need to check that this really does turn $P/\gamma$ into a category. For this, we need to verify that identity morphisms are present, and that composition is associative.

\begin{proposition}
    $P/\gamma$ is a super category.
\end{proposition}
\begin{proof}
    There are eight cases to check depending on the parity of each of the three morphisms. The only non-trivial case to verify is associativity in the case of three odd morphisms. Using the fact that $\gamma_g = -1_{g\tensor g}$, it is easy to verify that $ \epsilon \tensor 1_g = 1_g \tensor \epsilon$ as morphisms $g \tensor g \tensor g \to g$.
    
    Using this, one calculates
    \[
    (f \circ_g h)\circ_g k = \att{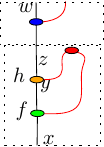}{1} = \att{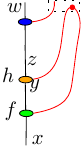}{1} = \att{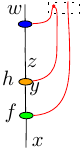}{1} = \att{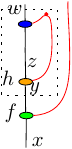}{1} = f\circ_g (h \circ_g k) 
    \]
    We leave the formal manipulations from the axioms of a strict monoidal category to the reader. 

    Identity morphisms are simply the identity morphisms in $P$. Thus far we have shown that $P/\gamma$ is a category. 

    Note that it is also clear that composition is an even map $\circ_g: Hom_{P/\gamma}( x\to y) \tensor Hom_{P/\gamma}( y \to z) \to Hom_{P/\gamma}(x \to z)$. Hence $P/\gamma$ is a super category. 

\end{proof}

\begin{theorem}
    Define tensor products on $P/\gamma$ piecewise as follows. Let $f: x\to y$ and $h: x' \to y'$ in $P/\gamma$. The tensor product is defined as follows depending on the grading of $f$ and $h$. 
    \begin{itemize}
        \item[$0\tensor 0$] $f\tensor_g h = g\tensor h$
        \item[$0 \tensor 1$] $f \tensor_g h = g\tensor h$
        \item[$1 \tensor 0$] $f \tensor_g h = (f \tensor 1_y) \circ (1_{x'} \tensor \sigma_y) \circ (1_{x'} \tensor h \tensor 1_g)$
        \item[$1 \tensor 1$] $f \tensor_g h = (f \tensor 1_y)\circ (1_{x'} \tensor \sigma_y) \circ (1_{x'} \tensor h \tensor 1_g) \circ (1_{x'} \tensor 1_{y'} \tensor \epsilon)$
    \end{itemize}

    This tensor product rule turns $P/\gamma$ into a strict monoidal supercategory. 
    
    Moreover, a super rigidity structure on $P/\gamma$ is inherited from $P$. The resulting category $P/\gamma$ is a pivotal supercategory. 
\end{theorem}

\begin{proof}
    The interchange rule is delightfully clear from a diagrammatic argument. We leave the formal verification to the reader, which simply comes down to unfolding the definitions of composition and tensor product in $P/\gamma$, the assumption that $\gamma_{g, g} = - I_{g\tensor g}$, and the axioms of a strict monoidal category. For two odd morphisms $f$ and $h$ (shown as $g$ in the diagram), we have 
    \[
    \att{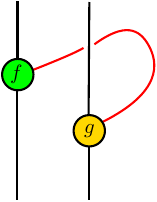}{0.5} = \att{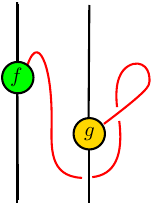}{0.5} = -\att{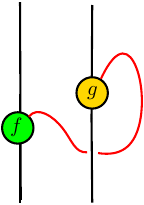}{0.5}.
    \]

    The cases for $f$ odd and $h$ even, $f$ even and $h$ odd, or for $f$ and $h$ both even, are immediate from drawing the diagrams. 

    A rigidity structure on objects in $P/\gamma$ (the objects are the same as in $P$) is inherited from $P$, and the evaluation maps and coevaluation maps are clearly even (since they are morphisms in $P$). 

    A 360 degree rotation of an even map $f: x \to y$ is the same as in $P$. On the other hand, suppose $f: x \to y$ is an odd map in $P/\gamma$. Then $f: x \to y \tensor g \in P^g$. 

    We now unfold the definitions of composition and tensor products. While this can be done purely formally from the axioms and the definition of the super monoidal product, the argument is much clearer diagrammatically. 

    \[
    \left(\att{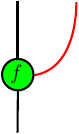}{1}\right)^{**} = \att{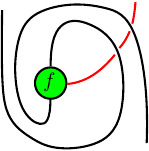}{1} = \att{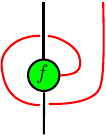}{1} = \att{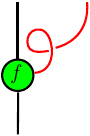}{1} = -\att{oddMorphism}{1}
    \]

\end{proof}

\begin{remark}
    Intuitively, the idea is that the morphism $f$ is a rod attached to a ribbon, the red strand, which is attached to a handle. As the rod rotates, the ribbon accumulates a twist. Note that this argument is actually nearly identical to the argument in \cite{qMacky} that the modules over a connected Etale algebra $A$ form a pivotal category exactly when $\theta_A = 1$. 
\end{remark}

\subsection{Some Young Diagram Combinatorics}\label{YoungDiagramCombinatorics}

\begin{definition}
    A \textbf{Young diagram} is a finite subset $\overline{\lambda} \sbs \N^2$ such that for any $(x, y)\in \overline{\lambda}$, $x - 1 > 0$ then $(x - 1, y) \in \sms{\lambda}$, and if $y - 1 > 0$, then $(x, y - 1)\in \sms{\lambda}$.
\end{definition}

There is a bijection between Young diagrams, and weakly decreasing (and eventually zero) sequences of non-negative integers. The bijection takes $\lambda = (\lambda, \lambda_2, \dots) \mapsto \overline{\lambda}$ where $\overline{\lambda} = \defset{(i,j)}{j \le \lambda_i}$. If there is no worry about confusion, then we will tend to identify Young diagrams with their sequence of row-lengths. The \textbf{length} of a Young diagram $\lambda$ is the number of rows (equivalently, the largest $i$ such that $\lambda_i > 0$). If $\lambda$ is a partition of $n$, we write $\lambda \dashv n$. 

\begin{definition}
    A \textbf{shifted Young diagram} is a finite subset $\overline{\mu} \sbs \N^2$ such that for any $(x, y)\in \sms{\mu}$, if $x -1 \ge y$ then $(x -1, y) \in \overline{\mu}$, and if $y - 1 > 0$, then $(x, y - 1)\in \sms{\mu}$. 
\end{definition}

Let $\mu = (\mu_1, \dots,\mu_k)$ be a strictly decreasing sequence of natural numbers. such that $\mu_1 + \mu_2 + \dots + \mu_k = n$. We say say $\mu$ is a strict partition of $n$, and write $\mu \dashv_s k$. For the \textbf{length} of $\mu$, we write $\ell(\mu) = k$

Similar to Young diagrams and weakly decreasing sequences of natural numbers, there is a bijection between strict partitions of $n$, and shifted Young diagrams with $n$ boxes. The \textbf{shifted Young diagram representing $\mu$} is the set of points 
\[
\overline{\mu} = \defset{(x,y)\in \N^2}{y \le x \le \mu_i + y, y \le n}
\]

As with Young diagrams, we will often identify a strict partition $\mu$ with its corresponding shifted Young diagram $\overline{\mu}$ unless the distinction needs to be emphasized. 

\begin{definition}
    Let $\mu = (\mu_1, \dots, \mu_k)$ be a strict partition. Then we can count \textbf{columns} of $\mu$,
    \[
    \mu^t_i = \#\set{(i,j)\in \overline{\mu}}
    \]
    and thus define the \textbf{transpose} of $\mu$, $\mu^t = (\mu^t_1, \dots, \mu^t_c)$ where $c\in \N\U \set{0}$. Note that $\mu^t$ is \textit{not} usually a strict partition. 
    
    For notational convenience, let $\mu_i = 0$ if $i > \ell(\mu)$, and similarly for $\mu^t$.
    
    Define the \textbf{doubling} of $\mu$, denoted $D(\mu)$, to be the (non-shifted) Young diagram with rows $D(\mu)_1 = \mu_1$ and $D(\mu)_i = \mu_i + \mu^t_{i - 1}$ for $i > 1$. 
\end{definition}

Doublings of shifted Young diagrams have a number of applications. See \cite{randomStrict} for both a brief background on some of the properties of doubling, as well as their applications. We also refer the reader to \cite[Section I.1]{MacDonald}

\begin{example}
    Below is an example of the \textit{doubling} of $\mu = (5, 4, 2, 1)$. 
    \[
    \mu = \ydiagram[*(red)]{0+5, 1 + 4, 2 + 2, 3 + 1}\rightsquigarrow
    D(\mu) = \ydiagram[*(red) ]
{0+5, 1 + 4, 2 + 2, 3 + 1}
*[*(green)]{5, 5, 4, 4, 4, 2}
    \]
\end{example}

Conventions vary on the definition of doubling. Some define the doubling of $\mu$ to be the transpose of what we have defined. We will actually have use for both $D(\mu)$ and $D(\mu)^T$ in many of our examples. Note that since transposing Young diagrams preserves essentially all combinatorial properties, results concerning one version of the doubling immediately transfer over to the transpose. In order to recognize the transpose of the doubling later on, we will need the following lemma. This fact is probably well known, but we state it here for completeness. 

\begin{lemma}
    Given a strict partition representing a shifted Young diagram $\mu = (\mu_1, \dots, \mu_n)$, $D(\mu)^T = (\mu_1 + \mu_1^T , \mu_2 + \mu^T_2, \dots)$ where we let $\mu_i = 0$ for $i \ge n$. 
\end{lemma}

Many of the combinatorial results in this paper can be reduced to counting paths on the Young graph. 

\begin{definition}
    The Young graph $Y$ is the graph with vertices all Young diagrams, and an edge $\lambda \leftrightarrow \mu$ whenever $\mu$ is obtained from $\lambda$ by adding or removing a box. 
\end{definition}

\begin{definition}
    Let $\lambda\sbs \N^2$ be a Young diagram with $|\lambda| = n$. A \textbf{standard Young tableaux} is a bijective labeling $T: \set{1, \dots, n} \to \lambda$ that is decreasing in rows and columns. An \textbf{oscilating tableaux} of length $k$ is a path on the Young graph from $\lambda$ to $(0)$. Let $f^\lambda$ denote the number of standard Young tableaux and $f^\lambda_k$ denote the number of oscilating tableaux of length $k$. 
    
    Let $\mu\sbs \N^2$ be a shifted Young diagram of size $n$. A \textbf{standard shifted Young tableaux} is a bijective labeling $T: \set{1, \dots, n} \to \mu$ that is decreasing in rows and columns. Let $g^\mu$ denote the number of standard shifted Young tableaux on $\mu$. 
\end{definition}

\begin{remark}
    Note that a standard Young tableaux can be thought of as encoding a set of ``instructions'' for constructing a path $\lambda\to (0)$ in the Young graph. Hence a standard tableaux for $\lambda \dashv n$ is essentially the same thing as a oscilating tableaux of length $n$. Note also that some texts define a standard tableaux to be increasing in rows and columns. Both definitions are essentially equivalent, since a change of labeling $i \mapsto n - i$ produces a bijection between standard tableaux of either definition. 
\end{remark}

We will need as a lemma the following theorem, first discovered by Schur in studying projective representations of $S_n$, and then proven constructively in \cite{schurIdentity} through a bijection similar to the Robinson Schensted correspondence. This bijection can also be seen directly from the representation theory of the Sergeev algebra. This is a lovely coincidence, since the Sergeev algebra turns out to also be directly related to the representation theory of the Clifford BMW algebra.

\begin{lemma}[Schur]\label{schurIdentity}
    For any strict partition $\mu$ of $n$, let $\ell(\mu)$ denote the length of the partition. Write $\mu \dashv_s n$ to denote a strict partition. Then \[
    \sum_{\mu \dashv_s n} 2^{n- \ell(\mu)}(g^\mu)^2 = n!\]
\end{lemma}

Obviously, this means that $\sum_{\mu \dashv_s n} 2^{-\ell(\mu)}(g^\mu)^2 = n!/2^n$.

Further, we will need the following lemma, which allows one to write the number standard Young tablauex of $D(\mu)$ for $\mu \dashv_s n$, in terms of the number of shifted standard Young tablauex of $\mu$. This was proven in \cite{randomStrict}, as well as a number of other places. We now recall it here.

\begin{lemma}\label{strictToDoubled}
For any $\mu\dashv_s n$ a strict partition,
    \[f^{D(\mu)} = 2^{-\ell(\mu)}(g^\mu)^2 {2n \choose n}\]
\end{lemma} 

We are now in a position to prove the following. 

\begin{proposition}
    \[
    \sum_{\mu\dashv_s k} f^{D(\mu)} = (2k - 1)!!
    \]
\end{proposition}
\begin{proof}
    By Lemma~\ref{strictToDoubled}, Lemma~\ref{schurIdentity}, and the basic identity $(2n)!/(2^nn!) = (2n - 1)!!$, we have
    \begin{align*}
        \sum_{\mu \dashv_s k} f^{D(\mu)} = \sum_{\mu \dashv_s k}2^{-\ell(\mu)} (g^\mu)^2 {2k \choose k} = 2^{-k}{2k\choose k} k! = (2k)!/(2^{-k}k)! = (2k - 1)!!.
    \end{align*}
\end{proof}

Using an explicit bijection similar to the Robinson Schensted correspondence, the following lemma was proven in \cite{oscilatingTabs}.

\begin{lemma}\label{standardToOScilating}
    For any $n, k \in \N \U \set{0}$, and for any $\lambda \dashv k$, the following identity holds,
    \[
    f_n^\lambda = (n - k - 1)!!{n \choose k} f^\lambda
    \]
\end{lemma}

Next, the following lemma can be proven by a direct manipulation of the definitions 
\begin{lemma}\label{gnarlyIdentity}
    For any $k \le n$, we have 
    \[
    {2n \choose 2k}(2n - 2k - 1)!! (2k - 1)!! = {n \choose k} ( 2n - 1)!!
    \]
\end{lemma}
\begin{proof}
    Using the identity $(2r - 1)!! = (2r)!/(2^rr!)$ and writing out the definitions of the binomial coefficients, we have
    \begin{align*}
    {2n \choose 2k}(2n - 2k - 1)!! (2k - 1)!! \\
    = \frac{(2n)!}{(2n - 2k)! (2k)!}\cdot \frac{(2n - 2k)!}{2^{n - k}(n - k)!} \cdot \frac{(2k)!}{2^kk!} \\
    = \frac{n!}{( n - k)!k!}\frac{(2n)!}{2^nn!} = {n \choose k}(2n - 1)!!
    \end{align*}
\end{proof}

We are now in a position to prove the main result of this section. 

\begin{theorem}
    The number of paths on the Young graph $Y$ of length $2n$, beginning from doubled strict partitions and ending at $\emptyset$, is $2^{n}(2n - 1)!!$. 
\end{theorem}

\begin{proof}
    If $\mu \dashv_s m$ for $m > n$, then $D(\mu) \dashv 2m > 2n$, and there are no paths of length $2n$ from $D(\mu)$ to $\emptyset$. Therefore the number of paths on $Y$ of length $2n$ beginning at doubled strict partitions and ending at $\emptyset$ is
    \begin{align*}
    \sum_{k = 0}^n \sum_{\mu \dashv_s k} f_{2n}^{D(\mu)} 
    = \sum_{k = 0}^n \sum_{\mu \dashv_s k}(2n - 2k-1)!!{2n \choose 2k} f^{D(\mu)}  \\
    \\
    = \sum_{k = 0}^n (2n - 2k - 1)!!{2n \choose 2k}\sum_{\mu \dashv_s k} f^{D(\mu)} \\
    = \sum_{k = 0}^n (2n - 2k - 1)!! (2k - 1)!! {2n \choose 2k} \\
    = \sum_{k = 0}^n {n \choose k}(2n - 1)!! = 2^n (2n - 1)!!
    \end{align*}
\end{proof}

In application, we will need slight variation of these identities. 
\begin{lemma}\label{evenShiftedYoung}
    The number of paths on the Young graph $Y$ of length $2n$, beginning from doubled strict partitions of even (odd) size is $2^{n- 1}(2n - 1)!!$.
\end{lemma}

Finally, we will have the occasion to use the following. 

\begin{corollary}\label{typeCPathCounting}
    The number of paths on $Y$ of length $2n$, beginning from the transpose of a doubled strict partition of even size, is $2^{n - 1}(2n - 1)!!$
\end{corollary}
\begin{proof}
    Transposing Young diagrams is an automorphism of the Young graph, from which the corollary follows immediately. 
\end{proof}

\subsection{Fusion Categories of Types B C and D}

We refer the reader to \cite{LieAlgebrasFiniteAffine} for a detailed introduction to Lie algebras in general, as well as \cite{LieTheoryForFusion} for a brief introduction to quantum groups.

Let $\fanC(\g, k)$ denote the level $k$ representations of the affine Lie algebra $\hat{\g}$ corresponding to a semisimple Lie algebra $\g$. 

Let $r \ge 2$. The category $\mathcal{C}(\mathfrak{so}_{2r + 1}, k)^{vec}$ is the subcategory of $\fanC(\so_{2r + 1}, k)$ generated by simple objects labeled by integral weights. The simple objects in these categories are labeled by Young diagrams $\lambda$ such that $\lambda_1' \le r$ and $\lambda_1 + \lambda_2 \le k$. Every simple object in this category appears in a tensor power of $V_{(1)}$

By the Kazhdan Lutsig equivalence, there is a pivotal braid equivalence
\[
\mathcal{C}(\mathfrak{so}_{2r + 1}, k)^{vec} \simeq \sms{Rep(U_q(\so_{2r + 1}))^{vec}}
\]
for $q = e^{\frac{2\pi i}{4(k + 2r - 1)}}$.


Similarly, $\mathcal{C}(\mathfrak{so}_{2r}, k)^{vec}$ is the subcategory of $\fanC(\so_{2r}, k)$ generated by simple objects labeled by integral weights.

The simple objects of $\mathcal{C}(\mathfrak{so}_{2r}, k)^{vec}$ are labeled by Young diagrams $\lambda$ with a choice of sign if $\lambda_1' = r$, such that $\lambda_1' \le r$ and $\lambda_1 + \lambda_2 \le k$.

By the Kazhdan Lutsig equivalence, there is a pivotal braid equivalence
\[
\mathcal{C}(\mathfrak{so}_{2r}, k)^{vec} \simeq \sms{Rep(U_q(\so_{2r}))^{vec}}
\]
for $q = e^{\frac{2\pi i}{4(k + 2r - 2)}}$. 

For types B and D, we will need the following fusion rule, which holds so long as $r \ge 2$:
\[
V_{(1)} \tensor V_{(1)} \cong V_{(2)} \oplus V_{(1, 1)} \oplus V_{(0)}.
\]

We will also need the fusion rule which eventually holds so long as $k, r$ are big enough relative to $|\lambda|$:
\begin{itemize}
    \item For every Young diagram $\lambda$, there exists some $L$ such that whenever $L < k, r$ we have 
    \[
    V_\lambda \tensor V_{(1)} = \bigoplus_{\mu \leftrightarrow \lambda} V_\mu
    \]
    where the sum is over all Young diagrams $\mu$ which are obtained from $\lambda$ by adding or removing a box, and where $V_\lambda, V_\mu, V_{(1)}$ are the simples labeled by their respective Young diagrams in either $\mathcal{C}(\so_{2r}, k)$ or $\mathcal{C}(\so_{2r + 1}, k)$.
\end{itemize}

Finally, $\mathcal{C}(\mathfrak{sp}_{2r}, k)$ has simple objects labeled by Young diagrams $\lambda$ such that $\lambda_1' \le r$ and $\lambda_1\le k$. 

For type C, we will need the following fusion rule:
\begin{itemize}
    \item For every Young diagram $\lambda$ with $\lambda_1 \le k, \lambda_1 \le r$, we have
    \[
    V_{(1)} \tensor V_\lambda \cong \bigoplus_{\mu} V_{\mu}
    \]
    where the sum is over all Young diagrams $\mu$ which are obtained from $\lambda$ by adding or removing a box, and where $\mu_1 \le k$ and $\mu_1' \le r$. In other words, the fusion graph relative to $V_{(1)}$ for the category $\mathcal{C}(\mathfrak{sp}_{2r}, k)$ is the truncation of the unoriented Young graph, where all of Young diagrams which are too big are simply removed.
\end{itemize}

More generally, all of the quantum group categories of types B C and D treated in this paper share a common stable Grothendiek ring, commonly denoted $Rep(O(\infty))$, although we won't need this fact. Remarkably, a this paper can be interpreted as showing that a similar phenomenon holds for the quantum subgroups of types B C and D treated in this section. Most likely, these are all of the $1$-supertransitive subgroups, although future work will be needed to verify this.



\subsection{Skein Theory for Fusion Categories of Types B, C and D}\label{BMWIntro}

\begin{definition}
    For $q, \theta\in \C$, define the $\C$-linear strictly monoidal, strictly pivotal category with symmetrically self dual objects $\N_0$, and with generating morphisms $\att{braid}{0.5}$, $\att{inverseBraid}{0.5}$ subject to the relations
    \[
    \text{(R2) }\att{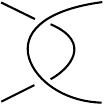}{0.5} = \att{doubleStick}{0.5} \text{ (R3) } \att{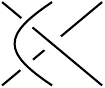}{0.5} = \att{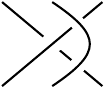}{0.5}
    \]
    \[
    \text{ (Dubrovnik) }\att{braid}{0.5} - \att{inverseBraid}{0.5} = (q - \inv{q})\left(\att{doubleStick}{0.5} - \att{capCup}{0.5}\right)
    \]
    \[
    \text{ (R1) }\att{twist}{0.5} = \theta \att{stick}{0.5}
    \]

    For a fixed $n\in \N_0$, the algebra denoted $BMW_n(q, \theta)$ is the endomorphism algebra $End_{BMW(q, \theta)} (n)$. 
\end{definition}

It is an easy consequence of the relations that $\att{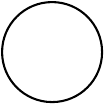}{0.5} = \delta(q, \theta) = \frac{\theta - \inv{\theta}}{q - \inv{q}} + 1$.

The representation theory and general structure of the BMW algebras has been very well-studied.

\begin{lemma}[Birma, Wenzl]\label{BMW Presentation}
    The BMW algebra $BMW_n(q , \theta)$ has presentation 
    \[
    BMW_n(q , \theta) = \gen{t_i, \inv{t_i}, e_i: i \in \set{1, \dots, n - 1}}/R
    \]
    where $R$ is the ideal generated by the relations
    \begin{enumerate}
        \item[W1] $t_i - \inv{t_i} = (q - \inv{q})(1 - e_i)$
        \item[W2] $t_i e_i = e_i t_i = \inv{\theta} e_i$.
        \item[W3] $e_i e_{i + 1} e_i = e_{i + 1}$ and $e_{i + 1}e_i e_{i + 1} = e_i$.
        \item[W4] $e_i^2 = \delta e_i$
        \item[W5] $t_i\inv{t_i} = 1$
        \item[W6] $t_i t_{i + 1} t_i = t_{i + 1} t_i t_{i + 1}$.
        \item[W7] $e_i e_j = e_j e_i$ when $|i - j| > 1$
        \item[W8] $e_i t_j = t_j e_i$ when $|i - j| > 1$
        \item[W9] $t_i t_j = t_j t_i$ when $|i - j| > 1$.
        \item[W10] $ e_{i + 1} t_i e_{i + 1} =  \theta e_{i + 1}$
        \item[W11] $ e_i t_{i + 1} e_i = \theta e_i $
    \end{enumerate}
\end{lemma}

\begin{remark}
    There are a few well known extensions of the BMW algebra. The affine BMW algebras play an important role in the representation theory of the BMW algebras, giving the Jucys Murphey elements. Certain finite dimensional quotients of these, called the cyclotomic BMW algebras, are also well studied. There are also the tied BMW algebras. Our family, the Clifford BMW superalgebras, appears to be an entirely new extension.
\end{remark}

The BMW relations are evaluable with $\dim BMW_n(q, \theta) = (2n - 1)!!$ for all $q, \theta \in \C$ with $q, q - \inv{q} \ne 0$ and $\theta \ne 0$. For generic $\theta,q \in \C$, $BMW_n(q, \theta)$ is semisimple \cite{Wenzl} with non-degenerate trace. It should be noted, however, that specializing to specific values of $q, \theta$ does not always retain semisimplicity or non-degeneracy of the trace.

The minimal idempotents in $BMW_n(q, \theta)$ are fully understood for generic $q$ and $\theta$. We will not need the full description here, as we only care about a few specific values for $q, \theta$, and $n$. 

\begin{proposition}[\cite{IdempotentsInBMW} Example 3.3,\cite{Wenzl}]
    For any $q, \theta\in \C$ such that $q + \inv{q}, q, \theta, 1 - q\inv{\theta}, 1 + \inv{q}\inv{\theta}\ne 0$, the following elements in $BMW_2(q, \theta)$ form a complete system of orthogonal idempotents:
    \begin{align*}
        &P_{(2)}= \frac{1}{q + \inv{q}}\left(\att{braid}{0.5} + \inv{q} \att{doubleStick}{0.5} + \frac{q - \inv{q}}{1 - q\theta} \att{capCup}{0.5}\right),\\
        &P_{(1,1)} =  -\frac{1}{q + \inv{q}}\left(\att{braid}{0.5} - q\att{doubleStick}{0.5} + \frac{q - \inv{q}}{1 +\theta/q} \att{capCup}{0.5}\right),\\
        &P_{(0)} = \frac{1}{\delta} \att{capCup}{0.5}
    \end{align*}
\end{proposition}

\begin{lemma}\label{eigenvaluesOfTheBraid}
    For any $q, \theta\in \C$ such that $q + \inv{q}, q, \theta, 1 - q\inv{\theta}, 1 + \inv{q}\inv{\theta}\ne 0$, $P_{(2)}$, $P_{(1,1)}$ and $P_{(0)}$ are eigenvectors for multiplication by $\att{braid}{0.5}$ with eigenvalues $q, -\frac{1}{q}, \inv{\theta}$ respectively.  
\end{lemma}

\begin{proof}
    The calculation here is routine, and left to the reader. 
\end{proof}


As shown by Wenzl and Tuba in \cite{Wenzl}, semisimple tensor categories with fusion rules similar to those of the orthogonal groups can be realized using the BMW categories. 
\begin{definition}
    Define for any $q \ne \pm 1$ the \textbf{quantum orthogonal category} \[O(N)_q = \begin{cases}
        Ab(\sms{BMW(q^2, q^{4r})}) & N = 2r + 1\\
        Ab(\sms{BMW(q, q^{2r-1}}) & N = 2r
    \end{cases}\]
\end{definition}

When $q$ is specialized to a root of unity, $O(N)_q$ has the same fusion rules as $Rep(O(N))$, while when for $q$ a root of unity the fusion ring becomes a quotient of the classical one with finite $\Z^+$ basis. More specifically, $q = e^{\frac{2\pi i}{4(k + 2r - 1)}}$ the category $O(2r + 1)_q$ is semisimple, and has $O(2r + 1)-O(k)$ fusion rules (using Daniel Copeland's terminology in \cite{DanielDude}). Similarly, for $q = e^{\frac{2\pi i}{4(k + 2r - 2)}}$ the category $O(2r)_q$ is semisimple, with $O(2r)-O(k)$ fusion rules. 

One can recover the categories $\mathcal{C}(\mathfrak{so}_{2r + 1}, k)^{vec}$ and $\mathcal{C}(\mathfrak{so}_{2r}, k)^{vec}$ from the categories $O(N)_q$ by de-equivariantizing. This is not surprising, since by Daniel Copeland's work in classifying ribbon categories with $SO(N)-O(k)$ fusion rules, which shows that the eigenvalues of the braid, together with the fusion rules, completely determine the category. One accomplishes this by adding in an isomorphism from the ``determinant'' representation, to the tensor unit. To this end, we define the following family of diagram categories. 

\begin{definition}\label{definitionOfSW}
    Given $q,\theta\in \C$ such that for any $m \in \set{1, \dots, N}$ $[m]_q\ne 0$ and $\theta + q^{1 - 2m}\ne 0$, the \textbf{special BMW algebra}, denoted $SW_{2r + 1}(q, \theta)$, is the $\C$-linear braided monoidal category with objects $\N_0$, which extends $BMW(q, i)$ by the additional arrow, and where the morphism $\att{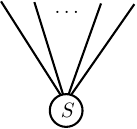}{0.5}$ is defined as the dual of $\att{extraGenerator}{0.5}$, and where these morphisms are subject to the relations 
    \[
    \text{(Braid)} \att{extraGeneratorBraided}{0.5} = \att{extraGeneratorUnbraided}{0.5} \hspace{0.5cm} \text{(Split)} \att{extraGeneratorProjection}{0.5} = P_{(1^N)}
    \]
    where $P_{(1^N)}$ is the minimal idempotent in $BMW(q, \theta)$ labeled by the Young diagram $(1^N)$.
\end{definition}

A few relations are immediately reduced from the relations. For one thing, note that in the BMW category, the idempotent $P_{(1^N)}$, when it exists, contains $P_{1^2}\tensor P_{(1)}^{\tensor (N -2)}$ with multiplicity $1$. Then $ \att{braid}{0.25} \circ P_{(1^2)} =P_{(1^2)}\circ \att{braid}{0.25} = -1/q P_{(1^2)}$. Putting these two things together, we have
\[
\att{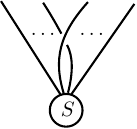}{0.5} = \att{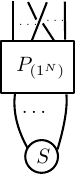}{0.5} = \att{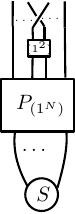}{0.5} = -1/q\att{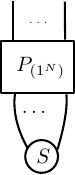}{0.5}.
\]

One can similarly deduce the relations
\[
\text{(Killed By Cup) }\att{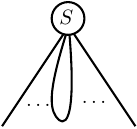}{1} = 0
\]

Note that by \cite{Wenzl} Lemma 7.6, there are restrictions on $q$ and $\theta$ in terms of $N$. If $q = e^{\frac{2\pi i}{4(k + \check{h})}}$ where $\check{h} = 2r - 1$ if $N$ is odd, and $\check{h} = 2r - 2$ if $N$ is even, then the parameters $q, \theta = q^{2r}$ satisfy the conditions necessary for the special BMW category to be defined. 

\begin{proposition}\label{specialBMWToC}\cite{Morrison}
    Let $r > 2$ For $q = e^{\frac{2\pi i}{4(k + 2r -1)}}$, one has a ribbon equivalence
    \[
    \mathcal{C}(\mathfrak{so}_{2r + 1}, k)^{vec} \simeq Ab\left(\sms{SW_{2r + 1}(q^2, q^{4r})}\right)
    \]
    where $\mathcal{C}(\mathfrak{so}_{2r + 1}, k)^{vec}$ has the unique pivotal structure that makes all of its quantum dimensions positive, which maps $V_{(1)}\mapsto 1\in Ab\left(\sms{SW_{2r + 1}(q^2, q^{4r})}\right)$. In particular, if $k = 2r - 1$, we make the change of variables $q = e^{\frac{2\pi i}{2(2r - 2)}}$, giving
    \[
    \mathcal{C}(\mathfrak{so}_{2r + 1}, 2r - 1)^{vec} \simeq Ab\left(\sms{SW_{2r + 1}(q, iq)}\right)
    \]
    Similarly, for $q = e^{\frac{2\pi i}{2(k + 2r - 2)}}$ one has a ribbon equivalence
    \[
    \mathcal{C}(\mathfrak{so}_{2r}, k)^{vec} \simeq Ab\left(\sms{SW_{2r}(q, q^{2r-1})}\right)
    \]
    which maps $V_{(1)}\mapsto 1\in Ab\left(\sms{SW_{2r}(q, q^{2r})}\right)$, and
    where $\mathcal{C}(\mathfrak{so}_{2r}, k)^{vec}$ has the unique pivotal structure that makes all of its quantum dimensions positive. In particular, if $k = 2r - 2$, one has
    \[
    \mathcal{C}(\mathfrak{so}_{2r}, 2r - 2)^{vec} \simeq Ab\left(\sms{SW_{2r}(q, iq)}\right)
    \]
\end{proposition}

Describing the categories $\mathcal{C}(\mathfrak{sp}_{2r}, k)$ is much simpler. 

\begin{proposition}\cite{Wenzl}
    For $k\in \N$, there is a braid equivalence
    \[
    \mathcal{C}(\mathfrak{sp}_{2r}, k) \simeq Ab(\sms{BMW(q, -q^{2r + 1})})
    \]
    where $q = e^{\frac{2\pi i}{4(r + k + 1)}}$, which maps $V_{(1)} \mapsto 1\in \mathcal{C}(\mathfrak{sp}_{2r}, k) \simeq Ab(\sms{BMW(q, -q^{2r + 1})})$. In particular, if $k = \check{h} = r + 1$, then there is a braid equivalence
    \[
    \mathcal{C}\left(\mathfrak{sp}_{2r}, r + 1\right) \simeq Ab\left(\sms{BMW\left(q, -\frac{i}{q}\right)}\right) 
    \]
    where $q = e^{\frac{2\pi i}{4(2r + 2)}}$.
    \end{proposition}

    \begin{remark}
        Note that $\mathcal{C}(\mathfrak{sp}_{2r}, k)$ has a pivotal structure such that the categorical dimension of the defining representation $V_{\Lambda_1}$ is positive. The problem is that this pivotal structure does not make $V_{\Lambda_1}$ symmetrically self-dual. One then changes the pivotal structure by multiplying the original pivotal structure by $-1$ for $V_{\lambda}$ when $\lambda$ has an odd number of boxes, and by $1$ for $V_\lambda$ when $\lambda$ has an even number of boxes. This allows one to use an unoriented skein theory for $\mathcal{C}(\mathfrak{sp}_{2r}, k)$, as in the above equivalence. See \cite{DanielDude} for the details of this process for any $\Z/2\Z$ graded category. 
    \end{remark}

    As mentioned above,
    $V_{(1)}^{\tensor 2} \cong V_{(0)} \oplus V_{(2)} \oplus V_{(1,1)}$ in any of the categories $C = \mathcal{C}(\mathfrak{so}_{2r + 1}, k)^{vec}$, $\mathcal{C}(\mathfrak{so}_{2r}, k)^{vec}$ or $\mathcal{C}(\mathfrak{sp}_{2r}, k)$. As a result, there exists three unique, primitive orthogonal idempotents in $End_\mathcal{C}(V_{(1)}^{\tensor 2})$ which sum to the identity map. Since $End_\mathcal{C}(V_{(1)}^{\tensor 2})$ is isomorphic to $BMW_2(q, \theta)$ where $q$ and $\theta$ are suitably chosen, it follows that these three idempotents must correspond to the idempotents $P_{(2)}$, $P_{(1,1)}$ and $P_{(0)}$. It should be pointed out that these idempotents do in fact project onto the simple objects labeled by the same Young diagrams, in any of the categories $ \mathcal{C}(\mathfrak{so}_{2r + 1}, k)^{vec}$, $\mathcal{C}(\mathfrak{so}_{2r}, k)^{vec}$ or $\mathcal{C}(\mathfrak{sp}_{2r}, k)$, so long as $k > 2$ and $r > 2$. 

\section{Conformal Embeddings}\label{ConformalEmbeddingsSection}

For any semisimple Lie algebra $\g$ and $k\in \N$, let $\mathcal{V}(\g, k)$ denote the WZW vertex operator algebra for $\g$ at level $k$. The language of vertex operator algebras and their representations provides a compact way of describing the category of level $k$ representations of the affine Lie algebra $\aff{\g}$. We will not dive into the details of VOAs. Due to the Kazhdan and Lutzig equivalence, the category of representations of $\V(\g, k)$, denoted $\mathcal{C}(\g, k)$, is braid equivalent to the category $\sms{Rep(U_q(\g))}$ where $q = e^{\frac{2\pi i}{4(k + \check{h})}}$.\cite{Kazhdan} Due to this equivalence, we will let $\mathcal{C}(\g, k)$ denote either category, and freely move between the two. 

Conformal embeddings of vertex operator algebras of the form $\mathcal{V}(\g, k) \sbs \V(\h, 1)$ give rise to a large class of examples of Etale algebras in the categories $\mathcal{C}(\g, k)$. Given a conformal embedding $\mathcal{V}(\g, k) \sbs \V(\h, 1)$, the restriction function $Res:\mathcal{C}(\h, 1) \to \mathcal{C}(\g, k)$ is lax monoidal, and preserves the braid. As a result, the unit map turns $A = Res(\mathbb{I})$ into an Etale algebra. 

Conformal embeddings are fully classified, and listed~\cite{WittGroup} Appendix A. A large class of conformal embeddings are of the form 
\[
\V(\g, \check{h}) \sbs \V(\mathfrak{so}_{\dim \g}, 1).
\]
In \cite{CainAndSnyder}, Edie-Michell and Noah Snyder use diagrammatic techniques to re-construct the category of modules $\mathcal{C}(\mathfrak{sl}_{r}, r^2 - 1)_A$ where $A$ is the Etale algebra coming from the conformal embedding $\V(\mathfrak{sl}_{r}, r + 1) \sbs \V(\mathfrak{so}_{r^2 + 1}, 1)$. Here, we consider conformal embeddings of the form
\begin{itemize}
    \item $\V(\mathfrak{sp}_{2r}, r + 1) \sbs \V(\mathfrak{so}_{r(2r + 1)}, 1).$
    \item $\V(\mathfrak{so}_{2r + 1}, 2r -1) \sbs \V(\mathfrak{so}_{r(2r + 1)}, 1)$
    \item $\V(\mathfrak{so}_{2r}, 2r - 2) \sbs \V(\mathfrak{so}_{r(2r - 1)}, 1)$
\end{itemize}
as well as the conformal embedding
\[
\V(\fsp_{2r + 2}, r) \sbs \V(\mathfrak{so}_{r(2r + 1)}, 1)
\]
which can be derived from the conformal embedding $\V(\mathfrak{sp}_{2r}, r + 1) \sbs \V(\mathfrak{so}_{r(2r + 1)}, 1)$ by level-rank duality.\cite{SymplecticLevelRank}

\begin{definition}
    For each of the above conformal embeddings $\V \sbs \V'$ listed above, let $C$ be the category of $\V$-modules. As noted, in each case we have the restriction functor $Res: \mathcal{C}(\so_{N}, 1) \to C$ for a suitable $N$.  Let $B = Res(V_{(1)})$ and let $A = Res(V_{(0)})$. The left unit map $\ell_{V_{(1)}}: \mathbb{I} \tensor V_{(1)}  \to V_{(1)}$ in $\mathcal{C}(\mathfrak{so}_{N}, 1)$ gives a module structure $Res({\ell_{V_{(1)}}}) : A\tensor B \to B$. Let $g$ denote $B$ as an $A$-module in $C$.
\end{definition}

Note that, for each of these conformal embeddings, we have a different category $C$, a different Etale algebra $A$, and a different $g$, with a different restriction functor $Res:\mathcal{C}(\so_N, 1) \to C$, all with a different $N$. Nonetheless, in order to not make our notation too cumbersome, whenever we write $A$ or $g$, it will be understood from context that these depend on the particular category in which $A$ lives.

Since $V_{(1)}$ as a $\mathbb{I}$ module is local, and since $Res$ is lax monoidal and braided, it follows that $g\in C_A$ is local. 

By Theorem 5.2 of \cite{qMacky}, the category of local modules over $A$ is completely understood. 
Restriction becomes a braid equivalence
\[
\mathcal{C}(\mathfrak{so}_{N}, 1) \to C_A^0
\]
where for each of the above conformal embeddings, where $N = \dim \g$.

As a result, we have the following
\begin{proposition}\label{gSimpleOrder2}
    For each of the above conformal embeddings, the object $g$ is simple, and invertible of order $2$. There is a half braid $\gamma$ on $g$ in $C_A$ which turns $g$ into a fermion. 
\end{proposition}
\begin{proof}
    Under the equivalence $C_A^{0} \simeq \mathcal{C}(\mathfrak{so}_N, 1)$, $g$ maps to $V_{(1)}$. We already have that $g$ is simple. The invertible simple objects in $\mathcal{C}(\mathfrak{so}_N, k)$ are fully understood for any $k\in \N$. 
    
    If $N$ is odd, by \cite{Fuchs} there are only two of invertible simple objects in $\mathcal{C}(\mathfrak{so}_{N},1)$, namely $V_{(1)}$ and $V_{(0)}$. The former clearly is invertible of order $2$, its inverse cannot be $V_{(0)}$. Since $k = 1$, $g = V_{(1)}$ must be invertible of order $2$.  
    If $N$ is even and not less than $4$, by \cite{Fuchs} there are four invertible simple objects in $\mathcal{C}(\mathfrak{so}_{N}, 1)$, and these are
    \[
    V_{(0)}, V_{(1)}, V_{\Lambda_{N/2 -1}}, V_{\Lambda_{N}/2}.
    \]
    The latter two of these have spinor highest weights with odd multiplicity. Tensoring $V_{(1)}$ with a spinor weight with odd multiplicity cannot result in $V_{(0)}$, from which it follows that $V_{(1)}$ must be its own inverse. 

    While there are many ways of verifying that the braid on $g$ turns $g$ into a fermion in $C_A^0$, the easiest is to use the formula for the twist $\theta_g$ on $g$ (which is the defining representation in $\mathcal{C}(\so_N, 1)$). The braid $g$ extends to a half braid $\gamma$ on $g$ in $C_A$ by Corollary 3.30 of \cite{WittGroup}.
\end{proof}

For conformal embeddings of the form $\V(\g, \check{h}) \sbs \V(\mathfrak{so}_{\dim g}, 1)$, the decompositions of the objects $A$ and $B$ are described in \cite{decompositionRules}. We repeat this here, translated into our notation, in the below proposition.

      Let $\g$ be any semisimple Lie algebra with root system $\Phi$ with a fixed choice of positive roots $\Phi^+$. 

      An abelian ideal of $\Phi^+$ is a subset $\Psi \sbs \Phi^+$ such that 
\begin{enumerate}
    \item $(\Psi + \Phi^+)\I \Phi^+ \sbs \Psi$
    \item $(\Psi + \Psi)\I \Phi^+ = \emptyset$. 
\end{enumerate}

\begin{definition}\cite{decompositionRules}
    Any set $\Psi\sbs \Phi^+$ satisfying condition (1) is called an \textbf{ideal}, and any set $\Psi$ satisfying both (1) and (2) is called \textbf{abelian ideal}. When there is no confusion about which root system and choice of positive roots is being discussed, let $\Sigma$ denote the set of all abelian ideals.
\end{definition}

\begin{remark}
    A choice of simple roots $\Delta$ induces a partial ordering $<$ on the $\Phi^+$, where $\alpha < \beta$ if and only if $\beta - \alpha$ is a non-negative integer combination of roots in $\Delta$. Condition (1) amounts to saying that $\Psi \sbs \Phi^+$ is a downward closed subset of $\Phi^+$ under the partial ordering $\le$. Hence ideals of $\mathfrak{b}$ are in bijection with downward closed subsets of $\Phi^+$.
\end{remark}

For any $\Psi \in \Sigma$, form the weight $\lambda_\Psi = \sum_{\alpha\in \Psi} \alpha$. As shown in \cite{decompositionRules}, this weight is dominant. 

By Theorem 3.9 of \cite{decompositionRules}, we have the following 
\begin{proposition}\cite{decompositionRules}
For any conformal embedding of the form $\V(\g, \check{h}) \sbs \V(\mathfrak{so}_{\dim \g}, 1)$, the objects $A$ and $B$ described above decompose as
    \[
A = \directSumOver{\Psi \in \Sigma\\\text{$|\Psi|$ even}} V_{\lambda_{\Psi}}, \hspace{0.5cm} B= \directSumOver{\Psi \in \Sigma\\\text{$|\Psi|$ odd}} V_{\lambda_{\Psi}}
\]
\end{proposition}

Note that for the root systems $B_r$, $D_r$, all of the weights $\lambda_\Psi$ are clearly integral weights. As a result, we have $A$ and $B$ as objects in $\mathcal{C}(\mathfrak{so}_{2r + 1}, 2r - 1)^{vec}$ and $\mathcal{C}(\mathfrak{so}_{2r}, 2r - 2)^{vec}$. It follows that $A$ is an Etale algebra in either of these categories, and $g$ is a local module in either of these categories. 

\begin{proposition}\label{CASemisimple}
    For $\mathfrak{g} = \mathfrak{so}_{2r + 1}$, $\mathfrak{sp}_{2r}$ and $\mathfrak{so}_{2r}$, and for any $k\in \N$, the categories $\mathcal{C}(\mathfrak{g}, k)$ are unitary. Moreover, $\mathcal{C}(\mathfrak{so}_{2r + 1}, k)^{vec}$ is unitary. Hence the categories $C_A$ are unitary for $C = \mathcal{C}(\mathfrak{so}_{2r + 1}, 2r - 1)^{vec}, \mathcal{C}(\mathfrak{sp}_{2r}, r + 1), \mathcal{C}(\mathfrak{so}_{2r}, 2r - 2)^{vec}$. 
\end{proposition}

\begin{proof}
    Note that if $\Psi\in \Sigma$, with $\Psi \ne \emptyset$, $\lambda_\Psi \ne V_{(0)}$. Obviously, $\emptyset\in \Sigma$ is an abelian ideal of even size. Hence $V_{(0)}$ appears exactly once in the decomposition of $A$, and hence for any of the categories listed above, $Hom_\mathcal{C}(A \to V_{(0)}) \cong \C$. In other words the algebra $A$ is connected (or, in the language of \cite{ConnectedAlgebraUnitary}, $A$ is Haploid). It was shown in \cite{AffineCatsUnitary} that $C = \mathcal{C}(\mathfrak{so}_{2r + 1}, 2r - 1), \mathcal{C}(\mathfrak{sp}_{2r}, r + 1), \mathcal{C}(\mathfrak{so}_{2r}, 2r - 2)$ are unitary categories. Also, $C = \mathcal{C}(\mathfrak{so}_{2r + 1}, 2r - 1)^{vec}, \mathcal{C}(\mathfrak{so}_{2r}, 2r - 2)^{vec}$ are easily seen to be unitary as tensor subcategories. Since $A$ is connected and seperable (it is an Etale algebra), $A$ is also rigid. It was shown in \cite{ConnectedAlgebraUnitary} that this condition ensures $A$ (in each of these categories) is a $Q$-system, and hence $C_A$ is a unitary category for each $C$ mentioned above.

\end{proof}

We will now make the decompositions of the objects $A$ and $B$ for the adjoint conformal embeddings of types B, C, and D even more explicit. For type C, this will establish a dimension formula for certain endomorphism algebras, which will allow us to establish the Deligne interpolation result listed in the introduction. As a result, our analysis of type C is the most important, while our analysis of types B and D is of secondary importance. 

It appears that our geometric descriptions of the abelian ideals is not new in principle. Although we cannot find every detail in the literature, similar patterns appear, for example, in Panyushev's work in \cite{Panyushev}. The correspondence in type C between subobjects of the Etale algebra $A$ and shifted tableaux via doubling may be known to experts, but we are unable to find any instance of it in the literature. The situation is similar for the partial correspondences in types B and D. 

\subsection{Type B}

The root system $B_r$ is the collection of vectors in $\Z^r \sbs \R^r$ with norm $1$ or $\sqrt{2}$. These are all sums of the standard unit vectors $e_i\in \R^r$.

Our choice of positive roots will be $\Phi^+ = \set{e_i}\cup \set{e_i \pm e_j}_{i < j}$, with choice of simple roots 
\[
\Delta = \defset{\alpha_i = e_i - e_{i + 1}}{1 \le i \le r - 1}\U \set{\alpha_r = e_r}.
\]
In this basis, we have
\[
e_i = \alpha_i + \dots + \alpha_r
\]
for $i < n$, 
\[
e_i - e_j = \alpha_i + \dots + \alpha_{j -1}
\]
for $1\le i < j\le r$, and 
\[
e_i + e_j = \alpha_i + \dots + \alpha_{j - 1}  + 2\alpha_j + \dots + 2\alpha_n
\]
for $1\le i < j\le r$.

A convenient way of describing the roots of $B_r$ is by using them to label the squares of a shifted Young diagram. From now on, let $r$ be a fixed integer $r \ge 2$. 

\begin{definition}
Let $T$ be the labeling of $\Phi^+$ given by the function $T: \Phi^+ \to \N^2$. 
\begin{align*}
T:& e_i \mapsto (r, i)\\
    & e_i - e_j \mapsto (2r + 1 - j, i) \\
    & e_i + e_j \mapsto (j - 1, i)\\
\end{align*}
\end{definition}

Note that $T$ can be viewed as a shifted Young Tableaux, labeling the shifted Young diagram corresponding to $\mu = (2r - 1, 2r - 3,\dots, 1)$ with the positive roots $\Phi^+$. We illustrate two examples, where $r = 2$ and $r = 3$,

\begin{align*}
r = 2 \rightsquigarrow
\scalebox{0.8}{
\ytableausetup{boxsize = 4em}$
\begin{ytableau}
e_1 + e_2 & e_1 & e_1 - e_2\\
\none & e_2
\end{ytableau}
$}
\hspace{1cm} r = 3 \rightsquigarrow \hspace{0.5cm} \scalebox{0.5}{
\ytableausetup{boxsize = 4em}$
\begin{ytableau}
e_1 + e_2 & e_1 + e_3 & e_1 & e_1 - e_3 & e_1 - e_2 \\
\none & e_2 + e_3 & e_2 & e_2 - e_3 \\
\none & \none & e_3
\end{ytableau}
$}
\end{align*}

\begin{lemma}\label{idealsTypeB}
    A subset $\Psi \sbs \Phi^+$ is an ideal if and only if $T(\Psi) \sbs \N^2$ is a shifted Young diagram. 
\end{lemma}
\begin{proof}
    Suppose that $\Psi$ is an ideal. We must show that $T(\Psi)$ is a shifted Young diagram. It is enough to show that for each $\alpha \in \Psi$, then any admissable squares in $T(\Phi^+)$ directly above and to the left are also in $T(\Phi^+)$ (the positive $y$ direction is ``down,'' the positive $x$ direction is to the ``right''). We must now analyze each possible form for $\alpha \in \Psi$.
    \begin{enumerate}
        \item[$\alpha = e_i$] Suppose $\alpha = e_i\in \Psi$. Then $T(\alpha) = (r, i)$. If $i = 1$, then $(r, i) = (r, 1)$ has no upward neighbors in $T(\Phi^+)$. If $i > 1$, then by definition of an ideal, we have $e_i + \alpha_i = e_i + e_{i - 1} - e_i= e_{i - 1}\in \Psi$, and hence $T(e_{i - 1}) = (r, i - 1)\in T(\Psi)$. 
        
        If $i = r$, then $T(e_i)$ has no neighbors to the left. If $i < r$, then since $e_{r}\in \Phi^+$ and $\Psi$ is an ideal, we have $e_i + e_{r} \in \Psi$, and $T(e_i + e_{r}) = (r-1, i)\in T(\Psi)$, which is the left neighbor of $T(e_i) = (r, i)$.

        In all cases, if $T(e_i)$ has a left or top neighbor, then that neighbor is also in $T(\Psi)$.
        \item[$\alpha= e_i + e_j$] Suppose $\alpha = e_i + e_j$ for $i < j$. Then $T(e_i + e_j) = (j - 1, i)$. 
        
        If $j - 1 = i$, then $T(e_i + e_j)$ has no left neighbors. Otherwise, we have $i < j - 1$, and hence $e_i + e_{j - 1} = e_i + e_j + e_{j -1} - e_j \in \Psi$, and $T(e_i + e_{j - 1}) = (j - 2, i)$ is the left neighbor of $T(e_i + e_j)$.

        \item[$\alpha = e_i - e_j$] Suppose $\alpha = e_i - e_j$ for $i < j$. Then $T(e_i - e_j) = (2r + 1 - j, i)$. 
        
        If $j = r$, then $e_i - e_r + e_r = e_i\in \Psi$, and $T(e_i) = (r, i) = (2r + 1 - r - 1,i)$ is the left neighbor of $e_i - e_j$. If $j < r$, then $e_j - e_{j + 1}\in \Phi^+$, and hence $e_i - e_{j + 1} = e_i - e_j + e_j - e_{j + 1} \in \Psi$. Note that $T(e_i - e_{j + 1}) = (2r + 1 - j -1, i)$, which is the left neighbor of $e_i - e_j$. In either case, the left neighbor of $e_i - e_j$ is also in $\Psi$. 

        If $i = 1$, then $T(e_i - e_j)$ has no top neighbors. If $i > 1$, then $e_{i - 1} - e_i \in \Phi^+$, and hence $e_{i - 1} - e_j = e_i - e_j + e_{i - 1} - e_i \in \Psi$. Note that $T(e_{i-1} - e_j)$ is the top neighbor of $T(e_i - e_j)$.
    \end{enumerate}

    It is now clear that if $\Psi$ is an ideal, then $T(\Psi)$ is a shifted Young diagram. 
    
    Now we must show that for $\Psi \sbs \Phi^+$, $T(\Psi)$ being a shifted Young diagram ensures that $\Psi$ is an ideal. Suppose that $T(\Psi)$ is a shifted Young diagram. In order to show that $\Psi$ is an ideal, it is enough to show that $(\Psi + \Delta)\I \Psi^+ \sbs \Psi$.
    
    Suppose $e_i\in \Psi$ for some $i$. Then the only simple roots $\alpha$ such that $e_i + \alpha\in \Phi^+$ are $\alpha = e_r$ (if $i < r$) and $\alpha = e_{i - 1} - e_i$ (if $i > 1$). If $i < r$, then $T(e_i) = (r, i)$ and $T(e_i + e_r) = (r - 1, i)$, which is the left neighbor of $T(e_i)$, and hence $e_i + e_r\in \Psi$.If $i > 1$, then $e_i + e_{i - 1} - e_i = e_{i -1}$, and $T(e_{i - 1}) = (r, i-1)$ which is the top neighbor of $T(e_i)$, and hence $e_i + e_j$. This shows that for roots of the form $e_i\in \Psi$, $(e_i + \Delta)\I \Phi \sbs \Psi$.

    A similar analysis works for roots of the form $e_i - e_j$ and $e_i + e_j\in \Psi$.
\end{proof}

\begin{corollary}
    A set $\Psi \sbs \Phi^+$ consisting of only of positive roots in the upper half of the diagram, that is, of the form $e_i + e_j$, automatically satisfies $(\Psi + \Psi) \I \Phi^+ = \emptyset$.
\end{corollary}
\begin{proof}
    This is immediately obvious from noting that $e_i + e_j + e_k + e_\ell$ either has norm $\sqrt{4}$, $\sqrt{5}$, or $\sqrt{8}$, and hence none of these can be positive roots. 
\end{proof}

Following after \cite{CainAndSnyder}, as shown in \cite{decompositionRules}, the summands for $A$ correspond to these abelian ideals in the following way.

\begin{lemma}
    A set $\Psi \sbs \Phi^+$ has $\sum_{\alpha \in \Psi} \alpha = x_1 e_1 +\dots +x_N e_N$ for unique $x_k \in \Z$. One then obtains a list of $N$ integers, $\lambda_\Psi = (|x_1|, |x_2|, \dots, |x_n|)$. This list of integers $\lambda_\Psi$ is weakly increasing in length, and hence defines a Young diagram. 
\end{lemma}

Unfortunately, while the ideals in $\mathfrak{sp}_{2r + 1}$ were given a satisfying geometric interpretation, it is much harder to give an nice description of the correspondence between abelian ideals and the Young diagram parameterizing the corresponding weight. Luckily, so long as $T(\Psi)$ has no more than $r$ columns, we have the following

\begin{lemma}\label{abelianIdealsToYoungDiagrams}
    Suppose $\Psi \sbs \Phi^+$ consists only of roots of the form $e_i + e_j$ (equivalently, $T(\Psi)$ has no more than $r- 1$ columns). Then $\lambda_\Psi = (\lambda_1, \dots, \lambda_r)$ is partition specified by the doubled Young diagram $\lambda_\Psi = D(T(\Psi))$. Also, $|\lambda_\Psi| \ge |T(\Psi)_1|$.
\end{lemma}
\begin{proof}
    The number of $e_1$'s appearing in $\Psi$ is the length of the first row of $T(\Psi)$, and hence $(\lambda_\Psi)_1 = T(\Psi)_1$. Let us now count the number of $e_i$'s for $1 < i \le r$. The labeling $T$ maps roots of the form $e_k + e_i \mapsto (i - 1, k)$, and maps roots of the form $e_i + e_k$ to $(k, i)$. Hence the number of $e_i$'s appearing in $T(\Psi)$ are either in the $i-1$'th column or $i$th row of $T(\Psi)$, and every such box contributes exactly one $e_i$. Hence there are exactly $T(\Psi)_i + T(\Psi)^t_{i - 1}$ $e_i$'s appearing in $T(\Psi)$. Hence $(\lambda_\Psi)_1 = T(\Psi)_i + T(\Psi)^t_{i -1} = D(T(\Psi))$ for $i > 1$. 

    To see that $|\lambda_\Psi| \ge |T(\Psi)_1|$, note that there are no positive roots with a negative $e_1$ component. 
\end{proof}

Let us now put this combinatorics to work. For the remainder of this section, let $C = \fanC(\mathfrak{so}_{2r + 1}, 2r - 1)^{vec}$ and let $C_A$ denote the category of $A$-modules.

\begin{lemma}\label{gInEnd2B}
    For $r \ge 3$, $\dim Hom_{C_A}(X\tensor X \to g) = 1$. 
\end{lemma}
\begin{proof}
    Note that $V_{\Lambda_1} \tensor V_{\Lambda_1} = V_\emptyset \oplus V_{(2)} \oplus V_{(1, 1)}$. 
   
    By the free forgetful adjunction (and using the tensorator), we have
    \begin{align*}
    Hom_{C_A}(X\tensor X \to g)\\
    \cong Hom_{C_A}(F(V_{\Lambda_1} \tensor V_{\Lambda_1}) \to g) \\
    \cong Hom_{C}(V_{\Lambda_1} \tensor V_{\Lambda_1} \to B)\\
    \cong Hom_{C}(V_\emptyset \to B) \oplus Hom_{C}(V_{(2)} \to B) \oplus Hom_{C}(V_{(1,1)}\to B)
    \end{align*}

    Note that if $\Psi\sbs \Phi^+$ is an abelian ideal, with $T(\Psi)$ having more than $r - 1$ columns, then $(\lambda_{\Psi})_1 > r > 2$. On the other hand, if $T(\Psi)$ has less than $r - 1$ columns, then $\lambda_{\Psi} = D(T(\Psi))$, so in particular $\lambda_{\Psi}$ must be labeled by a doubled Young diagram, and of a strict partition of odd size. Only $(1,1)$ fits this description, and hence only $V_{(1,1)}$ appears as a summand in $B$. Hence $\dim Hom_{C_A}(X\tensor X \to g) = 1$.
\end{proof}

We will also have use for a limiting formula for the dimensions of the endomorphism algebras in $C_A$ which stabilize as $ r\to \infty$. 

\begin{lemma}\label{stableDimB}
    For all $n \in \N$ there is some $L \in \N$ such that whenever $ r > L$, we have 
    \[
    \dim End_{C_A}(X^{\tensor n})  = 2^{n - 1}(2n - 1)!!
    \]
\end{lemma}

\begin{proof}
    Choose $L$ large enough so that for every abelian ideal $\Psi \in \Sigma$ such that $T(\Psi)$ has more than $L$ rows, $|\lambda_{\Psi}| > 2n$. Then by the free forgetful adjunction and self duality of $X$,
    \begin{align*}
    \dim End_{C_A}(X^{\tensor n}) = \dim End_{C}(V_{(1)}^{\tensor 2n} \tensor A \to V_{(0)}) \\
    = \#\left\{ \text{Paths on the fusion graph for $\mathcal{C}(\so_{2r + 1}, 2r - 1)$ relative to $V_{(1)}$ from summands of $A$ to $V_{(0)}$}\right\}.
    \end{align*}
    Choose $L$ large enough so that fusion graph for $\mathcal{C}(\so_{2r + 1}, 2r - 1)$ relative to $V_{(1)}$ agrees with the Young graph for all Young diagrams with size at most $2n$.
    
    If an abelian ideal $\Psi$ has $T(\Psi)$ with more than $L$ rows then no path on the fusion graph of length $2n$ from $V_{(0)}$ to $\lambda_\Psi$ is not on the Young graph. Moreover, since $r > L$, if $T(\Psi)$ has no more than $L$ rows then $\lambda_{\Psi} = D(T(\Psi))$. Hence $\dim End_{C_A}(X^{\tensor n})$ counts the number of paths on the Young graph from doubled shifted Young diagrams of length $2n$. By Lemma~\ref{evenShiftedYoung} this is $2^{n - 1}(2n - 1)!!$.
\end{proof}

The purpose of the above lemma will be to show that $End_{C_A}(X^{\tensor n})$ is eventually isomorphic to the even part of the Clifford BMW superalgebra as $r \to \infty$. 

Finally, the following lemma will be essential in proving the consistency of our labelings of the representations of the endomorphism algebras in the final section. Roughly speaking, this lemma shows that while $C_A$ is not $\Z/2\Z$ graded, it ``eventually is.''

\begin{lemma}\label{eventuallyZ2Graded}
    For all $n \in \N$ there exists some $r \in \N$ such that for all $m \le n$ 
    \[
    \dim Hom_{C_A} (X^{\tensor 2m + 1} \to A) = 0.
    \]
\end{lemma}
\begin{proof}
    For abelian ideals $\Psi \in \Sigma$ such that $T(\Psi)$ has less than $r$ columns, $\lambda_\Psi = D(T(\Psi))$. Now choose $r$ large enough so that the fusion graph for $C = \mathcal{C}(\so_{2r + 1}, 2r - 1)$ relative to $V_{(1)}$ is the Young graph up to level $2n + 1$. 

    If $T(\Psi)$ has at least $r$ columns, then $|\lambda_\Psi| \ge r$. Choose $r$ large enough so that $r > 2n + 1$. Then for any abelian ideal $\Psi\in \Sigma$ such that $T(\Psi)$ has at least $r$ columns, $\lambda_\Psi$ does not show up in $V_{(1)}^{2m + 1}$ for any $m \le n$. 

    By the free forgetful adjunction, $\dim Hom_{C_A} (X^{\tensor 2m + 1} \to A) \cong \dim Hom_C(V_{(1)}^{\tensor 2m + 1} \to V_{(0)})$ counts paths on the fusion graph from $V_{(0)}$ to simple summands of $A$ of length $2m + 1$. By our assumptions on $r$, every such path terminates at some $\lambda_\Psi$ for $\Psi$ an abelian ideal such that $T(\Psi)$ has no more than $r$ columns. Hence $\lambda_\Psi = D(T(\Psi))$ and in particular $|\lambda_\Psi|$ is even. By our assumptions on $r$, every path of length $2m + 1$ on the fusion graph for $C$ starting at $V_{(0)}$ is a path on the Young graph $(0) \to \dots \to \lambda_\Psi$ of length $2m + 1$. But the Young graph is a $\Z/2\Z$ graded graph, and hence no path of odd length from $(0)$ to an even Young diagram can exist. 
\end{proof}

\subsection{Type C}

Type C is the most important case to analyze, and it is also the nicest. 

A particular construction of the root system $C_r \sbs \R^r$ consists of all integer vectors of the form $e_i \pm e_j$ or $1 \le i, j \le r$, which can be obtained by taking the dual of $B_r$. Our choice of simple roots shall be
\[
\Delta = \set{\alpha_1 = e_1 - e_2, \dots, \alpha_{r-1} = e_{r -1} - e_r, \alpha_r = 2e_r}
\]

The positive roots $\Phi^+$ are easily enumerated as 
\[
\Phi^+ = \defset{e_i + e_j}{1\le i < j\le r}  \U \defset{1 \le e_i - e_j\le r}{ i < j} \U \defset{2e_i}{1\le i \le r}.
\]

\begin{definition}
Let $T$ be the labeling of $\Phi^+$ given by the function $T: \Phi^+ \to \N^2$. 
\begin{align*}
T:  & 2e_i \mapsto (i, i)\\
    & e_i + e_j \mapsto (j, i)\\
    & e_i - e_j \mapsto (2r - 1 - j, i) 
\end{align*}
\end{definition}

\begin{example}
    
\begin{align*}
r = 2 \rightsquigarrow
\scalebox{0.8}{
\ytableausetup{boxsize = 4em}$
\begin{ytableau}
2e_1& e_1 + e_2 & e_1 - e_2\\
\none & 2e_2
\end{ytableau}
$}
\hspace{1cm} r = 3 \rightsquigarrow \hspace{0.5cm} \scalebox{0.5}{
\ytableausetup{boxsize = 4em}$
\begin{ytableau}
2e_1 & e_1 + e_2 & e_1+ e_3 & e_1 - e_3 & e_1 - e_2 \\
\none & 2e_2& e_2 + e_3 & e_2 - e_3 \\
\none & \none & 2e_3
\end{ytableau}
$}
\end{align*}
\end{example}

The following lemma, while similar to the type B case, is essential for our Deligne interpolation result, which in our opinion is the main result of this paper. Due to this lemma's importance, we will include every detail of the following proof, even at the cost of some repetition from the type B case.

\begin{lemma}\label{idealsTypeC}
The ideals $\Psi\sbs \Phi^+$ are exactly the subsets $\Psi \sbs \Phi^+$ such that $T(\Psi)$ is a shifted Young diagram. Moreover, for every ideal $\Psi$, $\lambda_\Psi$ is a Young diagram with an even number of boxes. 
\end{lemma}
\begin{proof}
    The proof will be similar to the type B case. First suppose that $\Psi$ is an ideal. We must show that $T(\Psi)$ is a shifted Young diagram. Let us analyze each case for a type of root $\alpha \in \Psi$, and show that whenever $T(\alpha)$ has an up or left neighbor, that neighbor is also in $\Psi$. 

    \begin{itemize}
        \item[$\alpha = 2e_i$] If $i = 1$, then $T(2e_i) = (1,1)$, which has no upward or left-neighbors. If $i > 1$, then $T(2e_i)$ has no left-neighbors. On the other hand, $T(2e_i)$'s upward neighbor is $e_{i - 1} + e_i$. Note that $e_{i - 1} - e_i \in \Delta$, so since $\Psi$ is an ideal, we must have $2e_i + e_{i - 1} - e_i = e_{i - 1} + e_i \in \Psi$.
        \item[$\alpha = e_i + e_j$]. Suppose now that $\alpha = e_i + e_j$ where $i < j$. The left-most neighbor of $T(\alpha)$ is $e_i + e_{j - 1}$ (note that it may be the case that $i = j - 1$). Note that $j > i > 1$, hence $e_{j - 1}  - e_j\in \Delta$, and hence $e_i + e_j + e_{j - 1}  - e_j = e_i + e_{j -1} \in \Psi$. 

        If $i = 1$, then $\alpha = e_i + e_j$ has no upward neighbor. Otherwise, $\alpha$'s upward neighbor is $e_{i - 1} + e_j$. Note that $e_{i - 1} - e_i \in \Delta$, hence $e_{i} + e_j +e_{i - 1} - e_i = e_{i - 1} + e_j \in \Psi$.
        \item[$\alpha = e_i - e_j$] If $i > 1$, then $e_i - e_j$'s upward neighbor is $e_{i - 1} - e_j$, and since $e_{i - 1} - e_j \in \Delta$, $e_{i - 1} - e_j = e_{i - 1} - e_j +  e_i - e_j$. If $i = 1$, then $e_i - e_j$ has no upward neighbor. 

        If $j = r$, then the left neighbor of $\alpha = e_i - e_r$ is $e_i + e_r$, and since $2e_r\in \Delta$, $e_i + e_r = e_i - e_r + 2e_r\in \Psi$. 

        If $j < r$, then the left neighbor of $\alpha$ is $e_i - e_{j + 1}$. In this case, $e_j - e_{j - 1}\in \Delta$, and $e_i - e_{j + 1} = e_i - e_j + e_j - e_{j - 1} \in \Psi$.
    \end{itemize}

    This shows that for any $\Psi\sbs \Phi^+$ an ideal, $T(\Psi)$ is a shifted Young diagram. To show that every shifted Young diagram contained in $T(\Phi^+)$ corresponds to an ideal, suppose that $\mu \sbs T(\Phi^+)$ is a shifted Young diagram. One then observes that if $\alpha \in \Delta$ and $\beta \in \Psi$ with $\alpha + \beta \in \Phi^+$, then $T(\alpha + \beta)$ is always above or to the right of $T(\beta)$. So, since $\mu$ is a shifted Young diagram, it must be that $\alpha +\beta \in \Psi$, and hence $\defset{\beta \in \Phi^+}{ T(\beta) \in \mu}$ forms an ideal. 

    Finally, we must show that for each ideal $\Psi\sbs \Phi^+$, $|\lambda_\Psi|$ is even. To show this, we induct on the number of boxes in $T(\Psi)$ (or equivalently, the number of roots in $\Psi$). Note that if $|\mu| = \emptyset$, there is nothing to show.  Suppose that this holds for shifted Young diagrams $\mu\sbs T(\Phi^+)$ such that $|\mu| < n$. Suppose that $\mu$ is a shifted Young diagram of size $n$. If $|\mu| = 0$, then we are done. Otherwise, we can remove at least one box from $\mu$. We now analyze each case for the form of the root corresponding to this box. Let $\Psi\sbs \Phi^+$ be the ideal such that $T(\Psi) = \mu$, and let $\Psi'$ be the ideal obtained by removing $\alpha$. 

    Suppose first that the removed box corresponds to a root fo the form $\alpha = e_i - e_j$. Then $|\lambda_{\Psi'}|$ is even by the induction hypothesis. $\lambda_{\Psi}$ is obtained from $\lambda_{\Psi}$ by adding a box in the $i$th row, and removing a box in the $j$th row. Hence $|\lambda_{\Psi'}| = |\lambda_{\Psi}|$, and in particular $|\lambda_{\Psi}|$ is even. 

    Now suppose that the removed box is labeled by a root of the form $\alpha = e_i + e_j$. Then $\lambda_\Psi$ is obtained from $\lambda_{\Psi'}$ by adding a box in the $i$th row, and another in the $j$th row. Hence $|\lambda_{\Psi}| = |\lambda_{\Psi'}| + 2$, and $|\lambda_{\Psi}|$ is even. 

    Note that in the previous case, we did not assume that $i < j$. Hence the case where $\alpha = 2e_i$ has already been addressed. 

    This proves the lemma. 
\end{proof}

\begin{definition}
    Let $D = \set{e_i + e_j}_{i \le j}$. 
\end{definition}

Note that $D \sbs \Phi^+$ is an abelian ideal. It is trivial to show that if $\Psi$ is an abelian idea, and $\Omega \sbs D$ with $\Omega\sbs \Phi^+$ an ideal, then $\Omega$ is abelian.\footnote{In other words, abelian ideals are downward closed in the lattice of ideals.} It turns out that the converse is also true. Note that $\Psi \sbs D$ if and only if $T(\Psi)$ has no more than $r$ columns.

\begin{proposition}\label{sameAsDoubledC}
Under the correspondence $\Psi \mapsto \lambda_\Psi$, for any ideal $\Psi \sbs \Phi^+$ such that $T(\Psi)$ has no more than $r$ columns, $\Psi$ is abelian, and $\lambda_\Psi = D(T(\Psi))^t$. Hence the simple subobjects of $A$ are in bijection with shifted Young diagrams with no more than $r$ rows and columns, where the bijection is given by doubling and then transposing.
\end{proposition}
\begin{proof}
    Let $\Psi \sbs \Phi^+$ such that $T(\Psi)$ has no more than $r$ columns. Then every root in $\Psi\sbs D$ and hence $\Psi$ is abelian.
    
    Now notice that since every root in $\Psi$ is of the form $e_i + e_j$, we have that $(\lambda_{\Psi})_i$ is the number of times $e_i$ appears in $\Psi$. This can happen in two ways, either as $e_i+ e_j$ for $i \le j$, or for $e_j + e_i$ for $j \le i$. Under the labeling $T$, this corresponds to the $i$-th row, and the $i$-th column of $T(\Psi)$. Note that for the intersection of the $i$-th column and the $i$-th row, we have the root $2e_i$, which accounts for two $e_i$'s. Hence the number of $e_i$'s that appear in $\Psi$ is $T(\Psi)_i + T(\Psi)_i^T$ (the $i$-th column and the $i$-th row. In other words, $(\lambda_{\Psi})_i = T(\Psi)_i + T(\Psi)_i^T = D(T(\Psi))^T_i$ for all $i$, and hence $D(T(\Psi))^T = \lambda_{\Psi}$.

    The proposition then follows from Lemma~\ref{idealsTypeC}. 
\end{proof}

It is now time to put these results to use.

As before, let $X = F_A(V_{(1)})$ and let $g$ be as before.

\begin{lemma}\label{gInEnd2C}
For $r \ge 3$, $\dim Hom_{C_A}( X \tensor X \to g) = 1$. Hence there exists a unique minimal central idempotent $P_g\in End_{C_A}(X\tensor X)$ whose image is $g$. 
\end{lemma}
\begin{proof}
    Since $r \ge 3$, $V_{\Lambda_1}\tensor V_{\Lambda_1} = V_{(0)}\oplus V_{(2)}\oplus V_{(1,1)}$. By the free forgetful adjunction we have 
    \begin{align*}
        Hom_{C_A}( X \tensor X \to  g) \cong Hom_{C}(V_\emptyset \to B) \oplus Hom_{C}(V_{(2)} \to B) \oplus Hom_{C}(V_{(1,1)}\to B)
    \end{align*}
    If an abelian ideal $\Psi$ is not contained in $D$, then $T(\Psi)$ has more than $r$ rows, so that by Theorem~\ref{sameAsDoubledC} $|\lambda_{\Psi}| \ge 3$. Out of $(0), (2), (1, 1)$, only $(2)$ is the transpose of a Young diagram of odd size. This proves the lemma.
\end{proof}

Finally, we prove the following dimension formula, which will be essential for the proof of our Deligne interpolation result for the Clifford BMW category.

\begin{theorem}\label{dimensionOfCAThypeC}
    For $n < r/2$, we have 
    \[
    \dim End_{C_A}( X^{\tensor n}) = 2^{n - 1}(2n - 1)!!
    \]
\end{theorem}
\begin{proof}
    As in the type B case, we have 
    \[
    \dim End_{C_A}(X^{\tensor n}) = \#\set{\begin{array}{c}
\text{paths of length $2n$ from $V_{\lambda_\Psi}$ to $V_{(0)}$ in the fusion graph for $C$}\\
\Psi \in \Sigma\\
|\Psi| \equiv 0\pmod{2}
\end{array}}
    \]
     The fusion graph for $C$ is the truncation of the Young graph, where all Young diagrams with more than $r + 1$ columns and more than $r$ rows are removed. As a result, we care counting the number of paths of length $2n$ on the Young graph $Y$ from $(0)$ to $\lambda_{\Psi}$ which do not pass through Young diagrams with more than $r$ rows or $r + 1$ columns. Since $n < r/2 < r < r + 1$ all Young diagrams of size at most $n$ satisfy this condition. If a Young diagram $\lambda$ appears on a path from $(0)$  length $2n$, then $|\lambda| \le 2n$. Hence $\lambda_1, \lambda_1' \le 2n < r, r + 1$. As a result any path of length $2n$ starting and ending at $(0)$ is also in the fusion graph of $C$. If $\mu$ is a shifted Young diagram with more than $r$ rows, then $|D(\mu)^T| = |D(\mu)| = 2|\mu| \ge 2r > 4n$, and hence there cannot exist a path of length $2n$ from $D(\mu)^T$ to $(0)$ on the Young graph. As a result, every path on the Young graph from $D(\mu)^T$ to $(0)$ of length $2n$ is a path from $\lambda_\Psi$ to $(0)$ of length $2n$ in the fusion graph of $C$. 
     
     As a result, $\dim End_{C_A}(X^{2n})$ counts the number of paths of length $2n$ from $(0)$ to $D(\mu)^T$ for $\mu$ a strict partition of even size. By Corollary~\ref{typeCPathCounting}, $\dim End_{C_A}(X^{2n}) = 2^{n - 1}(2n - 1)!!$.
\end{proof}

\subsection{Type D}

The root system $D_r$ for $r > 3$ is realized in Euclidean space with the usual basis $\set{e_i}$ as
\[
\Phi = D_r = \defset{e_i - e_j}{1\le i \ne j\le r} \U \defset{e_i + e_j}{1\le i \ne j\le r}
\]
A convenient choice of positive roots is 
\[
\Phi^+ = \defset{e_i - e_j}{1\le i < j\le r} \U \defset{e_i + e_j}{1\le i < j\le r},
\]
with simple roots
\[
\Delta = \defset{e_i - e_j}{1\le i < j\le r}.
\]

Define the labeling of the shifted Young diagram $\mu = (2r - 2, 2r - 4, \dots)$, $T: \Phi^+ \to \mu$, by
\begin{align*}
    T(e_i - e_j) = (j, i)\\
    T(e_i + e_j) = (2r - 2 - j, i).
\end{align*}

Note that this arrangement differs slightly in shape from the type B and C cases. 

\begin{example}
For $r = 4$, \\
\begin{align*}
    \scalebox{0.5}{
\ytableausetup{boxsize = 4em}$
\begin{ytableau}
e_1 + e_2 & e_1 + e_3 & e_1 + e_4 & e_1 - e_4 & e_1 - e_3 & e_1 - e_2 \\
\none & e_2 + e_3 & e_2 + e_4 & e_2 - e_4 & e_2 - e_3 \\
\none & \none & e_3 + e_4 & e_3 - e_4
\end{ytableau}
$}
\end{align*}
\end{example}

We now state analog of Lemmas~\ref{idealsTypeC}~and~\ref{idealsTypeB}. Note that this proof follows exactly after the patterns of these analogous lemmas, so we only provide the broad strokes, which are routine. Again, let $D = \defset{e_i + e_j}{i < j}\sbs \Phi^+$, and note that $T(D)$ corresponds to the shifted staircase diagram $\delta = (r - 1, \dots, 1)$.

The proof of the following lemma is identical to the type B case, and mostly identical to the type C case.

\begin{lemma}\label{idealsTypeD}
    The ideals $\Psi \in \Sigma$ are in bijection with the shifted subdiagrams $\mu\sbs T(\Phi)$, where the bijection is given by $\Psi \mapsto T(\Psi)$. If $\Psi \sbs D$ is an ideal, then $\Psi$ is an abelian ideal. 
\end{lemma}

Next, we prove a weaker analog of Proposition~\ref{sameAsDoubledC} which is slightly stronger than the analogous results for type B.

\begin{lemma}
    Suppose $\Psi\in \Sigma$ is abelian and $\Psi \sbs D$. Then 
    \[
    \lambda_{\Psi} = D(T(\Psi))
    \]
    Moreover, if $\Psi \in \Sigma$ is any ideal, we have 
    \[
    |\lambda_{\Psi}| = |\lambda_{\Psi \I D}|
    \]
\end{lemma}

\begin{proof}
    The proof here is very similar to the type B and C cases. Given an abelian ideal $\Psi \sbs D$, the roots labeling a square $(i,j)\in T(\Psi)$ all have positive coordinates. Hence the $i$-th row of $\lambda_\Psi$ has length equal to the number of squares in $T(\Psi)$ labeled by a root with a non-zero $e_i$-coordinate. By definition of $T$, any such root is of the form $e_i + e_j$ since $\Psi \sbs D$, for $i < j$, or $e_k + e_i$ for $k < i$. Roots of the first form are found in the $i$-th row of $T(\Psi)$, and roots of the second are found in the $i-1$-th column. Since $T(\Psi)$ is a strict Young diagram, the $i$-th row and $i- 1$th column are disjoint. Hence $(\lambda_{\Psi})_i = T(\Psi)_i + T(\Psi)_i^t = D(T(\Psi))_i$.

    Next, we make a similar argument to the type B and C cases to show the last half of the lemma. The main point is that by adding a square labeled by a root not in $D$ to an ideal $\Psi$ (that is, a root of the form $e_i - e_j$), one removes a square from the $j$-th row, and adds a square to the $i$-th row, thus preserving the size of the Young diagram $\lambda_{\Psi}$.
\end{proof}

\begin{corollary}\label{TypeDZ2Graded}
    $C_A$ is a $\Z/2\Z$ graded category. 
\end{corollary}
\begin{proof}
    Since the simple summands of the Etale algebra $A$ corresponding to abelian ideals $\Psi \in \Sigma$, and since the doubling of any shifted Young diagram is even, and since $|\lambda_\Psi| = |\lambda_{\Psi\I D}|$, it follows that the Etale algebra is an even object in $C$. 
\end{proof}

\begin{remark}
Note that it can easily be shown that for any abelian ideal $\Psi$, the weight $\lambda_{\Psi}$ always corresponds to a Young diagram of even size. This is similar to the type C case, but not to the type B case. Since we will not need to alter the pivotal structure in the type D case, we have no need for this fact, although it is interesting in its own right. 
\end{remark}

In direct analogy to type B, we have the following. 
\begin{lemma}\label{gInEnd2D}
    Let $X = F_A(V_{(1)})$. Then 
    \[
    \dim Hom_{C_A}( X\tensor X \to g) = 1
    \]
  
\end{lemma}
\begin{proof}
The proof is identical to the analogous proposition for type B and C.
\end{proof}

Finally, the following will help us establish limiting fusion rules.

\begin{lemma}\label{stableDimD}
    For all $n \in \N$ there is some $L \in \N$ such that whenever $ r > L$, we have 
    \[
    \dim End_{C_A}(X^{\tensor n})  = 2^{n - 1}(2n - 1)!!
    \]
\end{lemma}
\begin{proof}
    Same as type B.
\end{proof}

\section{The Clifford BMW Supercategories}\label{CWSection}

\begin{definition}
    Define the ring $R = \C[q, \inv{q}, \frac{1}{q - \inv{q}}, \frac{1}{(1 + i)q + (i - 1)\inv{q}}, \frac{1}{(1 + i)\inv{q} + (i - 1)q}]$. 
\end{definition}

\begin{definition}
    The Clifford $BMW$ supercategory, denoted $CW(\sigma, \epsilon)$ is the $R$-linear pivotal supercategory which evenly extends $BMW(q, \sigma i q^{-\epsilon})$ by a single odd morphism $\att{oddDot}{0.5}: 1 \to 1$, subject to three additional relations
$$
\begin{array}{ccc}
   \text{(Rotate) } \left(\att{oddDot}{0.5}\right)^* = \att{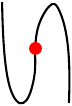}{0.5} = -\sigma i \att{oddDot}{0.5} & \text{(Half Braid) } \att{oddHalfBraid1}{0.5} = \att{oddHalfBraid2}{0.5} & \text{(Cancel) } \att{doubleDot}{0.5} = \att{stick}{0.5}
\end{array}
$$
\end{definition}

Recall that the BMW algebras, $BMW_n(q, \theta)$, are linear spans of $(2n - 1)!!$ diagrams, or equivalently, framed tangles with $n$ boundary points on a handle on top, and $n$ boundary points on a handle below. For any subset $S \sbs \set{1, \dots, n}$, let $v_S = \prod_{i\in S} v_i$ where the product is made over $S$ in increasing order. 

It will be convenient to fix the following terminology. Any instance of the morphism $\att{oddDot}{0.5}$ in a Clifford BMW diagram shall be referred to as a \textbf{spinor}, because it requires a 720 degree rotation for $v$ to return to $v$. If we ignore spinors, we obtain a BMW diagram. Strands in a Clifford BMW diagram refer to the strands of the resulting BMW diagram. 

The relation (Rotate) as well as the super interchange rule and the usual snake equation for caps and cups imply that any subdiagram of a Clifford BMW diagram can be deformed up to para-isotopy, as treated in \cite{LiuAndJaffe}. One can formalize this geometrically in a number of equivalent ways. We will refer to any application of super-interchange along with the relation (Rotate) as \textbf{super-isotopy}.

\begin{theorem}\label{smallSpanningSet}
    Pick some diagram basis $\mathcal{D} = \set{b_1, \dots, b_{(2n - 1)!!}}$ in the image of $BMW_n(q,\sigma iq^{-\epsilon)}$ inside $End_{CW(\sigma, \epsilon)}(n)$, such that none of the diagrams $b_i$ have closed loops. We may choose this basis so that no strand crosses another strand twice, and no strand crosses itself.

    Then the set
    \[
    \mathcal{B} = \set{v_S b_\alpha v_T} 
    \]
    where the set ranges over all subsets $S$ of the leftmost boundary points of the strands in $b_\alpha$, and where $T$ ranges over all of the leftmost boundary points of $b_\alpha$ which do not connect to the top, spans all of $End_{CW(\sigma, \epsilon)}(n)$. 
\end{theorem}

\begin{proof}
    The proof is by induction on the number of crossings. If there are no crossings, then any spinors in a diagram $D\in End_{CW(\sigma, \epsilon)}(n)$ can be freely moved around on its respective strand by para-isotopy, picking up multiples of $i$ and $-1$ for every 180 degree rotation and interchange. As a result, any two spinors on a single strand can be brought close to each other and canceled, up to multiplication by a unit in $R$ (a power of $i$). Noting that every closed loop with an odd number of strands evaluates to zero, we can remove any closed loops in the diagram. Then we can freely move each of the remaining instances of $v$ to their appropriate locations as prescribed by the statement of the theorem, again picking up powers of $i$.

    \begin{remark}
        The base case of our argument is implicit in Liu and Jaffe's argument in \cite{LiuAndJaffe} for para-TLJ, since the crossing free subcategory is a quotient of their $\Z/2\Z$ graded para-TLJ. 
    \end{remark}

    Now we suppose that every diagram $ D\in End_{CW(\sigma, \epsilon)}(n)$ with less than $N$ crossings can be evaluated. 
    
    Suppose first that the diagram in question has at least one closed loop, which we shall call $L$. Then using the Dubrovnik relation, bring $L$ behind every other strand in the diagram. This can be done up to a linear combination of diagrams with less than $N$ crossings, each with a scalar of $\pm (q - \inv{q})$ which is inverted in $R$. To each of these additional diagrams, we may apply the induction hypothesis. Having brought the loop $L$ behind every other strand, we now remove any twists in $L$ as follows. Suppose we have some connected interval $I$ on the loop $L$ which involves only one crossing $c$ of $L$ with itself. By the relation (Half Braid), isotopy, (R3) and (R2) we can ensure that $I$ does not cross any of the other strands in $D$. Using the relation (Half Braid) and (Rotate), move every spinor dot under the crossing $c$, freeing the section $I$ of any spinors. Then, applying the relation (R1), we have $D = \theta D'$ where $D'$ has at least one less crossing than $D$. Note that $\theta \in R^\times$. Now apply the induction hypothesis to $D'$. An example of the above process is shown below:\footnote{In diagrammatics for a supercategory, odd morphisms must be in general position with respect to their vertical coordinates. The left dot in the first diagram is below the right one.}
    \[
    \att{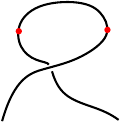}{0.5}  = \att{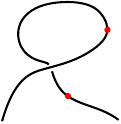}{0.5}  = -\sigma i\att{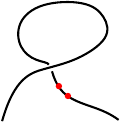}{0.5} = -\sigma i \theta \att{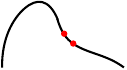}{0.5}
    \]
    In the case where the closed loop $L$ does not cross itself, use para-isotopy, (R3), (R2), and (Half Braid) to shrink $L$ down into a small neighborhood. $L$ then evaluates to either $\delta$ or $0$ depending on the parity of the number of spinors on it. 
    
    Now suppose there are no closed loops in the diagram. By the relation (Half Braid), any spinor can slide under any crossing. If, when moving the spinors to their appropriate locations, we wish to bring one over a crossing, we may apply the Dubrovnik relation to exchange the over-crossing for an under-crossing. This accumulates two terms with one less crossing, to which we may apply the induction hypothesis. Note that the Dubrovnik relation accumulates terms with leading coefficients as $\pm (q - \inv{q})$, which is inverted in $R$.

    Since each of the accumulated scalars at each step of the argument was a unit in $R$, it follows that every element in $End_{CW(\sigma, \epsilon)}(n)$ is an $R$-linear combination of elements in $\mathcal{B}$.
\end{proof}

If we omit the data of which strand passes over which at each crossing, letting $D \sim D'$ mean that $D = D' + \sum \zeta_\alpha D_\alpha$ with each $D_\alpha$ having less crossings than $D$, and with $D'$ having no more crossings, then the argument can be visualized by literally sliding the spinors to their appropriate locations and untangling the diagram, modulo diagrams with less crossings. For example, 
\[
\att{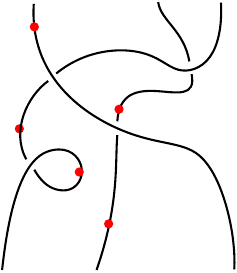}{0.5} \sim \att{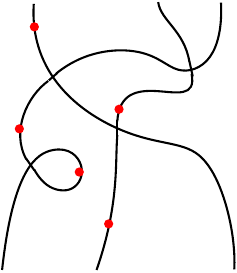}{0.5} \sim \att{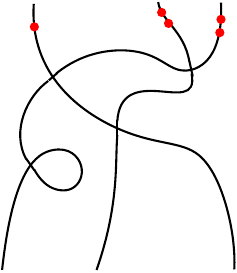}{0.5} \sim \att{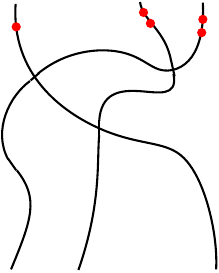}{0.5}
\]

Up to a scaling by powers of $i\in \C$, note that the spanning set above has an even more intuitive description. Taking a basis of BMW diagrams, we have a choice of whether to include an instance of $v$ on each of the strands. This makes the size of the spanning set obvious: there are $2^n$ choices for whether to add a spinor to each of the $n$ strands, and there are $(2n - 1)!!$ possible pairings of the boundary points.

The following observation will also make some of our proofs easier. 

\begin{corollary}\label{standardDiagrams}
    If $D\in Hom_{CW(\sigma, \epsilon)}(2n \to 0)$ is any diagram, then $D = \sum_\alpha \zeta_\alpha D_\alpha$ such that each $D_\alpha$ has at most as many crossings as $D$, and such that for each $\alpha$,
    \begin{enumerate}
        \item No strand $s$ in $D_\alpha$ crosses another strand $t$ twice. In particular, no strand $s$ crosses itself. 
        \item Each strand $s$ in $D_\alpha$ is decorated with at most one spinor, each occuring at the boundary. 
    \end{enumerate}
    For ease of reference in the next subsection, we say that a diagram is \textbf{standard} if the above conditions are met. 
\end{corollary}
\begin{proof}
    This follows directly the proof of Theorem~\ref{smallSpanningSet}, which was done explicitly by induction on the number of crossings. Having brought all of the spinors to the boundary, so that $D = c_\alpha \circ D'_\alpha$ where $D'_\alpha$ is a BMW diagram and $c_\alpha$ is a Clifford diagram. We can apply the Dubrovnik relation, R1, and isotopy, to re-arrange the BMW diagrams at the cost of a linear combination of terms with less crossings. 
\end{proof}

\begin{remark}
    For the analogous Hecke Clifford categories in \cite{CainAndSnyder}, the evaluation algorithm was relatively easy while working entirely in the even category. Due to the non-orientedness of the strands in the Clifford BMW category, none of their arguments generalize. It is still possible, however, to construct an explicit spanning set of size $2^{n - 1}(2n - 1)!!$ for each of the endomorphism algebras while only working in the even category. Before realizing that we could condense the fermion object $g$ before specialization, we originally had an incredibly complicated proof that a certain set of diagrams spans. The degree to which fermion condensation simplifies the proof is, in our opinion, incredible. 
\end{remark}

\subsection{The ring $R$, conjugation, and specialization}

The category $CW(\sigma, \epsilon)$ is defined over the ring $R$ for each choice of $\epsilon$ and $\sigma$. 

Although the categories $CW(\sigma, \epsilon)$ may be non-equivalent, these categories are closely related by automorphisms of the ring $R$. In fact, the ring $R$ was specifically designed for these symmetries to exist.

For any automorphism $\kappa\in Aut(R)$, one can define a new planar algebra $\kappa.CW(\sigma, \epsilon)$ where every relation in $\kappa.CW(\sigma, \epsilon)$ is rewritten according to the automorphism $\kappa$. 

A slightly more general way of viewing this is the following.

\begin{lemma}
    Let $S$ and $T$ be commutative rings. Let $f: T \to S$ be any ring homomorphism. Then $S$ has the structure of a right $T$-module with action for all $s\in S$ and $t\in T$ given by 
    \[
    s.t = sf(t).
    \]
    
    If $S$ and $T$ do not have characteristic $2$,\footnote{This is only needed for the super case} then for any strictly (super) pivotal and strictly monoidal $T$-linear monoidal (super) category $\fanC$, there is another strictly (super) pivotal and strictly (super) monoidal $S$-linear category $f.\fanC$ with the same objects, and with hom spaces 
    \[
    Hom_{f.\fanC}(x \to y) = S\tensor_T Hom_\fanC(x \to y)
    \]
    for all objects $x$ and $y$.
\end{lemma}
\begin{proof}
    This is purely formal. Since $S$ is a left $S$-module by multiplication, and a right $T$-module by the action given above, the tensor product $S\tensor_T Hom_\fanC(x \to y)$ is a left $S$-module. Composition of simple tensor morphisms $s_1\tensor \alpha: x\to y$ and $s_2 \tensor \beta: y \to z$ is given by 
    \[
    (s_2 \tensor \beta)\circ (s_1\tensor \alpha) = (s_1 s_2) \tensor (\beta \circ \alpha)
    \]
    and extended $S$-linearly.
    Since multiplication in $S$ is associative, composition is associative. It is easy to check that this is $S$-bilinear. 

    One defines the monoidal product similarly (although with admittedly overloaded tensor notation), on simple tensors as 
    \[
    (s \tensor \alpha) \tensor_{f.\fanC} (s' \tensor \beta) = (s' s)\tensor (\alpha \tensor_{\fanC} \beta). 
    \]
    and checks that the axioms are satisfied. 
\end{proof}

\begin{remark}
    If the homomorphism $f$ is an isomorphism, then $f$ gives a pivotal equivalence of categories, but not a linear one.
\end{remark}

The above notion generalizes both the idea of Galois conjugates of linear categories defined over a field, as well as specialization for the various continuous families of tensor categories defined over localizations of polynomial rings on some set of variables. For example, the Deligne interpolation categories $Rep(GL_t), Rep(O_t), Rep(S_t)$ are defined over the ring $\C[t, \inv{t}]$. For any $q \ne 0$, there is a specialization map $ev_q:\C[t, \inv{t}] \to \C$, giving a functor from each of these categories for a generic parameter to a specialization. 

There are essentially two distinct types of situations in which freeness of the hom spaces in a category over a ring $T$ translates to freeness of a new ring $S$ under the action of a homomorphism $T \to S$. We will have occasion to use both in this paper. Unlike the above lemma, the following criteria are purely module theoretic, and have nothing to do with monoidal categories. 

\begin{lemma}\label{specializationAndConjugation}
    Let $T$, $S$ $f$ be as in the previous lemma. Suppose a module $V$ is free on a set $\fanB \sbs V$ over $T$. Then $S \tensor_T V$ is free on $1_S\tensor \fanB = \defset{1_S\tensor b}{b\in \fanB}$ if either of the above conditions hold
    \begin{enumerate}
        \item $f$ is an isomorphism
        \item $S\sbs T$ is a subring and $f$ is a retraction onto $S$, (i.e. $f(s) = s$ for all $s\in S$). 
    \end{enumerate}
\end{lemma}
\begin{proof}
    The proofs of these criteria are immediate from the definitions, especially the first one. For the second, if $\sum_\alpha s_\alpha (1_S \tensor b) = 0$, then since $s_\alpha \in T$ for all $\alpha$, and since $ \sum_\alpha s_\alpha (1_S \tensor b) = 1_S \tensor \left(\sum_\alpha s_\alpha b_\alpha\right) = 0$ implies that $\sum_\alpha s_\alpha b_\alpha = 0$, we must have $s_\alpha = 0$ for each $\alpha$ since $V$ was free as a $T$-module on the set $\fanB$. 
\end{proof}

This is why interpolations over infinitely many representations work to prove existence of a one parameter family of categories at all values, although interpolation is ironically unnecessary for proving the existence of each of Deligne's interpolation categories. 

Now we outline the relevant ring homomorphisms for the categories $CW(\sigma, \epsilon)$. 

First, there is clearly the automorphism $\phi\in Aut(R)$ which maps $i \mapsto -i$. Note that for any value of $\sigma\in \set{\pm 1}$, $\phi$ maps $\theta = \sigma iq^{-\epsilon} \mapsto -\theta$.

There is another automorphism $\psi\in Aut(R)$ which makes the substitution $q \mapsto -1/q$ and which fixes $\C \sbs R$. Note that $\C[q]$ embeds into $R$ in at least two ways, one by $q\mapsto q$ and another by $q \mapsto -1/q$. Let $\sms{\psi}$ be the latter embedding. Note that $\sms{\psi}$ fixes $\C$. 

Next, observe that inverting all elements in the multiplicative closure of \[S = \set{q, q - \inv{q}, (1 + i)q + (i - 1)\inv{q}, (1 + i)\inv{q} + ( i - 1)q}\] is the same as inverting all polynomials in the multiplicative closure of the set \[S' = \set{q, q^2 - 1, q^2 + i, q^2 - i}.\] Hence the localization of $\C[q]$ by $q$ followed by the localization of $\C[q,1/q]$ by $S$ is just the localization of $\C[q]$ by $S'$. By the universal property of localizations, a ring homomorphism $\sms{\psi}:\C[q] \to R$ induces a unique ring homomoprhism $R \to R$ if $\sms{\psi}(S') \sbs R^\times$, such that ${\psi}|_{\C[q]} = \sms{\psi}$. 
Substituting $q$ for $-1/q$ gives 
\[
q^2 - 1 \mapsto -\frac{q^2 - 1}{q^2} \in R^\times
\]
\[
q^2 + i \mapsto i \frac{q^2 - i}{q^2} \in R^\times
\]
\[
q^2 - i \mapsto -i \frac{q^2 + i}{q^2}\in R^\times
\]

This establishes the existence of a map $\psi: R \to R$ which maps $q \mapsto -1/q$. Moreover, by our description of $\sms{\psi}$, it follows that for any polynomial $p(q)$, $\psi(p(q)) = p(-1/q)$. Hence for any rational function $p(q)/r(q) \in R$, we have $\psi(p(q)/r(q)) = p(-1/q)/r(-1/q)$. From this it is clear that $\psi^2$ is the identity, and hence that $\psi$ is an automorphism of $R$. It is also clear that $\psi$ fixes constants, and in particular fixes $i$.

The automorphisms $\phi, \psi$ generate a group isomorphic to $\Z/2\Z \times \Z/2\Z$. It is easy to verify that $\phi.CW(\sigma, \epsilon) = CW(-\sigma, \epsilon)$ and $\psi.CW(\sigma, \epsilon) = CW(\sigma, -\epsilon)$. This shows that all of the Clifford BMW categories are conjugate under automorphisms of the ring $R$.

We summarize the findings of this section in the following lemma. 
\begin{lemma}\label{specializationForCW}
    There exist automorphisms $\phi, \psi \in Aut(R)$ such that $\phi.CW(\sigma, \epsilon) = CW(-\sigma, \epsilon)$ and $\psi.CW(\sigma, \epsilon) = CW(\sigma, -\epsilon)$. Moreover, for all $q\in \C$ such that $q^8 \not\in \set{0, 1}$, there exists specialization maps $ev_q: R \to \C$. 
\end{lemma}

\[\begin{tikzcd}
	{CW(1, 1)} & {CW(-1, 1)} \\
	{CW(1, -1)} & {CW(-1, -1)}
	\arrow["\phi", squiggly, from=1-1, to=1-2]
	\arrow["\psi"', squiggly, from=1-1, to=2-1]
	\arrow["\psi", squiggly, from=1-2, to=2-2]
	\arrow["\phi"', squiggly, from=2-1, to=2-2]
\end{tikzcd}\]

\subsection{Evaluability of the Special Clifford BMW category}\label{SCW}.

We now define certain diagrammatic extensions of the Clifford BMW categories, which we will use to present the quantum subgroups treated in this paper of types B and D.

\begin{definition}
    Let $q^{8} \ne 1, 0$ and $N > 2 \in \N$ and $\theta = iq$, and $\epsilon\in \set{\pm 1}$, such that the idempotents $P_\lambda$ in $BMW(q, \theta)$ are defined when $\lambda^t_1 +\lambda^t_2 \le N$. Define the \textbf{special Clifford BMW category}, denoted $SCW_N(q)$, to be $\C$-linear pivotal category with symmetrically self-dual objects $\N$, which extends $CW_q(1, -1)$ by the even generator 
    \[
    \att{extraGenerator}{0.5}\in Hom_{SCW_N(q)}(N \to 0)
    \]
    subject to the relations
    \[
    \text{(Braid)} \att{extraGeneratorBraided}{0.5} = \att{extraGeneratorUnbraided}{0.5} \hspace{0.5cm} \text{(Split)} \att{extraGeneratorProjection}{0.5} = P_{(1^N)}
    \]
    where $P_{(1^N)}$ is the idempotent in $BMW(q, \theta)$ as in section~\ref{BMWIntro}, as well as the relation which we impose only when $N= 4k$ for some integer $k$,
    \[
    \text{(Pedals) } \att{dotsOnPedels}{1} = \lambda
    \]
    where $\lambda = \sqrt{\att{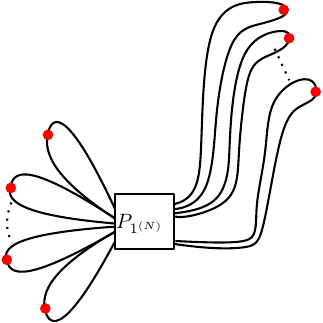}{1}}$.

    We also require  the additional relation (Sphericality) for every endomorphism diagram $f: m \to m$,
\[
\text{ (Sphericality) } tr_L(f) = tr_R(f)
\]
where $tr_L$ denotes the left trace defined with caps and cups, and wehre $tr_R$ denotes the right trace, again defined with caps and cups.
\end{definition}

Note that by super-isotopy and the relation (Split), one can show that, if there is any scalar $\lambda \in \C$ such that $\att{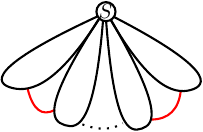}{1} = \lambda$ such that the resulting pivotal super-category is not zero, the one must have $\lambda = \pm \sqrt{\att{stupidNumber}{0.5}}$. This relation looks much more natural if one works in the even subcategory. Nonetheless, everything other than this relation is more natural in the full supercategory.

\begin{remark}
    We suspect that $\lambda = \pm 1$, but we currently do not know how to prove it. 
\end{remark}

Let $SCW'_N(q)$ be the category generated by the same morphisms as $SCW_N(q)$, with the same relations, only that $\lambda$ is replaced with $-\lambda$. Then sending all of the Clifford BMW generators in $SCW_N(q)$ to the corresponding generators in $SCW_N'(q)$, while also sending $S \mapsto -S$, gives an equivalence of categories, where one does the exact same thing to produce an inverse.

Note that there are many values of $q$ for which the Special Clifford BMW category can be defined. This does not mean, however, that the Special Clifford BMW categories are non-zero at these values.  

As stated in the introduction, this category will serve as a diagrammatic presentation for the categories $C(\so_{N}, \check{h})_A$ where $A$ is the Etale algebra coming from the conformal embedding $\V(\so_N, \check{h}) \sbs \V(\so_{\dim \so_N},1)$. In order for this category to be useful to us, we need to show that it is evaluable, that is, that every endomorphism of the unit object is a scalar multiple of the identity.

    \begin{lemma}~\label{EvenNumberOfS}
        If $D\in Hom_{SCW_{N}(q)}(m \to n)$ is a diagram such that $D$ has an even number of $N$-valent vertices $S$, then $D$ is a linear combination of Clifford BMW diagrams. 
    \end{lemma}
    \begin{proof}
        
    Any two $N$-valent vertices, they can be brought close together by repeated applications of the relation (Braid). Upon applying the relation (Split), this diagram becomes a linear combination of Clifford BMW diagrams. As a result, we simply apply the Clifford BMW category's evaluation algorithm in this case (Theorem~\ref{smallSpanningSet})
    \end{proof}

\begin{theorem}\label{evaluationAlgorithmSpecialCliffBMW}
    For any admissible choice of $q^{8}\ne 1, 0$ and $N\in \N$ such that the idempotent $P_\lambda$ in $BMW(q, iq)$ is defined whenever $\lambda_1^t + \lambda_2^t < N$,\footnote{That is, $q$ and $N$ such that the Special Clifford BMW category is defined} every element $f\in End_{SCW_N(q)}(0)$ is a scalar multiple of the identity morphism. 
\end{theorem}

\begin{proof}
    First, let us consider the case where $N$ is odd. In that case, any diagram $D\in End_{SCW_N(q)}(0)$, upon deleting the red strands representing projections (while keeping their attaching points) is a graph consisting of $4$-valent vertices representing crossings, as well as $N$-valent vertices. Since the total degree is even, it follows that there must be an even number of $N$-valent vertices in this graph.

    On the other hand, if $N$ is even, we must be a little bit more careful. Let $D\in End_{SCW_N(q, \lambda)}(0)$. In this case, for any pair of $N$-valent vertices, we can apply the same process as above, until we are left with a linear combination of diagrams with at most one $N$-valent vertex. 

    Note that the observation that a diagram $D$ in $Hom_{SCW_{N}(q)}$

    Hence we need only show how to evaluate diagrams of the form \[S\circ D = \att{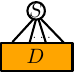}{1}\] for $D: 0 \to 2r$ a Clifford BMW diagram.
    
    \begin{lemma}
        Let $D: 0\to 2r$ be a Clifford BMW diagram. Then $S\circ D$ is a linear combination of diagrams with no crossings and one $2r$-valent vertex. 
    \end{lemma}

    \begin{proof}
      The proof is by induction on the number of crossings in $D$. By Corollary~\ref{standardDiagrams} we may assume that $D$ is standard. Let $c$ be any crossing in $D$. We may choose $c$ such that there is a path from $c$ to the boundary on one of the strands involving $c$, which does not meet any other crossings. Now select two paths $\gamma_1$ and $\gamma_2$ from $c$ to the boundary on distinct strands $s$ and $t$ respectively (the two strands involved at the crossing $c$ must be distinct by the assumption that $D$ is standard). Assume without loss of generality that $s$ crosses under $t$ at $c$, because otherwise we can apply the Dubrovnik relation to exchange the over-crossing for an under crossing at the cost of two terms with less crossings, to which we may apply the induction hypothesis. 
      
      Apply the relation Dubrovnik repeatedly, creating terms involving diagrams with less crossings, such that the paths $\gamma_1$ and $\gamma_2$ cross under at any crossings when traveling from the boundary to $c$. If the path $\gamma_1$ has any spinors on it, we may apply the relation (Half Braid) to bring $v$ under every crossing on $\gamma_1$, and then under at $c$, thus freeing $\gamma_1$ from any spinors. Now apply the Dubrovnik relation again at $c$, so that $t$ is below $s$ at $c$. If $\gamma_2$ has any spinor decorations, repeatedly apply the relation (Half Braid) to bring this decoration under $c$, thus freeing $\gamma_2$ from all spinor decorations. 

      Hence $S\circ D'$ can be written in the form 
        \[
        D'\circ S = \att{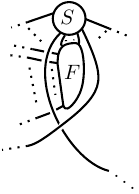}{1}
        \]
        for some sub-diagram diagram $F$. Using the fact that $S$ absorbs braids at the cost of scalars, we perform the following diagrammatic manipulations which hold up to scalar multiples,
        \[
        \att{handleCrossing1}{1} = \kappa_1\att{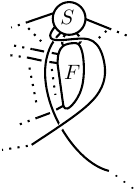}{1} = \kappa_1\att{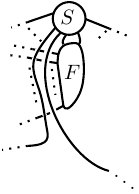}{1} = \kappa_2\att{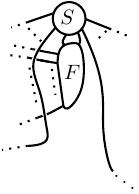}{1} = D'''
        \]
        where the scalars $\kappa_1$ and $\kappa_2$ come are powers of $-1/q$. Call the resulting Clifford BMW diagram $D''$, so that $S\circ D'  = \kappa S\circ D''$. 

        The resulting diagram $D''$ has exactly one less crossing than $D'$, so we may apply the induction hypothesis. 
    \end{proof}

    \begin{example}
        Pretty much the entire idea behind the above lemma is captured in the following example. 
        \begin{align*}
        \att{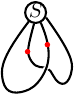}{1} = \att{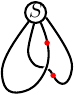}{1} = \att{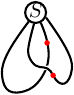}{1} + (q - \inv{q})\left(\att{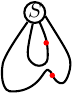}{1} - \att{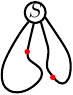}{1}\right)\\
        = -\att{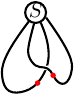}{1} + (q - \inv{q})\left(\att{crossingElimination3}{1} - \att{crossingElimination4}{1}\right) = \frac{1}{q}\att{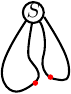}{1} + (q - \inv{q})\left(\att{crossingElimination3}{1} - \att{crossingElimination4}{1}\right)
        \end{align*} 
    \end{example}

    We will use the following auxiliary definition, borrowed from \cite{OkadaDiagrams}

    \begin{definition}
        Let $D: 0 \to 2n$ be a Clifford BMW diagram with no crossings. A strand $s$ is said to \textbf{contain} a strand $t$ if the region bounded by $s$ and the top boundary contains $t$. The \textbf{depth} of such a diagram is the maximal number of strands contained in another strand. The \textbf{depth} of an individual strand $s$ is the number of strands containing $s$. 
    \end{definition}

    \begin{example}
        Note that the actual combinatorial information concerning containments is not affected at all by spinor decorations. 

        In the diagram below, red arrows are drawn to indicate which strand contains which. 
        \[
        \att{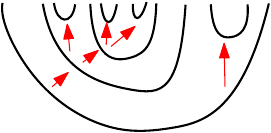}{0.5}
        \]
    \end{example}

    \begin{remark}
        This notion of containment plays an important role in giving a diagrammatic description of the Okada algebras in \cite{OkadaDiagrams}.
    \end{remark}

    \begin{lemma}
        Let $D$ be Clifford BMW diagram with no crossings. Then $S\circ D = \kappa S \circ D'$ for some $\kappa$ in $\C$ and for some Clifford BMW diagram $D'$ with no crossings, such that each strand has at most one spinor decoration, and such that $D'$ has depth zero.
        
        In other words, $S \circ D$ is of the form
        \[
        \att{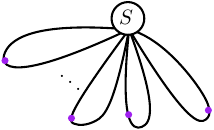}{1}
        \]
        where the purple dots are used to indicate where a spinor decoration may be placed. 
    \end{lemma}
    \begin{proof}
    Using super-isotopy, and (Cancel), we may bring any two spinor dots on a single strand close, and cancel them. As a result, we may assume that $D$ has only one spinor dot per strand. 

    Now we induct on the depth of a diagram $D$ paired with the number of strands of maximal depth, in lexicographic order. If the depth of a diagram $D$ is zero, we are done. If $s$ is a strand of maximal depth with a strand $t$ immediately containing $s$, then we may assume that the left boundary point $x_s$ of $s$ is adjacent to the left boundary point $x_t$ of $t$. Otherwise, some other strand $r$ is between the $x_s$ and $x_t$ with left boundary point immediately adjacent to $x_t$. Since $r$ is also contained in $t$, it follows that $r$ also has maximal depth, and hence we may simply consider $r$ instead of $t$. 

    Using the relations (Half Braid), (R2), and the fact that $S$ absorbs braids at the cost of a scalar, we can then calculate

    \begin{align*}
    S \circ D = \att{DepthReduction1}{1} \\= \att{DepthReduction2}{1} = \kappa_1 \att{DepthReduction3}{1} \\= \kappa_1\kappa_2\att{DepthReduction4}{1} = \kappa_1\kappa_2 S \circ D'
    \end{align*}
    where $\kappa_1$ and $\kappa_2$ are the scalars accumulated from $S$ absorbing a braid.

    The resulting diagram $D'$ either has one less strand of maximal depth, or the depth of $D'$ is less than that of $D$.
\end{proof}

We are now ready to finish the proof of Theorem~\ref{evaluationAlgorithmSpecialCliffBMW} in the case where $N = 2r$. By the previous lemma, the proof reduces to showing how to evaluate a diagram of the form 
    \[
       S \circ D = \att{pedalForm}{1}
        \]
        as a scalar multiple of the identity. 

There are two cases. Either each strand is decorated with exactly one spinor, or at least one strand is missing a spinor. 

In the former case, we either have $ 2r = 4k $ for some integer $k$, or we don't. If we don't, then $D$ has an odd number of dots, contradicting the fact that $D$ is in the Clifford BMW category. Otherwise, we may directly apply the relation the relation (Pedal), evaluating $S \circ D$ as $\lambda = tr\left(\att{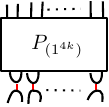}{0.25}\right)$. 

On the other hand, if on of the strands in $D$ is missing a dot, we may apply the relation (Killed By Cup), so that $S\circ D = 0$. 
\end{proof}

\section{Connecting the Clifford BMW categories and the module categories}

In this section, we will finally tie together the results of the previous two sections. In section~\ref{ConformalEmbeddingsSection}, we established dimension formulas for endomorphism algebras in $C_A$, which stabalized at $2^{n - 1}(2n - 1)!!$ as $ r\to \infty$. In section~\ref{CWSection}, we established an upper bound of $2^{n -1}(2n - 1)!!$ for the dimensions of endomorphism algebras in the CW categories. It turns out that this is no coincidence!!\footnote{double exclamation marks for double factorials!!}

In this section, we will construct quotients of the (S)CW algebras inside the categories of modules $C_A$ outlined in section~\ref{ConformalEmbeddingsSection}. We will show that, in each case, each of these subcategories is full and dominant. As a result, the (S)CW skein theory is, in principle, enough to fully understand the categories of modules $C_A$. This in turn will establish two continuous families of CW categories for which the upper bound of $2^{n -1}(2n - 1)!!$ for the dimensions of endomorphism algebras is also a lower bound. 

For each case, we will use the following matrix. Rather than re-defining this matrix in each section, it will be more efficient to introduce it first here.

Fix some $q\in \C$ such that $q^8 \ne 1, 0$, and $\epsilon \in \set{\pm 1}$. Consider the arrows 
\[
f_1= \att{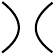}{0.5} ,f_2 =\att{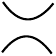}{0.5}, f_3=\att{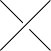}{0.5} ,f_4 =\sigma i \delta^{-2}\att{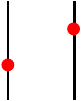}{0.25}, f_5 = \sigma i \delta^{-2}\att{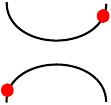}{0.25}, f_6= f_3 \circ f_4
\]
Then the Gramm matrix $A$ of span of these endomorphisms for the inner product $\inner{f_i}{f_j} = tr(f_i \circ f_j)$ is
\[A(\sigma, \epsilon)= \left[\begin{matrix}\delta^{2} & \delta & \delta \theta & 0 & 1 & \frac{1}{\beta}\\\delta & \delta^{2} & \frac{\delta}{\theta} & 1 & 0 & \frac{1}{\theta}\\\delta \theta & \frac{\delta}{\theta} & \frac{\delta \left(\delta q \theta - q^{2} + \theta^{2} \left(q^{2} - 1\right) + 1\right)}{q \theta} & \frac{1}{\beta} & \beta & - \frac{q}{\theta} + \frac{1}{q \theta} + \frac{q}{\beta} - \frac{1}{\beta q}\\0 & 1 & \frac{1}{\beta} & 1 & \frac{\beta^{2} \theta^{2}}{\delta} & \frac{\theta}{\delta}\\1 & 0 & \beta & \frac{\beta^{2} \theta^{2}}{\delta} & 1 & \frac{1}{\beta \delta \theta^{2}}\\\frac{1}{\beta} & \frac{1}{\theta} & - \frac{q}{\theta} + \frac{1}{q \theta} + \frac{q}{\beta} - \frac{1}{\beta q} & \frac{\theta}{\delta} & \frac{1}{\beta \delta \theta^{2}} & \frac{- q^{2} + \theta^{2} \left(q^{2} - 1\right) + 1}{\beta^{2} \delta q \theta^{3}}\end{matrix}\right]\]
where, as usual, $\delta = \frac{\theta - \inv{\theta}}{q - \inv{q}} + 1$, $\theta = \sigma iq^{-\epsilon}$, and $\beta = \epsilon q^\epsilon$, and where the entries are understood as scalar multiples of the identity morphism on the unit object.

Computing determinants we have 
\[
\det A(1, -1) = - \frac{8 i \left(q - i\right) \left(q + i\right) \left(q^{2} - i\right) \left(q^{4} + i\right)^{3}}{q^{4} \left(q - 1\right)^{4} \left(q + 1\right)^{4}}
\]
and
\[
\det A(-1, 1) = \frac{8 i \left(q - i\right) \left(q + i\right) \left(q^{2} + i\right) \left(q^{4} - i\right)^{3}}{q^{4} \left(q - 1\right)^{4} \left(q + 1\right)^{4}}
\]
as scalar multiples of the identity morphism on the unit object.

This allows us to prove the following lemma. 
\begin{lemma}\label{constrantsOnHom2to2}
    For any $q\in \C$ such that $q^8 \ne 0, 1$, \[\dim End_{\sms{CW_q(1, -1)}}(2), \dim End_{\sms{CW_q(-1, 1)}(2)} \in \set{0, 6}\] so long as $q^{16} \ne 1$. If $End_{CW_q(\sigma, \epsilon)}(0) \ne 0$, then $\dim \dim End_{\sms{CW_q(-1, 1)}(2)} = 6$. Moreover, for $N >2$ and for any $q$, either $\dim End_{\sms{SCW_N(q)}}(2) \ge 6$ or $End_{SCW_N(q)}(0) = 0$, so long as $q^{16} \ne 1$.
\end{lemma}
\begin{proof}
If $\sum_i \zeta_i f_i = 0$ for $\zeta_i \ne 0$ for some $i$, then by linearity of compositions and traces, we have $\sum_i \zeta_i A_{ij} = 0$ for all $j = 1, \dots, 6$. This gives a linear dependence of the columns in $A(\sigma, \epsilon)$. So long as $\dim End_{CW_q(\epsilon, \sigma)}(0) > 0$, the above calculations show that $A$ is non-degenerate, since its determinant has roots only at forbidden values of $q$, or sixteenth roots of unity. 

All of this transfers over to the semisimplification as well, because the semisimplification functor is pivotal. So long as $SCW_N(q)$ is not the zero category, the same argument establishes linear independence of the morphisms $f_1, \dots, f_6$, both in $End_{SCW_N(q)}(2)$ and in $End_{\sms{SCW_N(q)}}(2)$.
\end{proof}

\begin{remark}
    We will show in this section that in all relevant cases, $End_{SCW_N(q)}(2)$ is at least six dimensional, and that $End_{\sms{SCW_N(q)}}(2)$ is exactly six dimensional. Nonetheless, it is not at all clear that $End_{SCW_N(q)}(2)$ is six dimensional. The additional $N$-valent vertex does seem to allow for diagrams $2 \to 2$ which are not evaluable in terms of Clifford BMW diagrams. These additional morphisms also turn out to be negligible.
\end{remark}

\begin{lemma}[Twist Lemma]
    Suppose there is a super pivotal extension of the category $BMW(q, \theta)$ by an extra odd generator $V = \att{oddDot}{0.5}: 1 \to 1$ satisfying 
    \[
    V^* = \lambda V
    \]
    for some $\lambda \in \C$, and such that $V$ satisfies the relation (Half Braid) of the Clifford BMW supercategories. Suppose that these relations do not collapse so that $V = 0$. 
    Then $\lambda = -\sigma i$
    where
    \[
    \sigma = i\frac{\inv{\theta} + \theta}{q - \inv{q}} .
    \]
    In particular, $\theta$ is determined up to two possibilities by $\sigma$, and $\sigma$ is entirely determined by $\theta$. There is some $\epsilon \in \set{\pm 1}$ such that $q = \sigma i q^{-\epsilon}$. 
\end{lemma}
\begin{proof}
\begin{align*}
\theta^{-1} \att{oddDot}{1.5} = \att{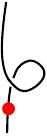}{1} =-\sigma i\att{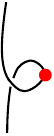}{1} \\
= -\sigma i\att{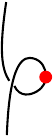}{1} + \sigma i(q - \inv{q}) \left( \att{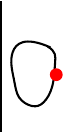}{1}  - \att{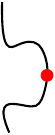}{1}\right)= -\att{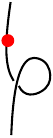}{1} - \sigma i ( q- \inv{q}) \att{oddDot}{1.5} \\
= (-\theta - \sigma i(q - \inv{q})) \att{oddDot}{1.5}.
\end{align*}

Since $V: 1 \to 1$ isn't zero, we have 
\[
\sigma = i\frac{\inv{\theta} + \theta}{q - \inv{q}} .
\]
Solving for $\theta$, we then have $\theta = \sigma i q^{\pm 1}$. 
\end{proof}

This lemma will be useful for defining functors from the Clifford BMW supercategory.

\subsection{$C_A$ for type C and Deligne interpolation}. 

The goal of this section is to prove the following theorem.

\begin{theorem}\label{bigResultC}
    For $r\ge 1$ with $q = e^{\frac{2\pi i}{8r + 8}}$ there exists a full and dominant even superpivotal functor 
    \[
    \Phi: CW_q(-1, 1) \to \mathcal{C}(\mathfrak{sp}_{2r}, r + 1)_A)/\gamma
    \]
    which descends to a full and faithful dominant functor 
    \[
    \sms{\Phi} : \sms{CW_q(-1, 1)[0]} \to \mathcal{C}(\mathfrak{sp}_{2r}, r + 1)_A
    \]
    thus giving an even superpivotal equivalence 
    \[
    Ab(\sms{CW_q(-1, 1)}) \simeq (\mathcal{C}(\mathfrak{sp}_{2r}, r + 1)_A)/\gamma
    \]
    as well as a pivotal equivalence
    \[
    Ab(\sms{CW_q(-1, 1)}) \simeq \mathcal{C}(\mathfrak{sp}_{2r}, r + 1)_A)
    \]
    where $\gamma$ is the half braid on the fermion object $g\in C_A^0$.
\end{theorem}

\begin{proof}

Let $r \ge 1$ be fixed, and let $q$ be as above. Throughout the remainder of this section, we shall let $C$ denote a strictly-pivotal model of the category $\mathcal{C}(\mathfrak{sp}_{2r}, r + 1)$, with the appropriate pivotal structure. The details of \textit{which} pivotal structure will be provided later. Similarly, let $C_A$ denote the category of $A$-modules over $C$. 

\subsubsection{Defining the functor}

The pivotal structure on $C$ which makes the defining representation $V_{(1)}$ symmetrically self dual is not the pivotal structure to which Proposition~\ref{freeModuleFunctorPivotal} applies. We need to verify that, with this change in pivotal structure, the functor $F_A: C \to C_A$ is still pivotal. As noted in the preliminary section, this is equivalent to the condition that $\theta_A = 1$ where $\theta$ is the twist induced from the braid and the pivotal structure. This is where our the observation from subsection~\ref{ConformalEmbeddingsSection} that every simple subobject of $A$ is indexed by a Young diagram of even size. 

\begin{proposition}
    Let $G: BMW(q, -i/q)\to C$ be the functor which, after taking the semisimplification and Cauchy completion, gives an equivalence $G: Ab(\sms{BMW(q, -i/q))}\to C$. When $C$ is equipped with the unique pivotal structure such that $Ab(\sms{G})$ is a pivotal functor, the free module functor $F_A : C \to C_A$ remains pivotal. 
    
\end{proposition}

\begin{proof}

    The category $C$ is $\Z/2\Z$ graded by the ring\footnote{$K(C)$ denotes the Grothendiek ring of $C$, which is a $\Z^+$ based ring with the above $\Z/2\Z$ grading} homomorphism $\deg: K(C) \to \Z[\set{\pm 1}]\cong\Z[\Z/2\Z]$, given by 
    \[
    \deg(V_\lambda) = (-1)^{|\lambda|}.
    \]

    Let $\phi$ be the unique pivotal structure on $C$ such that the objects of $C$ have positive quantum dimensions. As shown by Copeland in \cite{DanielDude}, there are exactly two spherical structures on a $\Z/2\Z$ graded category which are compatible with the grading. Explicitly, these are $\phi$ and $\phi'$ where 
    , \[
    \phi'_{V_\lambda} = \deg(V_{\lambda})\phi_{V_{\lambda}} = (-1)^{|\lambda|}\phi_{V_\lambda} \]

    In order for the equivalence $G: BMW(q, -i/q) \to C$ to be pivotal, $C$ must be equipped with the pivotal structure $\phi'$, as can be seen by observing that the trace of the single strand in $BMW(q, -i/q)$ is negative. 

    By Lemma~\ref{idealsTypeC}, for every ideal $\Psi\sbs C_r$, $|\lambda_\Psi| \equiv 0 \pmod{2}$. In particular, for any abelian ideal $\Psi$ of even size, one has $ |\lambda_\Psi| \equiv 0 \pmod{2}$. Hence every summand in $A$ must appear in an even tensor power of $V_{(1)}$. As a result, $A$ has a well-defined grading, with $\deg(A) = 1$. Hence $\phi'_A = \phi_A$. 

    By Theorem 2.9 of \cite{qMacky}, with ribbon structure $\theta$ induced from the pivotal structure $\phi$ and the braiding generated by the braiding on the BMW category, we have $\theta_A = I_A$. The change in pivotal structure gives a different twist $\theta'$ on $C$. Luckily, since $(\phi')_A = \phi_A$, we also have $\theta'_A = \theta_A = I_A$. 

\end{proof}

The local module $g\in C_A$ is a fermion with half braid $\gamma$ by Proposition~\ref{gSimpleOrder2}. We are working with a strictly pivotal and strictly model of $C_A$. We can form the fermion condensation $C_A/\gamma$, which is a pivotal supercategory. 

Now that we have a pivotal functor $\Psi: BMW(q, -i/q) \to C\to C_A$, we can compose this with the canonical inclusion of $C_A$ into the even part of $(C_A)/\gamma$.

We now have the task of extending $\Psi$ by the obvious functor $BMW(q, -q^{2r + 1}) \to CW_q(-1, 1)$ which maps each BMW diagram to the corresponding BMW diagram in $CW_q(-1, 1)$, to an even super-pivotal functor $\Phi$, as shown in the following commutative diagram, 
\begin{equation}\label{extensionOfBMWC}
\begin{tikzcd}
	{CW_q(-1, 1)} \\
	{BMW(q, -i/q)} & {(C_A)/\gamma}
	\arrow["\Phi", from=1-1, to=2-2]
	\arrow[from=2-1, to=1-1]
	\arrow["\Psi"', from=2-1, to=2-2]
\end{tikzcd}
\end{equation}

There is exactly one extra generator in $CW_q(-1,1)$ which is not in $BMW(q, -i/q)$, which is the spinor dot $\att{oddDot}{0.5}$.

Note now that $Hom_{C_A/\gamma}(X \to X) \cong Hom_{C_A}(X \to X) \oplus Hom_{C_A}(X \to X \tensor g)$, where the last direct summand is the space of odd morphisms. We have by the free forgetful adjunction and self duality of $X$ that $Hom_{C_A}(X \to X \tensor g) \cong Hom_{C_A}(X \tensor X \to g) \cong Hom_C(V_{(1)} \tensor V_{(1)} \to B) \cong \C$. As a result, $End_{C_A/\gamma}(X)\cong Cl_1$, and $X$ is a type $Q$ simple. Up to a sign choice, there is exactly one odd endomorphism $v: X \to X$ in $(C_A)/\gamma$ with the property that $v^2 = 1$. It will not matter which sign choice we have.

Note that by definition of an odd morphism in the fermion condensation, the odd automorphism $v: X \to X$ in $(C_A)/\gamma$ is literally a morphism $X \to X \tensor g$ in $C_A$.

Now define $\Phi\left(\att{oddDot}{0.5}\right) = v$. 

Because all of the relations involving the images of $BMW$ are already satisfied in $C_A$ by diagram~\ref{extensionOfBMWC}, we need only show that $v$ satisfies all of the mixed relations with the images of the morphisms in $BMW(q, -i/q)$ under $\Phi$ with $v$. Since $v^2 = 1$, the relation (Cancel) in $CW(-1, 1)$ is already satisfied. Since $Hom_{(C_A)/\gamma}(X \to X) = \C^{1 | 1}$, and since taking duals is an even functor, and since $X$ is symmetrically self dual (since $V_{(1)}$ is symmetrically self dual with our chosen pivotal structure and $F_A$ is pivotal), it follows that $v^* = \lambda v$ for some $\lambda \in \C$. Since $(C_A)_g$ is super-pivotal and $v$ is odd, $v^{**} = -v$, and hence $\lambda^2 = -1$. So, $\lambda = \pm i$. We need to pin down which sign choice.

Everything (including the sign choice above) amounts to verifying the relation (Half Braid). 

Because the free module functor $F_A: C \to C_A$ is central, the braid on $V_{(1)}$ (which is the image of the braid in $BMW$ under the functor $\Psi$) becomes a half braid on $X$ in $C_A$. For convenience, let $\kappa_X = F_A(c_{V_{(1)}, V_{(1)}}) = \Psi(\att{braid}{0.25}) = \Phi(\att{braid}{0.25})$. More generally, let $\kappa_W$ denote the half braid of $X$ with $W$ which extends the image of the braid under $F_A$, for any object $W$ in $C_A$.

By definition of composition of morphisms and tensor products in the fermion condensation, we have 
\[
\kappa_X \circ_g (I_X \tensor_g v) = (\kappa_X\tensor I_g) \circ (I_X \tensor v).
\]
Abusing notation and writing $\kappa_X$ as a braid, we have by naturality of the half braid for both half braids $\kappa$ on $X$ and $\gamma$ on $g$ that
\[
(\kappa_X\tensor I_g) \circ (I_X \tensor v) =
\att{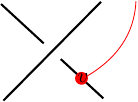}{1} = \att{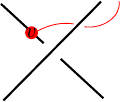}{1} = (I_X\tensor \gamma_X) \circ (v\tensor I_X) \circ \kappa_X
\]
Finally, by definition of the tensor product $\tensor_g$ and composition $\circ_g$ in the fermion condensation $(C_A)/\gamma$, we have 
\[
(I_X\tensor \gamma_X) \circ (v\tensor I_X) \circ \kappa_X = (v \tensor_g I_X) \circ_g \kappa_X.
\]

Putting this all together, we have 
\begin{align*}
\Phi\left(\att{braid}{0.5}\right) \circ_g \left(I_X \tensor_g \Phi\left(\att{oddDot}{0.5}\right) \right) = \kappa_X \circ_g (I_X \tensor_g v) \\
= (v \tensor_g I_X) \circ_g \kappa_X = \left(\Phi\left(\att{oddDot}{0.5}\right) \tensor_g I_X \right) \circ_g \Phi\left( \att{braid}{0.5}\right)
\end{align*}

Hence the relation (Half Braid) is satisfied by the image of $\att{oddDot}{0.5}$ and $\att{braid}{0.5}$ under $\Phi$. 

Now that we have established the relation (Half Braid), we have identitified a subcategory of $C_A$ which extends the relations for $BMW(q, -i/q)$ by an odd morphism satisfying $V^* = \lambda V$ for some $\lambda \in V$. It follows that $\sigma = -1$. This establishes the relation (Rotate).

By the universal property of free pivotal supercategories, it follows that there exists a unique functor $\Psi : CW(-1, 1) \to (C_A)/\gamma$ up to natural isomorphism, which maps $\att{braid}{0.5} \mapsto \kappa_X$ and $\att{oddDot}{0.5} \mapsto v$. 

Since odd generating morphisms were sent to odd morphisms, and even generating morphisms were sent to even morphisms, $\Phi: CW(-1, 1) \to (C_A)/\gamma$ is an even functor. 

As a result, $\Phi$ restricts to the even part, thus giving giving a pivotal functor $\Psi[0]: CW(-1, 1)[0] \to C_A$ making the diagram commute:
\[
\begin{tikzcd}
	{CW_q(-1, 1)} & {C_A/\gamma} \\
	{CW_q(-1, 1)[0]} & {C_A}
	\arrow["\Phi", from=1-1, to=1-2]
	\arrow[hook', from=2-1, to=1-1]
	\arrow["{\Phi[0]}", from=2-1, to=2-2]
	\arrow[hook', from=2-2, to=1-2]
\end{tikzcd}
\]

\subsubsection{Faithfulness} 

By Theorem~\ref{smallSpanningSet}, it follows that if the relations in $CW_q(-1, 1)$ are consistent, then $CW_q(-1, 1)$ has simple unit, and hence $CW_q(-1, 1)[0]$ has simple unit. At this point, we have found a category satisfying the relations of the Clifford BMW category inside $(C_A)/\gamma$, which is manifestly non-zero. This proves that the relations for $CW_q(-1, 1)$ are consistent, and that the identity on the unit object $0$ is non-zero. As a result, we conclude that $CW_q(-1, 1)[0]$ and $CW_q(-1, 1)$ have simple units. 

The category $C_A$ is unitary. The pivotal category $CW_q(-1, 1)[0]$ is spherical. The functor $\Psi[0]$ is pivotal. All of this allows us to use Proposition~\ref{semisimplification} to conclude that $\Psi[0]$ descends via the semisimplification functor to a faithful functor $\sms{\Phi[0]}:\sms{CW_q(-1, 1)[0]} \to C_A$. Moreover, by Proposition~\ref{superDescent} we have a faithful functor $\sms{\Phi} :\sms{CW_q(-1, 1)} \to (C_A)/\gamma$, and $\sms{\Phi}[0] = \sms{\Phi[0]}$.

\subsubsection{Fullness}

Thus far, we have shown the existence of a dominant, faithful functor $\sms{\Phi[0]}:CW_q(-1, 1)[0] \to C_A$. Our final task is to show that $\sms{\Phi[0]}$ is full. Using this, we will also show $\sms{\Phi}$ is also full. 

Following the standard procedure laid out in section 2 of \cite{CainAndSnyder}, we can take the Cauchy completions and semisimplifications of $BMW(q, -i/q)$ and $CW_q(-1, 1)[0]$. Doing so gives us the leftmost square of the following diagram, which commutes up to natural isomorphisms, where the unlabeled arrows (other than the equivalence) come from the obvious functor $BMW(q, -i/q) \to CW_q(-1, 1)[0]$ or the inclusions into the abelianizations, and where $F_1$ is defined so that the diagram commutes. 

\begin{tikzcd}
	&& {Ab(\overline{CW_q(-1, 1)[0]})} && \\
	{Ab(\overline{BMW(q, -i/q)})} &&&& {C_A} \\
	\\
	&& C
	\arrow["{Ab(\overline{\Phi[0]})}", from=1-3, to=2-5]
	\arrow[from=2-1, to=1-3]
	\arrow["{F_1}"', from=4-3, to=1-3]
	\arrow["{Ab(\overline{G})^{-1}}", from=4-3, to=2-1]
	\arrow["{F_A}"', from=4-3, to=2-5]
\end{tikzcd}

\begin{lemma}\label{galoisCorrespondenceC}
    The category $Ab(\sms{CW_q(-1, 1)})[0]$ is a unitary fusion category. There exists a subalgebra $A' \sbs A$ and a monoidal equivalence $ K:Ab(\sms{CW_{2r + 1}(q, iq, -1)[0]}) \to C_{A'}$ such that the diagram commutes up to natural isomorphism, 
\begin{equation}\label{GaloisDiagramC}
\begin{tikzcd}
	C & {Ab(\sms{CW_q(-1, 1)[0]})} & {C_A} \\
	& {C_{A'}}
	\arrow["{F_1}", from=1-1, to=1-2]
	\arrow["{F_A}"', curve={height=-30pt}, from=1-1, to=1-3]
	\arrow["{F_{A'}}"', from=1-1, to=2-2]
	\arrow["{Ab(\sms{\Phi[0]})}", from=1-2, to=1-3]
	\arrow["K"', from=1-2, to=2-2]
\end{tikzcd}\end{equation}
and such that $A'$ is either a pointed subalgebra, or $A$ itself.
\end{lemma}
\begin{proof}
    The functor $F_1$ is a composition of an equivalence and the induced map coming from the extension $BMW(q, -i/q) \to CW_q(-1, 1)[0]$, and hence is a composition of dominant functors. Moreover, by construction, this induced functor is faithful. Hence $F_1$ is a composition of two dominant, faithful, monoidal functors, and hence $F_1$ is dominant, faithful, and monoidal. 

    As proven in the last sub-subsection, $Ab(\sms{CW_q(-1, 1)[0]})$ has simple unit, and is semisimple. Since $F_A$ is dominant, and $F_A = Ab(\sms{\Phi[0]})\circ F_1$, we must have that $Ab(\sms{\Phi[0]})$ is dominant (and it is obviously monoidal). As a result, the hypotheses of Proposition~\ref{Galois} are satisfied, and the lemma follows immediately. The subalgebras of $A$ are fully understood, and known to be either pointed or $A$ itself \cite{CainAndSnyder}. 
\end{proof}

\begin{lemma}\cite{Fuchs}
    The only invertible simple objects in $\mathcal{C}(\mathfrak{sp}_{2r}, k)$ are $V_{k\Lambda_{r}}= V_{(k^r)}$ and $V_{(0)}$. In particular, the only non-unit invertibles object in $\mathcal{C}(\mathfrak{sp}_{2r}, r + 1)$ is $J = V_{((r + 1)^r)}$. For a weight $\lambda = x_1 \Lambda_1 + \dots + x_r \Lambda_r$ with $\sum_{i}x_i \le k$ and with $x_0 = k - \sum_{i = 1}^n x_i$, one has
    \[
    J\tensor V_{\lambda} = V_{x_{r -1}\Lambda_1 + \dots + x_1 \Lambda_{r - 1} + x_0 \Lambda_r}
    \]
    In particular, for $r > 3$ and $k \ge 2$,
    \[
    J \tensor V_{(1,1)} = V_{((r - 1)^{r - 2}, 1, 1)}, \hspace{0.5cm} J\tensor V_{(2)} = V_{((r + 1)^{r - 1}, r - 1)}.
    \]
\end{lemma}
\begin{proof}
    The coxeter labels $a_i$ for the root system $C_r$ are $a_i = \begin{cases}
        2 & i < r\\
        1 & i = r
    \end{cases}$. By Fuchs' classification, the invertible simple objects in $\mathcal{C}(\mathfrak{sp}_{2r}, k)$ are $V_{(0)}$ and $V_{k\Lambda_{(i)}}$ where $a_i = 1$ \cite{Fuchs}. That $J$ is the only invertible simple object follows immediately. The rule for tensoring with $J$ is given in Section 3.1 of \cite{Gannon}. The particular examples in the statement of the lemma are just routine exercises in converting from the $e_i$ coordinates into the $\Lambda_i$ coordinates, applying Terry Gannon's rule, and then converting back into the $e_i$ coordinates. 
\end{proof}

\begin{corollary}
    Any subalgebra $A' \sbs A$ is of the following forms:
    \[
    A' = A, \hspace{0.5cm} A' = J \oplus V_{(0)}, \hspace{0.5cm} A' = V_{(0)}.
    \]
    If $r > 3$, then for any subalgebras $A' \ne A$ of $A$, we have 
    \[
    dim Hom_{C_{A'}}(F_{A'}(V_{(1)})^{\tensor 2} \to F_{A'}(V_{(1)})^{\tensor 2}) \le 4
    \]
    
\end{corollary}

\begin{proof}
    If $A'= V_{(0)}$, the corollary is immediate. 

    If $A' = J \oplus V_{(0)}$, then by the free forgetful adjunction, and our description of $J\tensor V_{(1,1)}$, $J\tensor V_{(2)}$ and $J$, 
    \begin{align*}
    &Hom_{C_{A'}}(F_{A'}(V_{(1)})^{\tensor 2} \to F_{A'}(V_{(1)})^{\tensor 2})\\ &\cong Hom_\mathcal{C}(V_{(1)}^{\tensor 2} \to (J \oplus V_{(0)}) \tensor V_{(1)}^{\tensor 2})\\
    &\cong Hom_\mathcal{C}(V_{(0)} \oplus V_{(1,1)}\oplus V_{(2)}
    \to V_{(0)} \oplus V_{(1,1)}\oplus V_{(2)} \oplus V_{(r + 1)^r} \oplus V_{((r - 1)^{r - 2}, 1, 1)} \oplus V_{((r + 1)^{r - 1}, r - 1)}).
    \end{align*}
    If $r = 2$, then $V_{((r - 1)^{r - 2}, 1, 1)}  =V_{(1, 1)}$, while $|(r + 1)^r| = 9$ and $|((r + 1)^{ r - 1}, r - 1)| = 4$. So, in the case where $r = 2$, we have 
    \[
    \dim Hom_{C_{A'}}(F_{A'}(V_{(1)})^{\tensor 2} \to F_{A'}(V_{(1)})^{\tensor 2}) \le 4 < 6.
    \]
    On the other hand, if $r > 2$, we have 
    \[
    \dim Hom_{C_{A'}}(F_{A'}(V_{(1)})^{\tensor 2} \to F_{A'}(V_{(1)})^{\tensor 2}) \le 3 < 6.
    \]
\end{proof}

\begin{remark}
    It is interesting to note that, unlike in the type B case, $A$ really does occasionally have a non-trivial pointed subalgebra, namely $J \oplus V_{(0)}$. The object $J$ corresponds to the abelian ideal which in turn corresponds to the strict Young diagram $\delta = (r, r - 1, \dots, 1)$. This abelian ideal is of even size when either $4 | r$ or $4 | r + 1$.  
\end{remark}

For $ r\ge 2$, $q = e^{\frac{2\pi i}{8r + 8}}$ is a not a primitive $16$th root of unity. Since we have represented the category $CW_q(-1, 1)[0]$ non-trivially in $C_A$, it follows by Lemma~\ref{constrantsOnHom2to2} that $\dim End_{\sms{CW_q(-1, 1)[0]}}(2) = \dim End_{Ab(\sms{CW_q(-1, 1)[0]})}(2) = 6$.

We are now in a position to prove fullness of the functor $Ab(\sms{\Phi[0]})$. If this functor is not full, then there is a monoidal equivalence $K:Ab(\sms{CW_q(-1, 1)[0]})\to C_{A'}$ for $A'$ a proper subalgebra, and such that the diagram~\ref{GaloisDiagramC} commutes up to natural isomorphism. This gives us an isomorphism of endomorphism algebras
\[
End_{Ab(\sms{CW_q(-1, 1)[0]})}(2) \cong End_{C_{A'}}(F_{A'}(V_{(1)}^{\tensor 2})).
\]
But the right hand side has dimension $6$, while the left hand side has dimension at most $4$. Since this cannot happen, we conclude that the functor $Ab(\sms{\Phi[0]}) :Ab(\sms{CW_q(-i/q, 1)[0]})\to C_A$ is full. 

Since $\sms{\Phi[0]}$ was dominant, it follows that $Ab(\sms{\Phi[0]}) $ is essentially surjective, faithful, full and hence an equivalence. 

To see that $\sms{\Phi}$ is full, note that since $End_{C_A/\gamma}(X) \cong Cl_1$, $X$ is a type $Q$ simple, and hence $g \tensor X \cong X$. As a result, for every $N \ge 1$, we have an $End_{C_A}(X^{\tensor N}) \cong Hom_{C_A} (X^{\tensor N} \to g\tensor X^{\tensor N})$
which maps for all $f \in End_{C_A}(X^{\tensor N})$ as
\[
f \mapsto (v \tensor I_X^{\tensor N})\circ f.
\]
Since $v$ is in the image of the functor $\Phi$ by definition, we have that $End_{C_A}(X^{\tensor N})_1 \sbs \Phi(End_{CW_q(-1, 1)})(n))$.

This shows that $\Phi$ is full. As a result, $\sms{\Phi}$ is also full and faithful. It remains to see that $\sms{\Phi}$ is dominant. This is obvious from the fact that the objects in $C_A$ are literally the objects in $C_A/\gamma$, and every object in $C_A$ is contained in a tensor power of $X$. 

Since $\sms{\Phi}$ is full, faithful, even, and dominant, $Ab(\sms{\Phi})$ is full, faithful, and essentially surjective. Hence $Ab(\sms{\Phi})$ is an equivalence which is also superpivotal.

This concludes the proof of Theorem~\ref{bigResultC}.
\end{proof}

\subsubsection{Deligne Interpolation}

As noted, the Clifford BMW categories are all conjugate to one another under the automorphisms of the ring $R$. As a result, we need only carry out Deligne interpolation for $CW(-1, 1)$, from which our explicit dimension formulas will follow for all values of $\sigma$ and $\epsilon$.

\begin{remark}
    Interpolation results are a common phenomenon in the world of tensor categories. It should be noted, however, that there are essentially two kinds of interpolation results which often get blurred together. In cases of the classical interpolation categories such as $Rep(O_t)$, $Rep(S_t)$ or $Rep(GL_t)$, the linear independence of a canonical basis can be established by directly constructing each of these categories. In these cases, the interesting thing is not so much that these categories exist, but that they reveal regularities in a discrete family of categories, to such an extent that these categories can be defined continuously. On the other hand, there are other interpolation results which also serve as the main proof of the linear independence of a canonical basis at a continuous set of parameters. The following theorem essentially fits into both camps, as does the main results of \cite{CainAndSnyder}. For us, interpolation establishes the existence of a continuous family, while it will also help us to get an idea of the internal structure of the categories $C_A$.
\end{remark}

\begin{theorem}\label{dimensionOfCW}
    For all $\sigma, \epsilon\in \set{\pm 1}$ and for all $n\in \N$, $End_{CW(\sigma, \epsilon)}(n)$ is a free $R$-module of rank $2^{n}(2n - 1)!!$. Moreover, for all $q\in \C$ such that $q^{8}\ne 0, 1$ and for all $n\in \N$ we have
    \[
    \dim End_{CW_q(\sigma, \epsilon)(n)} = 2^{n} (2n - 1)!!
    \]
\end{theorem}
\begin{proof}

    For $r \ge 2$ and  $r > 2n$, Theorem~\ref{dimensionOfCAThypeC} gives $\dim End_{C_A}(X^{\tensor n}) = 2^{n - 1}(2n - 1)!!$. Since $g\tensor X \cong X$, we have $\dim End_{C_A/\gamma}(X^{\tensor n}) = 2^{n}(2n - 1)!!$. So, since $Ab(\sms{CW_q(-1, 1)}) \simeq C_A/\gamma$, combining this with Theorem~\ref{smallSpanningSet} yields
    \begin{align*}
    2^{n }(2n - 1)!! &= \dim End_{Ab(\sms{CW_{q_r}(-1, 1)})}(n)) \le \dim End_{CW_{q_r}(-1, 1)}\le 2^{n }(2n - 1)!!
    \end{align*}
    where $q_r = e^{\frac{2\pi i}{4(2r + 2)}}$, and hence $\dim End_{CW_{q_r}(-1, 1)} = 2^{n - 1}(2n - 1)!!$. From this it follows that for $q_r = e^{\frac{2\pi i}{4(2r + 2)}}$, $ r\ge 3$, and $r > 2n$, the spanning set $\mathcal{B}$ from Theorem~\ref{smallSpanningSet} is a basis for $End_{CW_{q_r}(-1, 1)}$. Suppose by way of contradiction that we had some linear dependence in $R$ of the spanning set $\fanB$, $\sum_{d\in \fanB} f_d(q) d$ where $f_d(q)$ is some rational function in $R$. Then by linear independence of $\fanB$ for specializations at $q_r$, we have $f_d(q_r) = 0$ for all $r \ge 3\in \N$ and for all $d \in \fanB$. The only rational function with infinitely many zeros is zero, hence $f_d = 0$ for each $d\in \fanB$. As a result, $\fanB$ is linearly independent over $R$. 

    By Lemma~\ref{specializationForCW}, for any $q\in\C$ such that $q^8\ne 1, 0$ we have a specialization map $ev_q: R \to \C$. Note that this is a retraction of $R$ onto $\C$. By Lemma~\ref{specializationAndConjugation}, it follows that $\fanB$ is linearly independent for all specializations to $q\in \C$ such that $q^8 \ne 1, 0$. 

\end{proof}

\subsection{$C_A$ for type B}

The goal of this section is to prove the following theorem.
\begin{theorem}\label{bigResultB}
    Let $q = e^{\frac{2\pi i}{4(4r - 2)}}$. There exists a full and dominant even superpivotal functor 
    \[
    \Phi: SCW_{2r + 1}(q) \to (\mathcal{C}(\so_{2r+1}, 2r - 1)^{vec}_A)/\gamma
    \]
    which descends to a full and faithful dominant functor 
    \[
    \sms{\Phi} : \sms{SCW_{2r + 1}(q)} \to (\mathcal{C}(\so_{2r+1}, 2r - 1)^{vec}_A)/\gamma
    \]
    thus giving an even super-pivotal equivalence
    \[
    Ab(\sms{SCW_{2r + 1}(q)}) \simeq (\mathcal{C}(\so_{2r+1}, 2r - 1)^{vec}_A)/\gamma
    \]
    which restricts to a pivotal equivalence 
    \[
    Ab(\sms{SCW_{2r + 1}(q)[0]}) \simeq \mathcal{C}(\so_{2r+1}, 2r - 1)^{vec}_A
    \]
\end{theorem}

\begin{proof}

The proof outline will be similar to the type C case, only now we have to deal with an extra $2r + 1$-valent vertex in our skein theory.

Throughout the remainder of this section let $r \ge 3\in \N$ be fixed, and let $C = \mathcal{C}(\mathfrak{so}_{2r + 1}, 2r - 1)^{vec}$. Note that for $r = 2$, $ \mathcal{C}(\mathfrak{so}_{2r + 1}, 2r - 1)^{vec} \simeq \mathcal{C}(\mathfrak{sp}_{2r }, r + 1)$ and that the Etale algebras corresponding to their respective conformal embeddings are equal. We will be working with a stricly pivotal and strictly monoidal model of $C_A$. Note that $g\in C_A$ admits a half braid $\gamma$ which turns it into a fermion, for every rank. Forming the fermion condensation $C_A/\gamma$ gives a striclty monoidal and strictly pivotal supercategory.

\subsubsection{Defining the functor}

By Proposition~\ref{specialBMWToC} we have a braided pivotal equivalence $Ab(\sms{SW_{2r + 1}(q)}) \simeq C$. Precomposing with the abelianization and semismiplification functors, we have a pivotal braid functor $SW_{2r + 1}(q) \to C$. Composing with the free module functor gives a pivotal central functor $\Psi[0]: SW_{2r + 1}(q) \to C_A$. Composing with the canonical even inclusion $C_A \hookrightarrow C_A/\gamma$ gives a pivotal central functor $\Psi: SW_{2r + 1}(q) \to C_A$. 

Our task is now to extend $\Psi$ to $CW_{2r + 1}(q)$ by specifying what to do with the odd generator $\att{oddDot}{0.5}$, so that the diagram commutes:

\begin{tikzcd}
	& {SCW_{2r + 1}(q)} \\
	{SW_{2r + 1}(q)} & {C_A/\gamma}
	\arrow["\Phi", from=1-2, to=2-2]
	\arrow[from=2-1, to=1-2]
	\arrow["\Psi"', from=2-1, to=2-2]
\end{tikzcd}

As before, by our decomposition of the algebra $A$ and the object $B$, and by the free forgetful adjunction, and by self duality of $X$
\[
Hom_{C_A}(X \to X \tensor g) \cong Hom_{C_A}(X \tensor X \to g) \to \cong Hom_{C}(V_{(1)} \tensor V_{(1)} \to B) \cong \C. 
\]
As a result, in the fermion condensation $C_A/\gamma$ we have 
\[
Hom_{C_A/\gamma}(X \to X) \cong \C^{1 | 1}. 
\]
Hence there is some odd automorphism $v: X \to X$ in $C_A/\gamma$ such that $v^2 = I_X$. 

Because $2r + 1$ is odd, we have no mixed relations between the $2r + 1$-valent vertex $S$ and $v$. The relation (Half Braid) follows from the fact that $F_A$ is central in the same way as it did for type C. 

Since $End_{C_A/\gamma}(X)_1 \cong \C$ and $X$ is self dual we must have $v^* = \lambda v$ for some $\lambda \in \C$. Hence $\lambda^2 = -1$ and $\lambda = -\sigma i q$ for some $\sigma \in \set{\pm 1}$. By the Loop Lemma we have $\sigma = 1$. This establishes the relation (Rotate). 

All of this together establishes the existence of the superpivotal functor $\Phi$. Since even (odd) generating morphisms were sent to even (odd) generating morphisms, the functor $\Phi$ is even. Restricting to even parts gives a pivotal functor $\Phi[0] : SCW_{2r + 1}(q)[0] \to C_A$ such that the diagram commutes:
\[
\begin{tikzcd}
	{SCW_{2r + 1}(q)} & {C_{A}/\gamma} \\
	{SCW_{2r + 1}(q)[0]} & {C_A}
	\arrow["\Phi", from=1-1, to=1-2]
	\arrow[hook, from=2-1, to=1-1]
	\arrow["{\Phi[0]}", from=2-1, to=2-2]
	\arrow[hook, from=2-2, to=1-2]
\end{tikzcd}.
\]

Note that since $F_A$ is dominant and every object in $C$ is in a tensor power of $V_{(1)}$, and hence every object of $C_A$ is in a tensor power of $F_A(V_{(1)}) = X$, it follows that $\Phi[0]$ is dominant. Since the objects of $C_A/\gamma$ are literally the objects of $C_A$, the same goes for $\Phi$. 

\subsubsection{Faithfulness}

By Theorem~\ref{evaluationAlgorithmSpecialCliffBMW}, and since we have found a non-trivial representation of $SCW_{2r + 1}(q)$ in $C_A$, it follows that the unit objects in $SCW_{2r + 1}(q)[0]$ and $SCW_{2r + 1(q)[0]}$ are simple. 

By Proposition~\ref{CASemisimple}, $C_A$ and $C_A/\gamma$ are unitary. Moreover, $\Phi$ is pivotal, and $SCW_{2r + 1}(q)$ and its even part are spherical. By Proposition~\ref{semisimplification} (Proposition 1.9 of~\cite{CainAndSnyder}), there exists a faithful functor $\sms{\Phi}$ such that the diagram commutes, 
\[\begin{tikzcd}
	{SCW_{2r + 1}(q)} \\
	{\overline{SCW_{2r + 1}(q)}} & {C_A/\gamma}
	\arrow[from=1-1, to=2-1]
	\arrow["\Phi", from=1-1, to=2-2]
	\arrow["{\overline{\Phi}}"', from=2-1, to=2-2]
\end{tikzcd}\]
where the leftmost arrow is the quotient functor onto the semisimplification. Note that the functor $\sms{\Phi}$ is dominant because $\Phi$ is dominant. 

The restriction of the $\sms{\Phi}$ to even parts is $\sms{\Phi}[0] = \sms{\Phi[0]}$, and is faithful and dominant. 

\subsubsection{Fullness}

As in the type $C$ case, taking Cauchy completions and semisimplifications gives a chain of pivotal functors
\[\begin{tikzcd}
	C & {Ab(\overline{SW_{2r + 1}(q, iq)})} & {Ab(\sms{SCW_{2r + 1}(q)})[0]} & {C_A}
	\arrow["\simeq", from=1-1, to=1-2]
	\arrow[from=1-2, to=1-3]
	\arrow["{Ab(\overline{\Phi}[0])}", from=1-3, to=1-4]
\end{tikzcd}\]

Note that this composition of functor maps $\att{braid}{0.25}\mapsto F_A(\att{braid}{0.25})$ and so on. Since $C$ is given a presentation in terms of these generating morphisms, it follows that this composition of functors is isomorphic to $F_A$, again by the universal property of free pivotal categories. 



Note that since $\sms{\Phi[0]}$ is faithful, so is $Ab(\sms{\Phi}[0])$. Similarly, since $\Phi[0]$ is dominant, so is $Ab(\sms{\Phi}[0])$. Since $C_A$ is unitary (Proposition~\ref{CASemisimple}) and $Ab(\sms{\Phi})$ is faithful, we have that $Ab(\sms{SCW_{2r + 1}(q)[0]})$ is unitary. In particular, $Ab(\sms{SCW_{2r + 1}(q)[0]})$ is semisimple. From this we can apply Proposition~
\ref{Galois} to prove the following,

\begin{lemma}\label{galoisCorrespondenceB}
    The category $Ab(\sms{SCW_{2r + 1}(q)})$ is finitely semisimple. There exists a subalgebra $A' \sbs A$ and a monoidal equivalence $ K:Ab(\sms{SCW_{2r + 1}(q)}) \to C_{A'}$ such that the diagram commutes up to natural isomorphism, 
\[\begin{tikzcd}
	C & {Ab(\sms{SCW_{2r + 1}(q)[0]})} & {C_A} \\
	& {C_{A'}}
	\arrow["{F_1}", from=1-1, to=1-2]
	\arrow["{F_A}"', curve={height=-30pt}, from=1-1, to=1-3]
	\arrow["{F_{A'}}"', from=1-1, to=2-2]
	\arrow["{Ab(\sms{\Phi}[0])}", from=1-2, to=1-3]
	\arrow["K"', from=1-2, to=2-2]
\end{tikzcd}\]
such that $A'$ is either pointed, or $A$ itself.
\end{lemma}
\begin{proof}
    By an argument exactly analogous to the proof of Lemma~4.26 of \cite{CainAndSnyder}, the category $Ab(\sms{SCW_{2r + 1}(q)[0]})$ is a unitary fusion category. Following the same procedure as in the proof of Lemma~4.26 in \cite{CainAndSnyder}, all proper subalgebras of $A$ are pointed. Using the Galois correspondence between intermediate functors which factor the free module functor, and subalgebras of $A$, as described in Proposition~\ref{Galois}, there must exist an equivalence $K$ which makes the above diagram commute up to natural isomorphisms. 
\end{proof}

\begin{lemma}\cite{Fuchs}\cite{Gannon}

The only proper subalgebra of $A$ is the unit object $V_{(0)}$. 
\end{lemma}
\begin{proof}
    Invertible objects (there called simple currents) in $\mathcal{C}(\mathfrak{so}_{2r + 1}, k = 2r - 1)$ are classified in \cite{Fuchs}. By this classification there can only be one invertible object other than $V_{(0)}$, namely $J = V_{k {\Lambda_1}}$. In our case, since $k = 2r - 1$, it suffices to show that $V_{(2r - 1)}$ is not a subobject of $A$. Recall that the subobjects of $A$ are parametrized by sets $\Psi$ of positive roots which form an abelian subalgebra of the root system $B_r$. By our description in Lemma~ \ref{idealsTypeB} any ideal $\Psi\sbs B_r$ must contain the root $e_1 + e_2$. In order for the weight $\lambda_\Psi$ not to have an $e_2$ coordinate, it must be the case that $e_1 - e_2 \in \Psi$. But then $e_1 = (e_1 - e_2) + (e_1 + e_2) \in B_r$, and hence $\Psi$ is not abelian. 
    
\end{proof}

\begin{lemma}\label{homSetsBigB}
    For $q = e^{\frac{2\pi i}{4(4r - 2)}}$ $\dim End_{Ab(\sms{SCW_{2r + 1}(q)[0]})}(2) \ge 6 $
\end{lemma}
\begin{proof}

    Note that $q = e^{\frac{2\pi i}{4(4\cdot r - 2)}}$ is a $16r - 8$-th root of unity. But $16r - 8 \ge 24 > 16$ since $r \ge 2$. As a result, since we have represented $SCW_{2r + 1}(q)$ non-trivially and therefore shown the relations to be consistent, by Lemma~\ref{constrantsOnHom2to2} we have $ \dim End_{Ab(\sms{SCW_{2r + 1}(q)[0]})}(2) \ge 6 $. 
\end{proof}

We are now in a position to prove fullness. Suppose by way of contradiction that $Ab(\overline{\Phi[0]}): Ab(\sms{SCW_{2r + 1}(q)[0]}) \to C_A$ is not full. Then $F_A$ factors through $Ab(\sms{SCW_{2r + 1}(q)[0]})$, so by Lemma~\ref{galoisCorrespondenceB} we have $Ab(\sms{SCW_{2r + 1}(q)[0]}) \simeq C_{A'}$ for some proper subalgebra $A' \sbs A$, with an equivalence that sends $1^{\tensor 2} \mapsto F_{A'}(V_{\Lambda_1})^{\tensor 2}$. This gives isomorphisms of endomorphism algebras
\[
End_{Ab(\sms{SCW_{2r + 1} (q)[0]}) }(2) \cong End_{\sms{SCW_{2r + 1} (q)[0]}}(2) \cong End_{C_{A'}}(F_{A'}(X^{\tensor 2})).
\]
But this cannot happen, as by the previous lemma the dimension of the first endomorphism algebra is at least $6$, while the dimension of the last must be $3$. This concludes the proof of Theorem~\ref{bigResultB}.

\end{proof}

\subsection{$C_A$ for type D}

Let $C = \mathcal{C}(\mathfrak{so}_{2r}, 2r - 2)^{vec}$ with $q = e^{\frac{2\pi i}{4(4r - 4)}}$. 

\begin{theorem}\label{bigResultD}
    Let $r \ge 4$. Let $q = e^{\frac{2\pi i}{4(4r - 4)}}$ 
    There exists a full and dominant even superpivotal functor 
    \[
    \Phi: SCW_{2r}(q) \to (\mathcal{C}(\so_{2r}, 2r - 2)^{vec}_A)/\gamma
    \]
    which descends to a full and faithful dominant functor 
    \[
    \sms{\Phi} : \sms{SCW_{2r}(q)} \to (\mathcal{C}(\so_{2r}, 2r - 2)^{vec}_A)/\gamma
    \]
    thus giving an even super-pivotal equivalence
    \[
    Ab(\sms{SCW_{2r}(q)}) \simeq (\mathcal{C}(\so_{2r}, 2r - 2)^{vec}_A)/\gamma
    \]
    which restricts to a pivotal equivalence 
    \[
    Ab(\sms{SCW_{2r}(q)[0]}) \simeq \mathcal{C}(\so_{2r}, 2r - 2)^{vec}_A
    \]
\end{theorem}

\subsubsection{Defining the Functor}

Throughout the remainder of this subsection, let $C = \mathcal{C}(\mathfrak{so}_{2r}, 2r - 2)^{vec}$. By Theorem 2.9 of \cite{qMacky} we get that $\theta_A = I_A$. By Theorem 1.17 of \cite{qMacky}, $C_A$ is endowed with a pivotal structure such that the free module functor $F_A:C \to C_A$ is pivotal. We will work with a model of $C_A$ that is strictly pivotal and strictly monoidal. The local module $g\in C_A$ is a fermion with half braid $\gamma$, so we can form the fermion condensation $C_A/\gamma$.

The objects of $SCW_{2r}(q)$ are the natural numbers including $0$. This makes it easy to define $\Phi$ on objects. Let $\Phi(n) = X^{\tensor n}$. 

Let $F_A: C \to C_A$ be the free module functor. By Proposition~\ref{specialBMWToC} there exists an equivalence $Ab(\sms{SW_{2r}(q)}) \simeq C$, which is a braid functor, and a pivotal functor. Precomposing this equivalence with the Cauchy completion functor $SW_{2r}(q) \to Ab(\sms{SW_{2r}(q)})$, we have a dominant, pivotal functor $SW_{2r}(q) \to C$. Composing with the free module functor $F_A$ gives a dominant and pivotal functor $\Psi[0] : SW_{2r + 1}(q) \to C_A$. Composing with the inclusion of $C_A$ into the $C_A/\gamma$ gives a dominant and pivotal superfunctor $\Psi: SW_{2r + 1}(q) \to C_A/\gamma$. Our task is to extend $\Psi$ to $\Phi$ through the obvious inclusion $SCW_{2r}(q) \hookrightarrow SCW(q)$, so that the diagram commutes:
\begin{equation}\label{extensionTypeD}
    \begin{tikzcd}
	{SCW_{2r}(q)} \\
	{SW_{2r}(q, iq)} & {C_A/\gamma}.
	\arrow["\Phi", from=1-1, to=2-2]
	\arrow[from=2-1, to=1-1]
	\arrow["\Psi"', from=2-1, to=2-2]
\end{tikzcd}
\end{equation}

By the same exact reasoning as the last two steps, $X$ in $C_A/\gamma$ is a type $Q$ simple, and we have an odd automorphism $v: X \to X$ so that $v^2 = 1$ and $v^* = -\sigma i v$ for some $\sigma \in \set{\pm 1}$. By exactly the same argument as in the type B and C cases, $v$ satisfies the relation (Half Braid) of the Clifford BMW category, with the image of $\att{braid}{0.25}$ under the functor $\Psi$. As a result, by the Twist Lemma, since $\att{twist}{0.5} = iq \att{stick}{0.5}$ we have 
\[
\sigma = 1.
\]
This establishes that $v$ satisfies the relation (Rotate) for $CW(1, -1)$, as well as the relations (Half Braid) and (Cancel). 

If $r$ is odd, there are no mixed relations involving $\att{oddDot}{0.5}$ and the generator $S$ in $SCW_{2r}(q)$. If $r$ is even, we have the relation (Pedals) to verify. Forming the diagram $\att{dotsOnPedels}{0.5}$ in $C_A/\gamma$ with $v$ replacing the odd dot and $F_A(S)$ replacing $S$, we have
\[
\att{dotsOnPedels}{1} = \lambda
\]
for some scalar $\lambda$ because the unit object in $C_A/\gamma$ is simple. 

Rotating 180 degrees and bringing the two $2r$-valent vertices reveals 
\[
\lambda^2 = \att{stupidNumber}{0.5},
\]
as shown in subsection~\ref{SCW}. This shows that the relation (Pedals) holds as well. 

As a result, mapping $\att{oddDot}{0.5}\mapsto v$ does indeed induce an even superpivotal functor $\Phi: SCW_{2r}(q) \to C_A/\gamma$ which extends $\Psi: SW_{2r}(q)$.

\subsubsection{Faithfulness}

By Theorem~\ref{evaluationAlgorithmSpecialCliffBMW}, the unit object $0$ in $SCW_{2r}(q)$ is simple, since the relations are consistent (having modeled them in $C_A$).

By Proposition~\ref{CASemisimple}, we know that $C_A$ is unitary. As noted before, $\Phi$ is pivotal. By Proposition~\ref{semisimplification} (Proposition 1.9 of~\cite{CainAndSnyder}), there exists a faithful functor $\sms{\Phi}$ such that the diagram commutes, 
\[\begin{tikzcd}
	{SCW_{2r}(q)} \\
	{\overline{SCW_{2r}(q)}} & {C_A/\gamma}
	\arrow[from=1-1, to=2-1]
	\arrow["\Phi", from=1-1, to=2-2]
	\arrow["{\overline{\Phi}}"', from=2-1, to=2-2]
\end{tikzcd}\]
where the leftmost arrow is the quotient functor onto the semisimplification. Note that the functor $\sms{\Phi}$ is dominant. We also get $\sms{\Phi}[0] : SCW_{2r}(q)[0] \to C_A$

\subsubsection{Fullness}

By the exact same argument as in the type B and C cases, we have 
\begin{lemma}\label{galoisCorrespondenceD}
    The category $Ab(\sms{SCW_{2r}(q)[0]})$ is a unitary fusion category. There exists a subalgebra $A' \sbs A$ and a monoidal equivalence $ Ab(\sms{SCW_{2r}(q)[0]}) \to C_{A'}$ such that the diagram commutes up to natural isomorphism, 
\[\begin{tikzcd}
	C & {Ab(\sms{SCW_{2r}(q)[0]})} & {C_A} \\
	& {C_{A'}}
	\arrow["{F_1}", from=1-1, to=1-2]
	\arrow["{F_A}"', curve={height=-30pt}, from=1-1, to=1-3]
	\arrow["{F_{A'}}"', from=1-1, to=2-2]
	\arrow["{Ab(\sms{\Phi[0]})}", from=1-2, to=1-3]
	\arrow["K"', from=1-2, to=2-2]
\end{tikzcd}\]
such that $A'$ is either pointed, or $A$ itself.
\end{lemma}

\begin{proposition}\cite{Fuchs}\cite{Gannon}
    There are exactly four invertible simple objects in $\mathcal{C}(\mathfrak{so}_{2r} , 2r - 2)$. These are 
    \[
    V_{(0)}, V_{(2r - 2)}, V_{(2r - 2)\Lambda_{r - 1}} = V_{((r - 1^r))^-} , V_{(2r - 2)\Lambda_r} = V_{((r - 1^r))^+}
    \]
    
    Tensoring with $V_{(2r - 2)}$ is given by the rule
    \[
    V_{(2r - 2)} \tensor V_{x_1 \Lambda_1 + \dots + x_r \Lambda_r} = V_{x_0 \Lambda_1 + x_2\Lambda_2 + \dots + x_{r - 2} \Lambda_{r - 2} + x_r \Lambda_{r - 1} + x_{r - 1} \Lambda_r}
    \]
    where $x_0 = 2r - 2  - x_1 - x_{r - 1} - x_r - 2\sum_{i = 2}^{r - 2} x_i$

    Tensoring with $V_{((r - 1^r))^+}$ depends on the parity of $r$. If $r$ is even, then 
    \[
    V_{((r - 1)^r)^+} \tensor V_{x_1 \Lambda_1 + \dots + x_r \Lambda_r} = V_{x_{r - 1} \Lambda_{1}  + x_{r - 2}\Lambda_{2} + \dots + x_1 \Lambda_{r - 1} + x_0 \Lambda_r}
    \]
    and if $r$ is odd, then
    \[
    V_{((r - 1)^r)^+} \tensor V_{x_1 \Lambda_1 + \dots + x_r \Lambda_r} = V_{x_{r} \Lambda_{1}  + x_{r - 2}\Lambda_{2} + \dots + x_1 \Lambda_{r - 1} + x_0 \Lambda_r}
    \]

    Similarly, tensoring with $V_{((r - 1)^r)^-}$ also depends on sign parity, since $V_{((r - 1)^r)^-} \cong V_{(2r - 2)} \tensor V_{((r - 1)^r)^+} \cong V_{((r - 1)^r)^+} \tensor V_{(2r - 2)}  $
\end{proposition}
\begin{proof}
    The classification of invertible simple objects in $\mathcal{C}(\mathfrak{so}_{2r}, 2r - 2)$ is given by Fuchs in \cite{Fuchs}. 

    The rule for tensoring with $V_{(2r - 2)}$ is found in Terry Gannon's paper \cite{Gannon}, section 3.4. 
\end{proof}

Now we apply these rules to perform the following computations. 

First, 
\[
V_{(2r - 2)} \tensor V_{(0)} = V_{(2r - 2)}
\]
\[
V_{(2r - 2)} \tensor V_{(1,1)} = V_{(2r - 3, 1)}
\]
\[
V_{(2r - 2)}\tensor V_{(2)} = V_{(2r - 4)}
\]
Regardless of the parity of $r$, we have
\[
V_{((r - 1^r))^+}\tensor V_{(0)} = V_{((r - 1^r))^+}
\]
\[
V_{((r - 1^r))^+}\tensor V_{(1,1)} = V_{((r - 1)^{r - 2}, 1, 1)^+}
\]
\[
V_{((r - 1^r))^+}\tensor V_{(2)} = V_{( r- 1)^{r - 1}, r - 3)^+}
\]

Next, we compute
\[
V_{((r - 1)^r)^-} \tensor V_{(0)} = V_{((r - 1)^r)^-} 
\]
\[
V_{((r - 1)^r)^-} \tensor V_{(1,1)} = V_{((r - 1)^{r - 2}, r - 2, r - 2)^-}
\]
\[
V_{((r - 1)^r)^-}\tensor V_{(2)} =V_{((r-1)^{r - 1}, r - 3)^-}
\]

Note that for $r \ge 4$, each of the Young diagrams labeling the irreducibles on the right has size greater than $2$.

\begin{lemma}
    For any proper subalgebra $A' \sbs A$, 
    \[
    \dim Hom_{C_{A'}}(F_{A'}(X)^{\tensor 2} \to F_{A'}(X)^{\tensor 2}) = 3.
    \]
\end{lemma}
\begin{proof}


    As noted above, any proper etale subalgebra $A' \sbs A$ must be pointed. As a result, a proper subalgebra must be a direct sum of the four invertible objects mentioned above. In other words, we have $A' = V_{(0)} \oplus S$ where $S$ is a direct sum of the non-trivial invertible simple objects. By the free-forgetful adjunction, we compute, 
    \[
    Hom_{C_A}(F_{A'}(X)^{\tensor 2} \to F_{A'}(X)^{\tensor 2}) \cong Hom_\mathcal{C}(V_{(1)}^{\tensor 2} \to  V_{(1)}^{\tensor 2} ) \oplus Hom_\mathcal{C}(V_{(1)}^{\tensor 2}  \to V_{(1)}^{\tensor 2} \tensor S)
    \]
    
    Since we are assuming $r \ge 4$, any of the irreducibles appearing in $V_{(1)}^{\tensor 2} \tensor S$ are labeled by Young diagrams larger than $2$. Hence $\dim Hom_\mathcal{C}(V_{(1)}^{\tensor 2}  \to V_{(1)}^{\tensor 2} \tensor S) = 0$.

    From this, fullness of $\Phi$ and $\Phi[0]$ follows by exactly the same argument as in the type B and C cases.
\end{proof}

\subsection{Level Rank Duality}

Thus far, we have represented specializations of $CW(1, -1)$ and $CW(-1, 1)$ inside quantum subgroups of type B C and D. We will now represent $CW(-1, -1)$ in quantum subgroups of type C. We will not address level rank duality for types B and D in this version of the preprint.

For a semisimple $\Z/2\Z$ braided tensor category $C$, recall that one can construct the braided tensor category $C^-$ with braiding given by 
\[
c^-_{X, Y} = \deg(X)\deg(Y)c_{X,Y}
\]
when $X$ and $Y$ are homogeneous with respect to the grading. So for example, if $X$ even and $Y$ is odd, then $C_{X,Y} = -C_{X,Y}$. The braiding can become nasty to describe directly if $X$ and $Y$ are not both homogeneous. Recall in addition that one can construct another braided tensor category $C^{rev}$, with braiding given by $c_{X, Y}^{-1}$.

It was proven in \cite[Theorem 4.1]{SymplecticLevelRank} that there exists a braid  equivalence $\mathcal{C}(\mathfrak{sp}_{2r}, k) \simeq (\mathcal{C}(\mathfrak{sp}_{2k}, r)^-)^{rev}$ (where $\mathcal{C}(\mathfrak{sp}_{2r})$ has the usual braid coming from the crossing diagram), which sends 
\[
V_{\lambda} \mapsto V_{\lambda^t}
\]
for any admissible Young diagram $\lambda$. For $k = \check{h} = r + 1$, we have a braid equivalence
\[
G: \mathcal{C}(\mathfrak{sp}_{2r}, r + 1) \to (\mathcal{C}(\mathfrak{sp}_{2r + 2}, r)^-)^{rev}.
\]

This equivalence carries the Etale algebra $A\in \mathcal{C}(\fsp_{2r}, r + 1)$ to an algebra obejct in $G(A) \in (\mathcal{C}(\mathfrak{sp}_{2r + 2}, r)^-)^{rev}$. Since the functor $G$ is a braid functor, it is clear that $G(A)$ is a connected Etale algebra in $(\mathcal{C}(\mathfrak{sp}_{2r + 2}, r)^-)^{rev}$, and hence also in $\mathcal{C}(\mathfrak{sp}_{2r + 2}, r)^-$. Now the fact that $A$ is homogeneous with respect to the grading on $\mathcal{C}(\fsp_{2r}, k)$ comes to the rescue once more. Recall from Theorem~\ref{idealsTypeC} that each irreducible subobject of $A$ is labeled by a Young diagram $\lambda$ with an even number of boxes. Hence $A\in \mathcal{C}(\fsp_{2r}, r + 1)_0$. Recall from Theorem~\cite[Theorem 4.1]{SymplecticLevelRank} that $G$ restricts to a braid equivalence $G_0: \mathcal{C}(\fsp_{2r}, r + 1)_0 \to \mathcal{C}(\fsp_{2r + 2}, r)_0^{rev}$. From this it is clear that $G(A) = G_0(A)$ is commutative with the usual braiding on $\mathcal{C}(\fsp_{2r + 2}, r)$, and hence that $G(A)$ is an Etale algebra in $\mathcal{C}(\fsp_{2r + 2}, r)$. 

It is also clear that $G(B) \in \mathcal{C}(\fsp_{2r + 2}, r)$ is a $G(A)$-module with the property that $G(B)\tensor G(B) \cong G(A)$ as $G(A)$-modules. 

Once again, in order to apply our diagrammatic techniques, we must equip $\mathcal{C}(\fsp_{2r + 2}, r)$ with the pivotal structure $\phi: Id \Rightarrow (\blank)^{**}$ such that there is a braided pivotal equivalence
\[
Ab(\sms{BMW(\xi_{4(2r + 2)}, -\xi_{4(2r + 2)}^{2(r + 1) + 1}}) \simeq \mathcal{C}(\fsp_{2r + 2}, r).
\]
Similarly, we have a pivotal equivalence (with the right choice of pivotal structure $\psi:  Id \Rightarrow (\blank)^{**}$ on $\mathcal{C}(\fsp_{2r}, r + 1)$, such that 
\[
Ab(\sms{BMW(\xi_{4(2r + 2)}, -\xi_{4(2r + 2)}^{2r + 1}}) \simeq \mathcal{C}(\fsp_{2r}, r + 1) 
\]

We must show that the functor $F_{G(A)} : \mathcal{C}(\fsp_{2r + 2}, r) \to \mathcal{C}(\fsp_{2r + 2}, r)_A$ is pivotal with this pivotal structure $\phi$. In order to show this, we must show that $\theta_{G(A)} = I_{G(A)}$, where \[
\theta_{G(A)} = (I_A \tensor ev_{G(A)^*}) \circ (I_A \tensor \phi_A^{-1} \tensor I_{G(A)^*}) \circ (c_{G(A), G(A)}\tensor I_{G(A)^*}) \circ (I_G(A) \tensor coev_{G(A)}).
\]

Note that with the pivotal structure $\psi$ on $\mathcal{C}(\fsp_{2r}, r + 1)$, we have already demonstrated that $\theta_A = I_A$. As a result, it suffices to show that the level rank duality functor $G: \mathcal{C}(\fsp_{2r}, r + 1) \to \mathcal{C}(\fsp_{2r + 2}, r)$ is pivotal. 

In order to show this, we realize the level rank duality functor by defining another functor 
\[
\Psi :BMW(\xi_{4(2r + 2)}, -\xi_{4(2r + 2)}^{2r + 1}) \to  BMW(\xi_{4(2r + 2)}, -\xi_{4(2r + 2)}^{2(r + 1) + 1})
\]
which maps 
\[
\att{braid}{0.5} \mapsto - \att{inverseBraid}{0.5}
\]
\[
\att{inverseBraid}{0.5} \mapsto -\att{braid}{0.5}
\]
and also sending the cap to the cap, and the cup to the cup.

Note that $BMW(\xi_{4(2r + 2)}, -\xi_{4(2r + 2)}^{2r + 1}) \to  BMW(\xi_{4(2r + 2)}, -\xi_{4(2r + 2)}^{2(r + 1) + 1})$ is presented as the planar algebra 
\[
\gen{\att{braid}{0.25}, \att{inverseBraid}{0.25} , \att{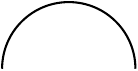}{0.25}, \att{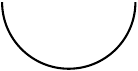}{0.25}}/R
\]
where $R$ is the tensor ideal generated by the relations imposed by the Kauffman relation, the Snake Equations, the second and third Reidemeister moves, as well as 
\[
(R1) \att{twist}{0.5} = \xi_{4(2r + 2)} \att{stick}{0.5}.
\]
All of these relations clearly hold, but (R1) does require a small computation. 

Since the map $\Psi$ satisfies the relations, it extends to a pivotal functor. By composing with the semisimplification functor, the Karoubi completion functor, and the addative completion functor, we obtain a pivotal functor 
\[
Ab(\Psi): BMW(\xi_{4(2r + 2)}, -\xi_{4(2r + 2)}^{2r + 1}) \to Ab(\sms{BMW(\xi_{4(2r + 2)}, -\xi_{4(2r + 2)}^{2(r + 1) + 1}}) 
\]

Composing equivlances, we obtain a pivotal equivalence 
\[
\tilde{G} :\mathcal{C}(\fsp_{2r}, r + 1) \to \mathcal{C}(\fsp_{2r + 2}, r)
\]

In order to show then that $G$ is pivotal, it suffices to show that $G \cong \tilde{G}$. To this end, we can exploit the classification of automorphisms of on $\mathcal{C}(\fsp_{2r}, r + 1)$. 

Composing, we have $\inv{G}\circ \tilde{G} \in Aut(\mathcal{C}(\fsp_{2r}, r + 1)) \cong \Z/2\Z$. Up to isomorphism, it follows from that $\inv{G}\circ \tilde{G}$, or $\inv{G}\circ \tilde{G} = I$, or $\inv{G} \circ \tilde{G}\cong J$, where $J$ is the simple current autoequivalence described in \cite[Section 5]{Autoequivalences} which maps $V_{(1)}\mapsto  V_{\Lambda_{r - 1} + r\Lambda_{r}}$. But note that $\inv{G}$ maps $V_{\lambda}\mapsto V_{\lambda^t}$ and in particular maps $V_{(1)} \mapsto V_{(1)} = V_{\Lambda_1}$, while by definition $\tilde{G}$ maps $V_{(1)} \mapsto V_{(1)}$. As a result, it follows that $\inv{G} \circ \tilde{G} \cong I$, hence $G \cong \tilde{G}$. Hence $G$ is a pivotal functor. 

Since the algebra $A$ lies entirely in the even part of the category $\mathcal{C}(\fsp_{2r}, r + 1)$, the algebra $G(A)$ must also lie in the even part of $\mathcal{C}(\fsp_{2r + 2}, r)$, and hence the braid of $G(A)$ with its dual is sent to the inverse braid (not the negative of the inverse braid). From this, one can see that $\theta_{G(A)} = I_A$, and that $G(A)$ is commutative with the usual braid structure on $\mathcal{C}(\fsp_{2r + 2}, r)$.

As a result, we have a pivotal functor 
\[
F_{G(A)}: \mathcal{C}(\fsp_{2r + 2}, r) \to \mathcal{C}(\fsp_{2r + 2}, r)_{G(A)}
\]

The pivotal equivalence $\mathcal{C}(\fsp_{2r}, r + 1) \to \mathcal{C}(\fsp_{2r + 2}, r)$ induces an equivalence on the categories of modules
\[
G_A:\mathcal{C}(\fsp_{2r}, r + 1)_A \to \mathcal{C}(\fsp_{2r + 2}, r)_{G(A)}
\]
such that the diagram commutes:
\[
\begin{tikzcd}
	{\mathcal{C}(\mathfrak{sp}_{2r}, r + 1)} & {\mathcal{C}(\mathfrak{sp}_{2r + 2}, r)} \\
	{\mathcal{C}(\mathfrak{sp}_{2r}, r + 1)_A} & {\mathcal{C}(\mathfrak{sp}_{2r + 2}, r)_{G(A)}}
	\arrow["G", from=1-1, to=1-2]
	\arrow["{F_A}", from=1-1, to=2-1]
	\arrow["{F_A}", from=1-2, to=2-2]
	\arrow["{G_A}", from=2-1, to=2-2]
\end{tikzcd}
\]

 Since the free module functor $F_A$ respects the $\Z/2\Z$ grading on $ \mathcal{C}(\fsp_{2r + 2}, r)$, as does the functor $G$, it follows that $G_A$ restricts to a braid equivalence on the even graded local modules, 
 $(G_A)_0: (\mathcal{C}(\fsp_{2r}, r + 1)_A^0)_0 \to (\mathcal{C}(\fsp_{2r + 2}, r)_{G(A)}^0)_0$. This sends the fermion object $g\in (\mathcal{C}(\fsp_{2r}, r + 1)_A^0)_0$ to another fermion object $(\mathcal{C}(\fsp_{2r + 2}, r)_{G(A)}^0)_0$. The image of the half braid $\gamma$ on $g$ in $\mathcal{C}(\fsp_{2r}, r + 1)_A$ is then sent to another half braid $\gamma' = G_A(\gamma)$ on $\mathcal{C}(\fsp_{2r + 2}, r)_{G(A)}$. It is not hard to use this to construct an explicit superpivotal equivalence on the fermion condensations,
 \[
 K:\mathcal{C}(\fsp_{2r}, r + 1)_A/\gamma \simeq \mathcal{C}(\fsp_{2r + 2}, r)_{G(A)}/\gamma',
 \]
 making the larger diagram commute:
 
\begin{tikzcd}
	{\mathcal{C}(\mathfrak{sp}_{2r}, r + 1)} & {\mathcal{C}(\mathfrak{sp}_{2r + 2}, r)} \\
	{\mathcal{C}(\mathfrak{sp}_{2r}, r + 1)_A} & {\mathcal{C}(\mathfrak{sp}_{2r + 2}, r)_{G(A)}} \\
	{\mathcal{C}(\mathfrak{sp}_{2r}, r + 1)_A/\gamma} & {\mathcal{C}(\mathfrak{sp}_{2r + 2}, r)_{G(A)}/\gamma'}
	\arrow["G", from=1-1, to=1-2]
	\arrow["{F_A}", from=1-1, to=2-1]
	\arrow["{F_A}", from=1-2, to=2-2]
	\arrow["{G_A}", from=2-1, to=2-2]
	\arrow[hook, from=2-1, to=3-1]
	\arrow[hook, from=2-2, to=3-2]
	\arrow["K", from=3-1, to=3-2]
\end{tikzcd}

 This equivalence sends the odd automorphism $v\in End_{\mathcal{C}(\fsp_{2r}, r + 1)_A }(F_A(V_{(1)}))$ to another odd automorphism $v' \in End_{\mathcal{C}(\fsp_{2r + 2}, r)_{G(A)}/\gamma}(F_{G(A)}(V_{(1)}))$. It also sends $F_{A}\left(-\att{inverseBraid}{0.25}\right)$ to $F_{G(A)}\left(\att{braid}{0.25}\right)$. 

 From the fact that $F_{A}\left(-\att{inverseBraid}{0.25}\right)$ and $v$ together generate the planar superalgebra generated by $F_A(V_{(1)})$ (as follows by the fullness part of Theorem~\ref{bigResultC}), it follows that their images under the equivalence $K:\mathcal{C}(\fsp_{2r}, r + 1)_A/\gamma \simeq \mathcal{C}(\fsp_{2r + 2}, r)_{G(A)}/\gamma'$ generate the planar superalgebra generated by $F_{G(A)}(V_{(1)})$. Note also that by commutativity of the above diagram, $F_{G(A)}(\att{braid}{0.25})$ is the image of $F_A(-\att{inverseBraid}{0.25})$ under the equivalence $G_A$.

That the relation (Cancel) for $CW_q(-1, -1)$ holds for $v'$ is immediate. From the fact that the functor $F_{G(A)}$ is central we have the relation (Half Braid) with the image of $\att{braid}{0.25}$.

From self duality of $F_{G(A)}(V_{(1)})$ and the the fact that $\att{twist}{0.25} = -\xi_{4(2r + 2)}^{2(r + 1) + 1}\att{stick}{0.25} = - i\xi_{4(2r + 2)}\att{stick}{0.25} $, it follows by the Twist Lemma that $v'^* = iv'$. 

All of this establishes the existence of a full functor $CW(-1, -1) \to \mathcal{C}(\mathfrak{sp}_{2r + 2}, r)_{G(A)}/\gamma'$, which completes to a superpivotal equivalence. 

This has proven the following theorem.

\begin{theorem}
    With $A'\in \mathcal{C}(\fsp_{2r + 2}, r)$ the Etale algebra coming from the conformal embedding $\V(\fsp_{2r + 2}, r) \sbs \V(\mathfrak{so}_{r(2r + 1)}, 1)$, there is an equivalence of superpivotal categories 
    \[
    Ab(\sms{CW_q(-1, -1)})\simeq \mathcal{C}(\fsp_{2r + 1}, r)_{A'}/\gamma
    \]
    which restricts to an equivalence of pivotal categories 
    \[
    Ab(\sms{CW_q(-1, -1)[0]}) \simeq \mathcal{C}(\fsp_{2r + 1}, r)_{A'}
    \]
\end{theorem}

\subsection{Level Rank Duality of Types B and D}

Simple adjustments of our methods almost certainly carry over to types $B$ and $D$ level rank duality. Nonetheless, some of the details become very annoying, so we will leave these cases out of this version of the paper.

\section{Stable Fusion Rules}\label{algebraStuff}

We now turn to the generic structure of the Clifford BMW superalgebras. Our approach is closely modeled after Wenzl's method for studying the structure of the Brauer algebras, which generalized to the BMW algebras. 

Recall that the Brauer algebras quotient onto the symmetric group algebras, and that the BMW algebras quotient on the Hecke algebras. In our case, the exact same pattern holds, only the role of the symmetric group or Hecke algebra is replaced with the Hecke Clifford superalgebras.

\subsection{Presenting The Clifford BMW superalgebras}

Recall the definition of the Clifford Hecke superalgebras. Our presentation here involves a slight change of variables to make it compatible with our description of the Clifford BMW supercategories. 

\begin{definition}
    The Hecke Clifford superalgebra $G_n(q)$ is the algebra generated by $T_1, \dots, T_{n - 1}$, $V_1, \dots, V_n$, subject to the following relations.
    \begin{itemize}
    \item[H] The $T_k$'s for $k = 1, \dots, n - 1$ satisfy the relations of the Hecke algebra $H_n(q)$,
    \begin{itemize}
        \item[H1] $T_k T_{k+ 1} T_k= T_{k + 1} T_k T_{k + 1}$
        \item[H2] $T_k^2 = (q- \inv{q}) T_k + 1$
        \item[H3] $T_k T_j = T_j T_k$ if $|k - j| > 1$
    \end{itemize}
    \item[C] The $V_\ell$'s for $\ell\in \set{1, \dots, n}$ satisfy the relations of the Clifford algebra $Cl(n)$, namely
    \begin{itemize}
        \item[C1] $V_\ell^2 = 1$ for all $\ell$, 
        \item[C2]  $V_\ell V_k = -V_k V_\ell$ for $\ell \ne k$
    \end{itemize}
\end{itemize}
    There are also two mixed relations which, geometrically, say that a $V_k$ can move under (but not over!) a strand, while the generators $V_k$ and $T_j$ commute if far enough apart.
    \begin{itemize}
        \item[M1] $V_k T_j = T_j V_k$ if $|k- j| > 1$
        \item[M2] $ V_{\ell + 1} T_\ell = T_\ell V_\ell$
    \end{itemize}
\end{definition}

In \cite{AffineSergeev}, the mixed relations given were $V_k T_k' = T_k' V_{k + 1}$. The resulting algebras are isomorphic by setting $T_k' = \inv{T_k}$, since $(V_{k + 1} T_k)^{-1} = T_k^{-1} V_{k + 1}$ and $(T_k V_k)^{-1} = V_k T_k^{-1}$. Note that this isomorphism also commutes with the natural inclusions $G_n(q) \sbs G_{n + 1}(q)$ given by $T_k \mapsto T_k$ and $V_\ell\mapsto V_\ell$ for $k \le n - 1$ and $\ell \le n$.

\begin{definition}
    Define the Clifford BMW superalgebras $CW_n(\sigma, \epsilon)$ by 
    \[
    CW_n(\sigma, \epsilon) = End_{CW(\sigma, \epsilon)}(n).
    \]
\end{definition}

\begin{theorem}
    The algebra $CW_n$ has presentation 
\[
CW_n\cong \gen{v_k, t_\ell, e_\ell,\ell \in \set{1, \dots, n - 1}, i\in \set{1, \dots, n}}/R
\]
where $R$ is the ideal generated by the following relations
\begin{itemize}
    \item[BMW] For $\ell \in \set{1, \dots, n - 1}$ the generators $t_\ell$ and $e_\ell$ satisfy the relations of the BMW algebra $BMW(q, \sigma i q^{-\epsilon})$
    \item[C] The $v_\ell$'s for $\ell\in \set{1, \dots, n}$ satisfy the relations of the Clifford algebra $Cl(n)$, namely
    \begin{itemize}
        \item[C1] $v_\ell^2 = 1$ for all $\ell$, 
        \item[C2]  $v_\ell v_k = -v_k v_\ell$ for $\ell \ne k$
    \end{itemize}
    \item[M] As well as the mixed relations
    \begin{itemize}
        \item[M1] $v_{\ell + 1} t_\ell = t_\ell v_\ell$
        \item[M2] $v_\ell e_\ell = \sigma i v_{\ell + 1}e_\ell $ and $e_\ell v_\ell = -\sigma i e_\ell v_{\ell + 1}$
        \item[M3] $v_\ell e_k = e_\ell v_k$ when $|\ell - k| > 1$
        \item[M3] $v_\ell t_k = t_k v_\ell$ when $|\ell - k| > 1$
    \end{itemize}
\end{itemize}
\end{theorem}

\begin{proof}

    By mapping 
    \[ v_k \mapsto 1^{\tensor k - 1}\tensor \att{oddDot}{0.25} \tensor 1^{\tensor n - k - 1}\], \[1^{\tensor k - 1} \tensor \att{capcup}{0.25}\tensor 1^{\tensor n - k - 2}\]
    and so forth, one easily constructs a homomorphism 
    \[
    \Psi: \gen{v_k, t_\ell, e_\ell|\ell \in \set{1, \dots, n - 1}, i\in \set{1, \dots, n}}/R \to CW_n
    \]

    From the basis theorem, it is not hard to see that this map $\Psi$ is surjective. If we can verify that $\dim \gen{v_k, t_\ell, e_\ell|\ell \in \set{1, \dots, n - 1}, i\in \set{1, \dots, n}}/R  \le 2^n (2n - 1)!!$, we will be done.

    Any word $w$ in the generators $v_k, t_\ell, e_\ell, \inv{t}_\ell$ for $k = 1,\dots, n$ and $\ell = 1, \dots, n - 1$ corresponds to a decorated tangle diagram $D(w)$. It is well known that any word in the generators $ t_\ell, e_\ell, \inv{t}_\ell $ subject to the BMW relations evaluates to a spanning set of size $(2n - 1)!!$, of words corresponding to certain standard tangle diagrams.

    For any word $w$ in the generators, $c(w)$ shall be the number of $t_k$ or $t_k^{-1}$, that is, the number of crossings in $D(w)$. 

    For any strand in $D(w)$, a spinor dot $v_k$ appearing on this strand can inductively be brought to any other point on this strand, at the cost of a scalar multiple and a linear combination of diagrams with less crossings. As a result, any two spinor dots on any closed loop in $D(w)$ can be brought close and canceled. Then, by applying the BMW relations and the relation $v_{\ell + 1} t_\ell = t_\ell v_\ell$, as well as the local commutativity and anticommutativity relations, a closed loop can be brought behind every strand up to a scalar multiple and a linear combination of words with less crossings and evaluated, either as zero if a dot remains on the closed loop after cancelation, or by picking up a multiple of $1/\delta$. 

    Having eliminated all closed loops, one can use the Dubrovnik relation as well as the relation $v_{k + 1} t_k = t_k v_k$, and so on, to bring every spinor dot to its appropriate location, thus factoring $w$, up to a linear combination of words with less crossings, as $v_S b v_T$ where $S$ and $T$ are subsets of $[n]$. We can do this in such a way that ensures that each strand has at most two spinor decorations. 

    Then apply the induction hypothesis to every accumulated term.

\end{proof}

\begin{remark}
    If one considers the subalgebra of $CW_n$ without crossings, one obtains a presentation for the endomorphism algebras for the $\Z/2\Z$-graded case of Liu and Jaffe's para-planar algebra in \cite{LiuAndJaffe}. The Clifford BMW supercategory contains this para-TLJ category. The author originally had studied this algebra in isolation, eventually stumbling upon Liu and Jaffe's work. Generically, the simple objects of the even subcategory are indexed by odd natural numbers $V_{(2n)}$ and signed even numbers $V_{(2n)^{\pm 1}}$. Choosing $q\in \C$ such that $\delta = \sqrt{2}(q + \inv{q}) = \sqrt{2}[2]_q$, one has quantum dimensions
    \[
    \dim V_{(2n)^{\pm}} = [2n + 1]_q
    \]
    and
    \[
    \dim V_{(2n + 1)} = \sqrt{2}[2n + 2]_q
    \]
    The proof of this dimension formula for para-TLJ graded over $\Z/n\Z$ for arbitrary $N$ was given in \cite{LiuAndJaffe}. The inductive argument for these dimensions, and for the branching rules, is in our opinion very beautiful.
\end{remark}

From the presentation of the Clifford BMW supercategory, the following quotient becomes clear.

\begin{corollary}
Let $I = (e_k)$ be the ideal generated by the cap cup. Then $I$ is a superideal, thus turning $CW_n/I$ into a superalgebra. 
    $CW_n/I\cong G_n$ by an even isomorphism, and hence 
    \[
    CW_n[0]/I[0] \cong G_n(q)[0]
    \]
    This isomorphism is given by even surjections
    \[
    \pi^{(n)}: CW_n \to G_n(q)
    \]
    which commute with the standard inclusions $CW_n \sbs CW_{n + 1}$ and $G_n(q) \sbs G_{n + 1}(q)$, that is, such that the diagram commutes:
\[\begin{tikzcd}
	{CW_n} & {G_n} & {CW_n[0]} & {G_n[0]} \\
	{CW_{n + 1}} & {G_{n + 1}} & {CW_{n + 1}[0]} & {G_{n + 1}[0]}
	\arrow["{\pi^{(n) }}", two heads, from=1-1, to=1-2]
	\arrow[hook, from=1-1, to=2-1]
	\arrow[hook, from=1-2, to=2-2]
	\arrow["{\pi^{(n)}_0}", from=1-3, to=1-4]
	\arrow[hook, from=1-3, to=2-3]
	\arrow[hook, from=1-4, to=2-4]
	\arrow["{\pi^{(n + 1)}}", from=2-1, to=2-2]
	\arrow["{\pi^{(n + 1)}_0}", from=2-3, to=2-4]
\end{tikzcd}\]
    
\end{corollary}

\begin{proof}

    That $I$ is a superideal follows from the fact that $e_n$ is an even element.
    
    Map $t_k\mapsto T_k$ and $v_k \mapsto V_k$, $e_k \mapsto 0$. To show that this function from the generators into $G_n(q)$ extends to a homomorphism, we need to check that all of the relations in $CW_n$ are satisfied by the images of the generators in $G_n(q)$. Noting that $BMW_n(q, \sigma i q^{-\epsilon})\sbs CW_n$ by mapping $t_k$ to $t_k$ and $e_k$ to $e_k$, and since $H_n(q) \sbs G_n(q)$ by mapping $T_k \mapsto T_k$, and since the map $t_k \mapsto T_k$ and $e_k \mapsto 0$ defines a homomorphism $BMW_n(q, \sigma i q^{-\epsilon}) \to H_n(q)$, it follows that every relation not involving $v_\ell$ for $\ell \in \set{1, \dots, n}$ is satisfied. Any relation involving $e_k$ and $v_k$ is automatically satisfied, because the corresponding words are sent to $0$. The Clifford relations $C1$ and $C2$ are satisfied by $V_1$ and $V_2$. The relations involving $T_k$ and $V_\ell$ are satisfied automatically in $G_n(q)$. 

    Since $CW_n$ is graded by the parity of the number of $v_i$'s appearing in a word in the generators, and since $G_n(q)$ is similarly graded by the number of $V_i$'s in a word in its generators, it is immediately obvious that this map is even with kernel $(e_k)$ for any $k$ (note that $e_i$ is obtained from $e_k$ by a simple application of the $t_i$'s). This map $\pi: CW_n \to G_n$ is therefore an even map with kernel $(e_k)$ for any $k$. As a result, this map descends to $\pi_0: CW_n[0]  \to G_n[0]$, with $\ker \pi_0 = (e_k) \I CW_n[0]$.

    The leftmost inclusion in the leftmost square in the statement of the corollary map $t_k\mapsto t_k$, $e_k \mapsto e_k$, $v_\ell \mapsto v_\ell$ for $k \in \set{1, \dots, n - 1}$ and $\ell \in \set{1, \dots, n}$. The rightmost inclusion maps $T_k \mapsto T_k$ for $k \in \set{1, \dots, n - 1}$ and $V_\ell \mapsto V_\ell$ for $ell \in \set{1, \dots, n}$. From this the commutativity of the rightmost square is clear. 

    Finally, the commutativity of the rightmost square in the statement of the corollary follows from the commutativity of the leftmost square, since each of the morphisms in the leftmost square are even homomorphisms of superalgebras. Taking even parts of algebras and restricting even morphisms gives a functor from the category of superalgebras and even homomorphisms to the category of algebras. 
\end{proof}

\subsection{Jones Basic Construction}

The observations of the previous subsection will help us to understanding the generic structure of the Clifford BMW algebras and superalgebras, because it will allow us to exploit the categorical trace and the theory of Jones basic constructions to inductively construct the Bratteli diagrams. First, we need to establish the generic non-degeneracy of the trace on the Clifford BMW superalgebras. 

\begin{lemma}\label{TraceNonDegenerate}
    For all $n \in \N$ and for $q\in \C$ with $q^8 \ne 1, 0$, the categorical trace on $CW_n$ is non-degenerate if and only if the categorical trace on $ CW_n[0]$ is non-degenerate. Moreover, the categorical trace is non-degenerate on both algebras for all but finitely many $q\in \C$.
\end{lemma}

\begin{proof}
    Suppose first that the categorical trace on $CW_n[0]$ is non-degenerate, and let $f\in CW_n$ be non-zero. Then either $f_0 \ne 0$ or $f_1 \ne 0$. 
    
    First suppose $f_0 \ne 0$. Then $f_0 \in CW_n[0]$. Since $\sms{CW_n} = CW_n$ from the fact that the trace is non-degenerate, $f_0\in CW_n[0]$ is not negligible, hence there exists some $g\in CW_n[0]$ such that $tr(f g) \ne 0$. But note that $g = g_0$, and hence $f_1 g_0$ is odd. Since odd morphisms have zero trace, we have $tr(f_1 g) = 0$, hence
    \[
    tr(fg) = tr(f_0 g) + tr(f_1 g) = tr(f_0 g) \ne 0.
    \]

    On the other hand, suppose that $f_0 = 0$, so that $f$ is odd. Then $v_0f$ is even, and non-zero. 

    We now prove the converse. Suppose that the categorical trace on $CW_n$ is non-degenerate. Suppose $f = f_0 \in CW_n[0]$ is non-zero. Then there exists some $g = g_0 + g_1\in CW_n$ such that $tr(fg) \ne 0$. Note that $f_0 g_1$ is odd, and hence $tr(f_0 g_1) = 0$. Hence
    \[
    tr(f_0 g_0) = tr(fg) - tr(f_0 g_1) = tr(fg) \ne 0.
    \]
\end{proof}

\begin{remark}
    Semisimplicity is preserved by fermion condensation in general. If $C$ is a braided semisimple tensor category, let $g\in C$ be a fermion object. Then any simple object $X\in C$ is a type $Q$ simple in $C/\gamma$ if $X \tensor g \cong X$ in $C$, and a type $M$ simple otherwise. As a result, since objects in the $C/\gamma$ are objects in $C$, we have $X\cong \bigoplus_i n_i X_i$ with $X_i$ finitely many simple objects. Each of these is either a type $Q$ or type $M$ simple. It is not hard to see from this that this forces $End_{C/\gamma}(X)$ to be a direct sum of matrix algebras either over $\C$, or over $Cl_1$.    
\end{remark}

\begin{corollary}
    $CW_n(\sigma, \epsilon)$ and $CW_n(\sigma, \epsilon)[0]$ are semisimple with non-degenerate categorical trace for all but finitely many values of $q\in \C$. 
\end{corollary}
\begin{proof}
    The categorical trace on $CW_n(-1, 1)$ is non-degenerate when $q = e^{\frac{2\pi i}{8r + 8}}$ and $n < r/2$, since $End_{C_A}(F(V_{(1)})^{\tensor n}) \cong CW_n(-1, 1)$ for these values of $q$. Hence the Gram matrix for the basis elements for both $CW_n[0]$ and $CW_n$ have non-zero determinate. Since these are rational functions in $q$, it follows in particular that these determinates are not zero over the ground ring $R$. Hence for all but finitely many specializations, both algebras have a non-degenerate trace, and hence both are semisimple for all but finitely many specializations. 

    This establishes that the categorical trace trace $tr:CW_n(-1, 1) \to R$ is non-degenerate before specialization, with the determinant of the Gram matrix $A_{-1, 1}(q)$ some non-zero rational function $f_{-1, 1}(q)\det A_{-1, 1}(q) \in R$. Since the permutations $q\mapsto -1/q$ and $i \mapsto -i$ both induce automorphisms of the ring $R$, it follows that the Gram matrix $A_{\sigma, \epsilon}(q)$ has determinant $f_{\sigma, \epsilon}(q)$ non-zero in $R$ for all sign choices $\sigma, \epsilon \in \set{\pm 1}. $
\end{proof}

Now we briefly review the Jones basic construction.

Suppose algebras $A$ and $B$ are semisimple algebras over $\C$, where $A \sbs B$. Suppose there exists a non-degenerate trace $\tau: B \to K$ such that the restriction of $\tau$ to $A$ is non-degenerate. Then a non-degenerate symmetric bilinear form on both $B$ and $A$ is defined by $\inner{x}{y} = \tau(xy)$ for $x, y\in B$. Since this form is non-degenerate and symmetric, there exists a unique orthogonal projection $\epsilon_A : B \to A$, called the \textbf{conditional expectation.}

Since $B$ is split semisimple, it is a direct sum of matrix algebras over $\C$, so that matrix transposition gives an isomorphism $J:B \to B^{op}$. This defines a right action of $\rho: B\hookrightarrow End_\C(B)$, thus embedding $B$ into $End_\C(B)$. The subalgebra $A$ acts on $B$ on the left by multiplication, giving an embedding $\lambda: A \hookrightarrow End_\C(B)$.

\begin{definition}
    The \textbf{Jones basic construction} $J$ for the extension $A \sbs B$ is the centralizer of $\lambda(A)\sbs  End_\C(B)$, $J = End_\C(B)^{\lambda(A)}$
\end{definition}

This centralizer has an explicit description in terms of the conditional expectation $\epsilon_A$. Note that $\rho(B)$ commutes with $\lambda(A)$.
As a result, the embedding $\rho: B \hookrightarrow J$ gives a three leveled tower of algebras
\[
A \sbs B \hookrightarrow J.
\]

Let $\set{V_\lambda}_{\lambda \in \Lambda^A}$ be all the non-isomorphic irreducible representations of $A$, and $\set{W_\mu}_{\mu \in \Lambda^B}$ be a set of non-isomorphic irreducible representations of $B$. Then for each $\mu \in \Lambda^B$ we have 
\[
Res (W_\mu)\cong  \bigoplus_{\lambda \in \Lambda^A} N_{\mu \lambda} V_\mu.
\]

\begin{theorem}[Wenzl, Jones]
    The algebra $J$ is semisimple. Let $p_1, \dots, p_k$ be a partition of unity in $A$. Then $p_1 \epsilon_A, \dots, p_k \epsilon_A$ is a partition of unity in $J$, so that the irreducible representations of $J$ are in bijection with those of $A$. Letting $\set{U_\lambda}_{\lambda \in \Lambda^A}$ be the irreducible representations, the restriction rules to $B$ are given by 
    \[
    Res(U_\lambda) \cong \bigoplus_{\mu \in \Lambda^B} N_{\lambda \mu V_\lambda}. 
    \]
    
    Moreover, $J = \rho(B) \epsilon_A \rho(B) = \gen{\rho(B), \epsilon_A}$.
\end{theorem}

In other words, the Bratelli diagram for the inclusion $B \hookrightarrow J$ is the transpose of the Bratelli diagram for $A \sbs B$. 

Now we recall some facts about irreducible modules over a graded algebra. If $A = A_0 \oplus A_1$ is a superalgebra, a supermodule is a superspace $V = V_0 \oplus V_1$ along with an action of $A$ on $M$ such that $A_i.M_j \sbs M_{i + j}$. If $A$ is semisimple, all of its modules are semisimple, including the graded ones. But it need not be the case that this decomposition respects the grading. This gives rise to two notions of irreducible modules: type M and type Q. A module $M$ is type $M$ if $End_{A}(V) \cong \C$, and type $Q$ if $End_A(V) \cong Cl_1$.

Now we will show that the classic result of the Jones basic construction can be ``upgraded'' to a super setting. By an \textbf{even trace} on a superalgebra $A$, we mean a trace $\tau: A \to \C$ such that $\tau(a) = 0$ if $a$ is odd.

\begin{lemma}\label{GradedJBM}
    Let $A \sbs B$ be semisimple superalgebras with a compatible even trace. Then $\epsilon_A: B \to A$ is even and the Jones basic construction $B\epsilon_A B$ contains $B$ as a superalgebra. If $\set{p_\lambda}$ is a set of orthogonal minimal even idempotents of $A$ which sum to $1$, then so is $\epsilon_A p_\lambda$ in $A$. We also have an isomorphism of super spaces for every minimal even idempotent $q_\mu$ of $B$,
    \[
    Hom_{A}( Bq_\mu \to Ap_\lambda) \cong Hom_{B\epsilon_A B} (B\epsilon_A B\epsilon_A p_\lambda \to Bq_\mu)
    \]
\end{lemma}

\begin{proof}
    The inner product $\inner{x}{y} = \tau(xy)$ makes $B_0$ and $B_1$ orthogonal since the trace $\tau$ is even. As a result, since $\epsilon_A: B \to A$ is the orthogonal projection onto $A$, it follows that $\epsilon_A(B_i) \sbs A_i$ for $i \in \Z/2\Z$. Hence $\epsilon_A$ is even. As a result, the Jones basic construction $J$, as a graded subalgebra of $End_\C(B)$, has $\epsilon_A$ even and $B \sbs J$ evenly.\footnote{The grading on $J$ is induced from the obvious grading on $End_\C(B)$, which in turn is induced from the grading on $B$. A linear map $f: B \to B$ has grading $|f|\in \Z/2\Z$ if $f(B_i) \sbs B_{i + |f|}$.}

    Since $\epsilon_A$ is even, if $p_\lambda\in A$ is a type $Q$ minimal even idempotent, then $p_\lambda = x_\lambda  + y_\lambda$ for $x_\lambda$ and $y_\lambda$ inhomogeneous minimal idempotents. Since $\epsilon_A$ is even, $\epsilon_Ax_\lambda$ and $\epsilon_A y_\lambda$ are also inhomogeneous, and hence $\epsilon_A p_\lambda$ is a type $Q$ minimal even idempotent in $B\epsilon_A B$. On the other hand, if $p_\lambda$ is a type $M$ minimal even idempotent, it is a minimal idempotent, and $\epsilon_A p_\lambda$ is both minimal and even, and hence type $M$. 

    We have the following isomorphism of super spaces, using the fact that $\epsilon_A A = \epsilon_A B \epsilon_A$, as well as the fact that multiplying $A$ by $\epsilon_A$ gives an even isomorphism $A\mapsto \epsilon_A A$, as well as the fact that $\epsilon_A$ is an idempotent,
    \begin{align*}
        Hom_{A}( Bq_\mu \to Ap_\lambda)\cong q_\mu B p_\lambda = q_\mu B A p_\lambda \cong q_\mu B A \epsilon_A p_\lambda \\= q_\mu B \epsilon_A B \epsilon_A p_\lambda = Hom_{B\epsilon_A B} ( Bq_\mu \to B\epsilon_A B\epsilon_A p_\lambda)
    \end{align*}
\end{proof}

If one has an even inclusion of semisimple superalgebras $A \sbs B$, any simple  supermodule $V_\lambda$ of $B$ has $Hom_A(V_\lambda \to W_\mu) = \C^{M_{\lambda \mu}|N_{\lambda \mu}}$. So, drawing a two colored graph and coloring the simple supermodules of both $A$ and $B$ as type $Q$'s or type $M$'s, one has a \textbf{super-Bratelli diagram}.

The above lemma shows that, just as the ordinary Bratelli diagram for $B \sbs J$ is a reflection of the Bratelli diagram for $A \sbs B$, if $A \sbs B$ are both super algebras with an even compatible trace, then $J$ is also a super algebra with the super-Bratelli diagram for $B \sbs J$ a reflection of the super-Bratelli diagram for $A \sbs B$.

\begin{remark}
    This clearly works just as well for arbitrary gradings by arbitrary abelian groups. There, one would have a type $Q_H$ for every subgroup. The same notion of a ``super-Bratelli diagram'' would generalize automatically.
\end{remark}

The Jones basic construction $J$ for a tower (graded or not) $A \sbs B$ with a compatible trace is also characterized by a ``universal property.''

\begin{theorem}[Jones, Wenzl]
    Suppose $C$ is a semisimple algebra over $\C$ such that $B\sbs C$, and such that there is an element $e\in C$ such that 
    \begin{enumerate}
    \item $e^2= e$
    \item $ e b e = \epsilon_A(b) e = e \epsilon_A(b) $
    \item the map $a \mapsto ae$ is an injective algebra homomorphism $A\to C$.
    \end{enumerate}
    Then $J \cong BeB$ and $C \cong J \oplus \tilde{C}$ for some subalgebra $\tilde{C}\sbs C$.  
\end{theorem}

In many of the most interesting examples of towers of algebras $A_0 \sbs A_1 \sbs \dots$ with an element $e_n\in A_{n+1}$ satisfying the above properties for the inclusion $A_{n - 1} \sbs A_{n}$, we have the identity 
\[A_n e_{n} A_n= A_{n + 1}e_n A_{n + 1}\]

This is the case for the TLJ algebras, the Brauer algebras, the BMW algebras, and many others. It is also the case for the Clifford BMW superalgebras.

We will need this well known property of the BMW algebras:

\begin{lemma}\label{BMWIdeal}
    For all $n\in \N$, the identity holds:
    \[
    BMW_{n + 1} e_n BMW_{n + 1} = BMW_n e_n BMW_n.
    \]
    In particular, 
    \[
    BMW_{n + 1} e_n = BMW_ne_n
    \]
    and 
    \[
    e_n BMW_{n + 1} = e_n BMW_n.
    \]
\end{lemma}

From now on, let us denote top boundary points in $CW_n$ or $BMW_n$ by barred integers $[\sms{n}] = \set{\sms{1}, \dots, \sms{n}}$ and bottom boundary points by unbarred integers $\set{1, \dots, n}$. For a BMW or Clifford BMW diagram $d$ in $CW_n$, and for any $\mathbf{i}\in [n] \U [\sms{n}]$ let $\Sigma_d(\mathbf{i})$ denote the strand connecting to the $\mathbf{i}$ boundary point. For any $\mathbf{i}, \mathbf{j}\in [n] \U [\sms{n}]$ write $\mathbf{i} \sim_d \mathbf{j}$ if $\Sigma_d(\mathbf{i}) = \Sigma_d(\mathbf{j})$. 

We will need a slight strengthening of Theorem~\ref{smallSpanningSet}. 

\begin{lemma}\label{basisChangeLemma}

    Let $\mathcal{D}$ be any basis of diagrams for $BMW_n$ such that every diagram in $\mathcal{D}$ has the property that no strand crosses another strand twice, and no strand crosses itself. For each diagram $d\in \mathcal{D}$, pick any set $\Omega_d \sbs [n]\U [\sms{n}]$ such that $\mathbf{i}\ne \mathbf{j}\in \Omega_d$ implies that $\Sigma_d(\mathbf{i}) \ne \Sigma_d(\mathbf{j})$.
    
    There is a basis $\mathcal{B}_{\mathcal{D}, \Omega}$ of $CW_n$ of elements of the form 
    \[
    c_T d c_S
    \]
    where $T = [\mathbf{n}]\I \Omega_d$ and $S = [n]\I \Omega_d$
 
\end{lemma}

\begin{proof}
    The proof is by induction on the number of crossings, and is identical to the proof of Theorem~\ref{smallSpanningSet}. Spinors can be freely moved around on the strands, at the cost of terms with less crossings and scalar multiples coming from super-isotopy. Any basis $\mathcal{D}$ for $BMW_n$ consisting of diagrams which have no self crossings or double crossings must also include the standard basis of Tiemperley Lieb Jones diagrams. As a result, we can exchange any crossings at the cost of terms with less crossings, eventually terminating at elements which are necessarily in the basis. 
\end{proof}

Another way of stating this lemma is that we have complete freedom to choose whichever diagram basis of $BMW$ we like (so long as there are no self crossings or double crossings) as well as a choice of which boundary point of each strand is allowed to have spinors. Every choice yields a basis for $CW_n$.

\begin{corollary}\label{superDuperSpecialChoice}
    For every diagram $d\in BMW_n$, at least one $i \in \set{n - 1, n}$ has $i \not\sim_d \sms{n}$. For every BMW diagram $d\in BMW_n$ let $i_d$ be a chosen element $i_d \in \set{n - 1, n}$ with this property. 

    For every $n\in \N$ there is a basis $\mathcal{D}$ of diagrams in $BMW_n$ such that for every diagram $d\in \mathcal{D}$, $\Sigma_d(i_d)$ is beneath every other strand of $d$. Then there is a basis $\mathcal{B}_{\mathcal{D}, \Omega}$ of the form described in the previous lemma, such that for all $d\in \mathcal{D}$, $\sms{n}\not \in \Omega_d$. 
\end{corollary}
\begin{proof}
    The first part of the corollary comes from the obvious fact that no strand can connect to three boundary points. A basis can then be chosen for $BMW_n$ with the prescribed property of the corollary, because every choice of over-under crossings for any Brauer diagram yields a basis for $BMW_n$. For each of the diagrams $d\in \mathcal{D}$, the strand $\Sigma_d(\sms{n})$ connects to another boundary point not equal to $\sms{n}$, so we can choose $\Omega_d$ so that $\sms{n}\not\in \Omega_d$. 
\end{proof}

\begin{lemma}
    Identifying $CW_n[0]$ as a subalgebra of $CW_{n + 1}[0]$, and similarly for $CW_n$, we have the equalities for all $n$,
    \begin{itemize}
        \item $CW_n e_n CW_n = CW_{n + 1} e_n CW_{n + 1}$
        \item $CW_n[0] e_n CW_n[0] = CW_{n + 1}[0] e_n CW_{n + 1}[0]$ when $n > 1$. 
    \end{itemize}
    
\end{lemma}

\begin{proof}

First, we will prove $CW_{n + 1}e_n = CW_n e_n$. Clearly $CW_ne_n \sbs CW_{n + 1}e_n$, so we need only prove the reverse inclusion. By the previous corollary, this reduces to showing that every diagram of the form
\[
\att{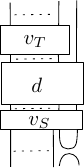}{1}
\]
where the diagram $d\in BMW_n$ satisfies the property that one of the bottom two boundary points connects to a strand which does not connect to the top left boundary point, and such that that strand crosses under every other strand in $d$. 

Either both $n, n + 1 \in S$, only one of $n, n + 1\in S$, or neither $n, n + 1\in S$. 

If both $n, n + 1\in S$, then the spinors can slide under the cap at the cost of scalar multiples and cancel. Thus the case where both $n, n + 1\in S$ reduces to the case where neither $n, n + 1\in S$. 

On the other hand, if $n + 1\in S$ but $n\not\in S$, the spinor dot at the $n + 1$ boundary point can slide along the cup at the cost of a scalar to the $n$ boundary point. Hence we really only have two cases: either $n \in S$ and $n + 1\not\in S$, or neither $n, n + 1 \in S$. 

In the case where neither $n, n + 1\in S$, we apply Lemma~\ref{BMWIdeal} to get
\[
\att{ideal1}{1} = \att{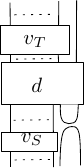}{1} = \att{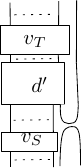}{1}\in CW_n e_n
\]
where $d'\in BMW_n$ is a linear combination of diagrams. 

On the other hand, suppose $n\in S$ but $n + 1\not\in S$. In the case where the chosen strand connecting to one of $n, n + 1$ on the bottom in $d$ which does not connect to $\sms{n + 1}$ connects to $\sms{1}$, and since the spinor elements in $v_S$ are ordered left to right:
\[
\att{ideal1}{1} = \att{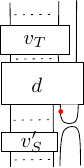}{1} = \att{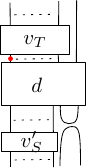}{1} = \att{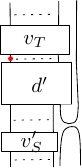}{1} \in CW_n e_n.
\]
where the relations (Half Braid) and (Rotate) were used repeatedly to bring the spinor along the strand $\Sigma_d(n)$ to the top boundary. It is true, but not necessary to our purposes, that the accumulated scalar from these relations is $1$. 

All of the other cases are similar. 

This establishes that $CW_{n + 1} e_n \sbs CW_n e_n$. Reflecting diagrams accross a horizontal axis defines a $\C(q)$ conjugate linear (mapping $i \to -i$ and fixing $q$) anti-automoprhism $\dag:CW_n \to CW_n$. From this we have $e_n CW_{n + 1} = e_n CW_n$. 

As a result, $CW_{n + 1} e_nCW_{n + 1} = CW_{n}e_n CW_n$.

Now we turn to the even part. By definition, we have 
\[
CW_{n + 1}[0] e_n CW_{n + 1}[0] \sbs (CW_{n} e_n CW_{n})[0]. 
\]
As a result, we need to show that $(CW_n e_n CW_n)[0] \sbs CW_n[0] e_n CW_n[0]$. 

Suppose we have $x e_n y$ for $x,y\in CW_n$ Clifford BMW diagrams. If $|x| = 0$ and $|y| = 0$, we already have $x e_n y\in CW_n[0] e_n CW_n[0]$. Otherwise, we have $|x| = 1$ and $|y| = 1$. Assuming then that $n > 2$, we have by local commutativity of $e_n$ and $v_1$ that
\[
x e_n y = xv_1^2 e_n y = x v_1 e_n v_1 y.
\]
Now notice that $|x v_1| = |x| + 1 = 0$ and similarly for $y$. Hence 
\[
x e_n y = x v_1 e_n v_1 y \in CW_n[0] e_n CW_n[0].
\]
\end{proof}
\begin{remark}
    The identity $(CW_n e_n CW_n)[0] \sbs CW_n[0] e_n CW_n[0]$ is not true when $ n = 1$. Consider $v_1 e_n v_1 \in (CW_n e_n CW_n)[0]$ which is not in $CW_n[0] e_n CW_n[0]$. 
\end{remark}

The previous lemma, as well as Corollary~\ref{superDuperSpecialChoice}, allow us to essentially ``steal'' all of the representation theory from the Hecke Clifford superalgebras and apply it to the Clifford BMW superalgebras. 

The representation theory for the Hecke Clifford algebras is well understood generically, with recent work establishing some results about its behavior at roots of unity. If $q$ is not a root of unity, then the Hecke Clifford algebras $G_n(q)$ are semisimple for all $n$. On the other hand, if $q$ is a root of unity, then $G_n(q)$ may not be semisimple. Nonetheless, if $q$ is a $4N$-th root of unity and $n < N$, it has been established by Wenzl and Edie-Michell that $G_n(q)$ and $G_n(q)[0]$ are in fact semisimple. 
The even part of $G_n(q)$, $G_n(q)[0]$, is the endomorphism algebra of Edie-Michell and Snyder's interpolation category $\mathcal{E}_q$. In order to understand the simple objects and fusion rules in $\mathcal{C}(\mathfrak{sl}_{r}, r + 1)_A$ where $A$ is the Etale algebra coming from the conformal embedding $\V(\mathfrak{sl}_{r}, r + 1) \sbs \V(\so_{\dim \mathfrak{sl}_{r}}, 1)$, Edie-Michell and Wenzl in \cite{CainAndHans} extended the work of Jones and Nazarov in \cite{AffineSergeev} to classify irreducible representations of the even subalgebra $G_n(q)[0]$. They also gave a partial description of the specializations of both $G_n(q)$ and $G_n(q)[0]$ to roots of unity.

We use the representations of both $G_n(q)$ and $G_n(q)[0]$. The even subalgebras have a slightly less clean description. 

\begin{definition}
    The Schur graph $\Sigma$ is the graph whose nodes are strict partitions, with an edge $ \lambda \to \mu$ if $\mu$ is obtained from $\lambda$ by the addition of a single box. The Schur graph is a leveled graph, with $\Sigma_n$ consisting of all strict partitions of size $n$. Write $\lambda \leftrightarrow \mu$ when either $\lambda \to \mu$ or $\mu \to \lambda$. 
\end{definition}

\begin{theorem}[\cite{CainAndHans}]
    For all but finitely many $q\in \C$ the algebra $G_n(q)[0]$ is semisimple. Whenever $G_n(q)[0]$ is semisimple, the irreducible representations of $W_{\lambda, \kappa}$ are indexed by strict partitions $\lambda \dashv_s n$ along with a sign choice $\kappa \in \set{1, (-1)^{\ell(\lambda) + 1}}$. For all but finitely many $q\in \C$, restriction is given by 
    \[
    Hom_{G_{n - 1}[0](q)} (W_{\lambda, \kappa}^{(n)} \to W_{\mu, \pi}^{(n - 1)}) \cong \begin{cases}
        \C^2 & \mu \to \lambda, \ell(\mu)\ell(\lambda) \mod 2 = 1\\
        \C & \mu \to \lambda, \ell(\mu)\ell(\lambda) \mod 2 = 0\\
        0 & \mu \not\to \lambda
    \end{cases}
    \]
\end{theorem}

\begin{theorem}[\cite{AffineSergeev}, \cite{CainAndHans}]\label{repsOfGn}
    For all but finitely many $q\in \C$, the algebra $G_n(q)$ is semisimple with non-isomorphic graded irreducible representations $V_\lambda$ indexed by strict partitions $\lambda \dashv_s n$, where $V_\lambda$ is type $Q$ if $\ell(\lambda)$ is odd and type $M$ if $\ell(\lambda)$ is even. Moreover, for all $\mu\dashv_s n - 1$,
    \[
    Hom_{G_{n -1}}(V_\lambda^{(n)} \to V_\mu^{n - 1}) \cong \begin{cases}
        \C^{2|2} & \mu \to \lambda, \ell(\mu)\ell(\lambda) \mod 2 = 1\\
        \C^{1 | 1} & \mu \to \lambda, \ell(\mu)\ell(\lambda) \mod 2 = 0\\
        0 & \mu \not\to \lambda
    \end{cases}
    \]

    In particular $G_n$ is semisimple when $q$ is a primitive $4N$-th root of unity and $n \le N$ or when $q$ is not a root of unity. The above fusion rule holds when $q$ is not a root of unity, as well as whenever $ n \le N$ if $q$ is a primitive $4N$-th root of unity.

    We also have the restriction rules
    \[
    Res_{G_n(q)[0]}^{G_n(q)}(V_\lambda) = \begin{cases}
        W_\lambda \oplus W_\lambda & \ell(\lambda) \equiv 1\pmod{2}\\
        W_{\lambda, -1} \oplus W_{\lambda, 1} & \ell(\lambda) \equiv 0 \pmod{2}
    \end{cases}
    \]
\end{theorem}
\begin{proof}
    In \cite{AffineSergeev} the type $Q$ supermodules $V_\lambda$ for $\ell(\lambda)$ odd were viewed as a direct ordinary modules of $G_n(q)$, $V_{\lambda, 1}$ and $V_{\lambda,-1}$. So, as an ungraded module, when $\ell(\lambda)$ is odd, $V_\lambda = V_{\lambda, 1} \oplus V_{\lambda,-1}$. 
    
    For all but finitely many $q\in \C$, it was shown in \cite{CainAndHans} that (as ordinary representations rather than graded ones) $Hom_{G_n(q)}(V_{\lambda, \pm 1}^{(n)} \to V_{\mu,(-1)^a}^{(n - 1)}) \cong \begin{cases}
        \C & \mu \to \lambda \\
        0 & \mu \not \to \lambda
    \end{cases}$ regardless of the parity of $\ell(\mu)$. As a result, if $\ell(\lambda)\ell(\mu) \equiv 1 \pmod{2}$ we have $\dim Hom_{G_{n-1}(q)}(V_\lambda^{(n)} \to V_{\mu}^{(n - 1)}) = 4$. Similarly, if $\ell(\lambda)\ell(\mu) \equiv 0\pmod{2}$ then $\dim Hom_{G_{n-1}(q)}(V_\lambda^{(n)} \to V_{\mu}^{(n - 1)}) = 2$. Regardless, since for all but finitely many $q\in \C$ we have that $G_n(q)$ and $G_{n - 1}(q)$ are semisimple, it follows from semisimplicity that
    \[
    Hom_{G_{n-1}(q)}(V_\lambda^{(n)} \to V_{\mu}^{(n - 1)}) \cong p_\lambda G_n(q) p_\mu
    \]
    as supserspaces, where $p_\lambda$ and $p_\mu$ are the minimal even idempotents in $G_n(q)$ and $G_{n - 1}(q)$ respectively such that $V_\lambda \cong G_n(q) p_\lambda$ and $V_\mu \cong G_{n - 1}(q) p_\mu$. Now note that $\dim G_n(q)[0]\cong \dim G_n(q)[1]$. So, since the idempotents $p_\lambda$ and $p_\mu$ are even, it follows that the even and odd dimensions of $Hom_{G_{n-1}(q)}(V_\lambda^{(n)} \to V_{\mu}^{(n - 1)})$ are equal. As a result, we have the graded restriction rule stated in the theorem. 

    The restriction rule from $G_n(q)[0]$ to $G_n(q)$ was given in \cite{CainAndHans} Lemma 3.16.
\end{proof}

Using this, we can prove the following theorem.

\begin{theorem}\label{RepsOfCW}
    For all but finitely many $q\in \C$, the graded irreducible representations of $CW_n$ are indexed by $\lambda \dashv_s n - 2f$ for $f = 0, \dots, \lfloor n/2 \rfloor$ where $V_\lambda^{(n)}$ is a type $M$ irreducible when $\ell(\lambda)$ is even, and where $V_\lambda^{(n)}$ is a type $Q$ irreducible if $\ell(\lambda)$ is odd, and whenever $n > 1$ we have 
    \[
    Hom_{CW_n}(V_\lambda^{(n)} \to V_\mu^{(n - 1)}) \cong \begin{cases}
        \C^{2|2} & \mu \leftrightarrow \lambda, \ell(\mu)\ell(\lambda) \mod 2 = 1\\
        \C^{1 | 1} & \mu \leftrightarrow \lambda, \ell(\mu)\ell(\lambda) \mod 2 = 0\\
        0 & \mu \not\leftrightarrow \lambda
    \end{cases}
    \]

    Similarly, for all but finitely many $q\in \C$, the irreducible representations $W_{\lambda, \kappa}$ of $CW_n$ are indexed by $\lambda \dashv_s n - 2f$ for $f = 0, \dots, \lfloor n/2 \rfloor$ along with a sign choice $\kappa \in \set{1, (-1)^{\ell(\lambda) + 1}}$, such that for all $\lambda \dashv_s n - 2f$ and $\mu\dashv_s n -1 - 2f $ we have 
    \[
    Hom_{G_{n - 1}[0](q)} (W_{\lambda, \kappa}^{(n)} \to W_{\mu, \pi}^{(n - 1)}) \cong \begin{cases}
        \C^2 & \mu \leftrightarrow \lambda, \ell(\mu)\ell(\lambda) \mod 2 = 1\\
        \C & \mu \leftrightarrow \lambda, \ell(\mu)\ell(\lambda) \mod 2 = 0\\
        0 & \mu \not\leftrightarrow \lambda
    \end{cases}
    \]

    Restriction rules for the inclusion $CW_n[0] \sbs CW_n$ are given by the same formula as for the inclusion $G_n(q)[0] \sbs G_n(q)$:
    \[
    Res_{CW_n[0]}^{CW_n}(V_\lambda) = \begin{cases}
        W_\lambda \oplus W_\lambda & \ell(\lambda) \equiv 1\pmod{2}\\
        W_{\lambda, -1} \oplus W_{\lambda, 1} & \ell(\lambda) \equiv 0 \pmod{2}
    \end{cases}
    \]
    for every $\lambda \dashv_s n - 2f$.

\end{theorem}
\begin{proof}
    Assume $q$ is such that $\delta\ne 0$. The categorical trace $\tau$ can be normalized into a compatible trace for the tower $\set{CW_n}_{n = 0}^\infty$ by 
    \[
    tr_{n}(x) = \tau|_{CW_n}(x)/\delta^n.
    \]
    Compatibility of this trace is clear from drawing the diagrams. Moreover, the Gram matrix for this trace is a non-zero rational function in $q$ generically by Lemma~\ref{TraceNonDegenerate}. As a result, for every integer $N\in \N$, $\tau_k$ is non-degenerate for each $k \le N$ at all but finitely many values of $q$. Unioning the zeros of each of the Gram matrices at levels up to $N$ gives a finite set of elements $q\in \C$ for which at least one of the traces 
    \[
    tr_0, tr_1, \dots, tr_N
    \]
    becomes degenerate. Appart from these values, $tr_k$ is non-degenerate on $CW_k$ and for $CW_k$ for all $k \le N$.
    
    The proof is by induction on $n$. For $n = 0$, $CW_n \cong \C$ which has a single irreducible representation $V_{(0)}^{(0)} = \C^{1 | 0}$. For $n = 1$, $CW_n^* \cong Cl_1$ which has a single irreducible representation of type $Q$ $V_{(1)}^{(1)}$. In the karoubi completion, the minimal even idempotent $P_{(1)}\in CW_n^*$ with $ImP_{(1)} = V_{(1)}^{(1)}$ is a type $Q$ simple, having an odd automorphism $v: P_{(1)}\cong P_{(1)}$. 

    Now suppose the theorem holds for all $k \le N$, and suppose $q$ is such that the trace $tr_k$ is non-degenerate for all $k \le N$. Then $CW_{N - 2}^*$ and $CW_{N - 1}$ split into direct sums of matrix algebras over $\C$. 

    Since we are assuming $tr_{N + 1}$ is non-degenerate, $CW_{N + 1}^*$ is semisimple. By semisimplicity, we have a splitting 
    \[
    CW_{N + 1} = CW_{N + 1} e_n CW_{N + 1} \oplus G_{N + 1}(q).
    \]

    The first summand is the Jones basic construction for the inclusion $CW_{N - 1}^* \sbs CW_{N}$ since $CW_{N + 1} e_n CW_{N + 1}^*$. For convenience let $J =  CW_{N + 1} e_n CW_{N + 1} $. By the induction hypothesis, the irreducible super-representations $V_\lambda^{(N - 1)}$ of $CW_{N - 1}$ are indexed by strict partitions $\lambda\dashv_s N - 1 - 2f$, where $V_\lambda^{(N - 1)}$ is type $Q$ if $\ell(\lambda)$ is odd and type $M$ otherwise. Moreover, the strict partitions $\mu \dashv_s N - 2f$ index the irreducible super-representations of $CW_N$, again by induction. We also have, again by the induction hypothesis, that \[Hom_{CW_{N - 1}}(V_\lambda^{(N - 1)} \to V_\mu^{N})
\cong \begin{cases}
        \C^{2|2} & \mu \leftrightarrow \lambda, \ell(\mu)\ell(\lambda) \mod 2 = 1\\
        \C^{1 | 1} & \mu \leftrightarrow \lambda, \ell(\mu)\ell(\lambda) \mod 2 = 0\\
        0 & \mu \not\leftrightarrow \lambda
    \end{cases}
    \] for any $\mu \dashv_s N - 2f$. By Lemma~\ref{GradedJBM} the irredicible super-representations of $J \sbs CW_{N+ 1}$ are in bijection with those of $CW_{N - 1}$ with type $Q$'s and type $M$'s preserved, and we have 
    \[
    Hom_{CW_N}(V_\lambda^{(N + 1)} \to V_\mu^{N}) \cong Hom_{CW_{N - 1}}(V_\lambda^{(N - 1)} \to V_\mu^{N}).
    \]

    We also have by Lemma~\ref{repsOfGn} that, for all but finitely many $q$ the irreducible representations of $G_{N + 1}(q)$ are indexed strict partitions $\lambda \dashv_s N + 1$. Since $G_{N + 1}(q)$ appears as a direct sum of $CW_{N + 1}$, the projection onto $G_{N + 1}(q)$ turns each of these into irreducible super-representations of $CW_{N + 1}$, obviously with the same types. Since $G_N(q)$ appears as a direct summand of $CW_N$ in a way which commutes with the inclusion $CW_N \sbs CW_{N + 1}$, it follows that the same restriction rules hold, namely, if $\lambda \dashv_s N + 1$ and $\mu \dashv_s N$, we have 
    \[
    Hom_{CW_N}(V^{(N + 1}_\lambda \to V_\mu^{(N)}) \begin{cases}
        \C^{2|2} & \mu \to \lambda, \ell(\mu)\ell(\lambda) \mod 2 = 1\\
        \C^{1 | 1} & \mu \to \lambda, \ell(\mu)\ell(\lambda) \mod 2 = 0\\
        0 & \mu \not\to \lambda
    \end{cases}.
    \]

    All of this demonstrates that the theorem holds at level $N + 1$, with the appropriate labelings of the irreducible super-representations, and with the appropriate restriction rules.

    An exactly similar argument works for the even parts, only without having to worry about gradings, after the base cases are established. For $n = 0$ and $n = 1$, $CW_n[0]= \C$. For $n = 2$, it is not true that $(CW_2 e_2 CW_2)[0]  = CW_2[0] e_2 CW_2[0]$. Nonetheless, for all but finitely many $q$ we still have that $CW_2[0]$ is semisimple. Moreover, $(CW_2 e_2 CW_2)[0]$ is a a two sided ideal of $CW_2[0]$, and it is easy to check that $\dim (CW_2 e_2 CW_2)[0] = 2$. The two idemotents $p_{(0), 1} = e_n$ and $p_{(0), -1} = v_n e_n v_n$ are easily seen to be orthogonal. As a result, $ p_{(0), 1} + p_{(0), -1}$ is the minimal central idempotent projecting onto $(CW_2 e_2 CW_2)[0]$. We have gives $CW_2[0]/(CW_2 e_2 CW_2)[0] \cong G_2(q)[0]$. Since $G_2(q)[0]\cong M_2(\C)$, it follows that we have one more minimal idempotent in $CW_2[0]$, $p_{(2)}$, appearing with multiplicity $2$ in $CW_2$. Restriction of the representation $W_{(2)}^{(2)}$ to $CW_1$ splits into two copies of $W_{(1)}^{(1)}$. Letting $p_{(2)}$ and $p_{(2)}'$ be orthogonal idempotents with $p_{(2)} + p_{(2)}' = z_{(2)}$ the minimal central idempotent projecting onto $G_2(q)$, we have a partition of unity:
    \[
    p_{(0), 1} + p_{(0), -1} + p_{(2)} + p_{(2)}' = 1
    \]
    Moreover, since $CW_1[0] = \C$, it is now clear that the Bratelli graph for the inclusion $CW_1[0] \sbs CW_2[0]$ is 
\begin{tikzcd}
	&& {(1)} \\
	{(0)^+} & {(0)^-} & {(2)}
	\arrow[no head, from=1-3, to=2-1]
	\arrow[no head, from=1-3, to=2-2]
	\arrow[equals,  from=1-3, to=2-3]
\end{tikzcd}

Since $CW_{n + 1}[0]e_n CW_{n + 1}[0] = CW_n[0] e_n CW_n[0]$, we may use the Jones basic construction to inductively build the Bratelli diagram for $CW_n[0]$, verifying that the restriction rules stated in the theorem hold. 

The restriction rules for $G_n(q)[0] \sbs G_n(q)$ transfer directly to the restriction rules for $CW_n[0] \sbs CW_n$ for all but finitely many $q$. 
\end{proof}

\subsection{Limiting Fusion Rules}

Now we use the results of the last subsection to prove the existence of a limiting fusion ring for all of the adjoint conformal embeddings (and the type C one from level rank duality). For the type C series, this is especially easy. For types B and D, our stabilizing dimension formulas allow us to show that, after semisimplifications, the extra vertex generators eventually have no effect on the endomorphism algebras. 

\begin{theorem}
    For each $r \in \N$ let $C_r\in \set{\fanC(\so_{N}, \check{h}), \fanC(\fsp_2r, r + 1), \fanC(\fsp_{2r + 2}, r)}$ and for $A_r$ the Etale algebra coming from the conformal embeddings in the table below. Let $\gamma_r$ be the half braid on the appropriate fermion object in every case. Let $\sigma$ and $\epsilon$ and $q$ be appropriately chosen as in the table:
\begin{center}
\begin{tabular}{c|c|c|c}
    Conformal Embedding & $q$ & $\epsilon$ & $\sigma$ \\\hline
    $\V(\fsp_{2r}, r + 1) \sbs \V(\so_{\dim \fsp_{2r}}, 1)$ & $e^{\frac{2\pi i}{4(2r + 2)}}$ & $1$ & $-1$  \\\hline
    $\V(\so_{2r + 1}, 2r - 1) \sbs \V(\so_{\dim \so_{2r + 1}}, 1)$ & $e^\frac{2\pi i}{4(4r - 2)}$ & $-1$ & $1$ \\\hline
    $\V(\so_{2r}, 2r - 2) \sbs \V(\so_{\dim \so_{2r}}, 1)$ & $e^{\frac{2\pi i}{4(4r - 4)}}$& $-1$& $1$   \\\hline
    $\V(\fsp_{2r + 2}, 2r) \sbs \V(\so_{\dim \fsp_{2r}}, 1)$ &$e^{\frac{2\pi i}{4(2r + 2)}}$ & $-1$ & $-1$ 
\end{tabular}
\end{center}
Then for every $n \in \N$ there exists some $r \in \N$ such that for all $m \le n$ we have 
    \[
    End_{C_A/\gamma}(X^{\tensor m}) \cong CW_m(\sigma, \epsilon).
    \]
The type $Q$ simple objects $V_\lambda$ appearing in tensor powers of $X^{\tensor m}$ for $m \le n$ are indexed by strict partitions of size $\le n$, where $V_\lambda$ is type $Q$ if $\ell(\lambda)$ is odd and type $M$ if $\ell(\lambda)$ is even. Moreover, for all $\lambda, \mu$ with $|\lambda|, |\mu|< n$, we have 
\[
Hom_{C_A/\gamma}(V_\lambda \tensor V_{(1)} \to V_\mu)
\cong \begin{cases}
        \C^{2|2} & \mu \leftrightarrow \lambda, \ell(\mu)\ell(\lambda) \mod 2 = 1\\
        \C^{1 | 1} & \mu \leftrightarrow \lambda, \ell(\mu)\ell(\lambda) \mod 2 = 0\\
        0 & \mu \not\leftrightarrow \lambda
    \end{cases}.
\]
Similarly, $End_{{C_r}_{A_r}/\gamma}(X^{\tensor m}) \cong CW_m(\sigma, \epsilon).$
\end{theorem}
End
\begin{proof}
    By Lemma~\ref{stableDimB} and Lemma~\ref{stableDimD} and by Theorem~\ref{dimensionOfCAThypeC} we have 
    \[
    \lim_{r \to \infty}\dim End_{C_r{_{A_r}}}(X\tensor k) = 2^{k - 1} (2k - 1)!!.
    \]

    As a result, we can choose $r$ large enough so that for all $m \le n$ we have $\dim End_{C_r{_{A_r}}}(X^{\tensor k}) = 2^{m - 1} (2m - 1)!!$. Since $X$ is a type $Q$ simple in every case, 
    \[
    \dim End_{C_r{_{A_r}}/\gamma}(X^{\tensor k}) = 2\dim End_{C_r{_{A_r}}}(X^{\tensor k}) = 2^{m} (2m - 1)!!
    \]
    for all $m \le n$. 

    For every $m \le n$, the Gram matrix for the categorical trace on $CW_n(\sigma, \epsilon)$ has only finitely many zeros. Choose $r$ large enough so that none of the relevant specializations are zeros of any of the Gram matrices for $m \le n$. 

    With this choice of $r$, $\dim \sms{CW_m(\sigma, \epsilon)} = 2^m (2m - 1)!!$ for all $m \le n$ and $\sms{CW_m(\sigma, \epsilon)} = CW_m(\sigma, \epsilon)$. As a result, since $\sms{CW_m(\sigma, \epsilon)}$ injects into $End_{C_r{_{A_r}}/\gamma}(X\tensor k)$ which has the same dimension, $CW_m(\sigma, \epsilon)\cong  End_{C_r{_{A_r}}/\gamma}(X\tensor k)$. Taking even parts, $End_{{C_r}_{A_r}}(X^{\tensor m}) \cong CW_m(\sigma, \epsilon)[0]$.

    Every minimal even idempotent $p_\lambda$ in $End_{C_r{_{A_r}}/\gamma}(X\tensor k)$ corresponds to a simple object (of type M if $p_\lambda$ is minimal, and type Q otherwise.) As a result, by Theorem~\ref{RepsOfCW} for each $m \le n$, the simple objects $V_\lambda$ contained in $X^{\tensor m}$ in ${C_r}_{A_r}/\gamma$ are indexed by strict partitions $\lambda \dashv_s m - 2f$, with $V_\lambda$ a type $Q$ if $\ell(\lambda)$ is odd and a type $M$ otherwise. If we had an isomorphism $V_\lambda \cong V_\mu$ in ${C_r}_{A_r}$ for distinct partitions, then the matrix blocks of the corresponding irreducible supermodules would collapse into a larger block. Therefore distinct strict partitions of $\lambda \dashv_s m - 2f$ correspond to non-isomorphic simples. 

    Similarly, distinct partitions of $\lambda \dashv_s m - (2f - 1)$ label non-isomorphic simples. 

    As a result, it follows that distinct strict partitions of size $\le n$ of the same parity correspond to non-isomorphic simple objects. It remains to be seen that there are no isomorphisms $\phi: V_\lambda \to V_\mu$ whenever $|\lambda|, |\mu| \le n$ and $|\lambda| - |\mu| \equiv 1 \pmod{2}$.

    If $C_r$ is of type $D$ or $C$, the categories ${C_r}_{A_r}$ are $\Z/2\Z$-graded, as is the fermion condensation ${C_r}_{A_r}/\gamma$, because the weight correspondng to every abelian ideal in these cases is zero graded. Therefore, the only potentially problematic case is type B. By Lemma~\ref{eventuallyZ2Graded} we may choose $r$ large enough such that for all $m \le n$ there are no non-zero morphisms 
    \[
    \phi\in Hom_{{C_r}_{A_r}}(X^{\tensor 2m + 1} \to A).
    \]
    As a result, in the fermion condensation, there are no non-zero even morphisms $X^{\tensor 2m + 1} \to A$ whenever $m \le n$. Because $X$ is a type $Q$ simple, it follows that the even and odd hom-spaces $X^{\tensor 2m + 1} \to A\in {CA_r}_{A_r}/\gamma$ are isomorphic as vector spaces by precomposing with the odd automorphism $v\tensor I_X^{\tensor 2m}$. As a result $Hom_{{C_r}_{A_r/\gamma}}(X^{\tensor 2m + 1} \to A) = 0$

    Suppose now by way of contradiction that for distinct strict partitions $\lambda \dashv_s m - 2f$ and $\mu \dashv_s m - 2k - 1$, we have an isomorphism $\phi: V_\lambda \to V_\mu$. Since $V_\lambda \sbs X^{\tensor m - 2f}$ and $V_\mu \sbs X^{\tensor n - 2k - 1}$, by semismiplicity we have the inclusion $\iota_\lambda : V_\lambda \hookrightarrow X^{\tensor m - 2f}$  and a projection $\pi_\lambda X^{m - 2f} \to V_\lambda$ such that $\pi_\lambda \circ \iota_\lambda = I_{V_\lambda}$, and similarly for $V_\mu$. We can then form the composition $\iota_\mu \circ \phi \circ \pi_\lambda: X^{\tensor m - 2f} \to X^{m - 2k - 1}$. This is non-zero, as we can precompose with $\iota_\lambda$ and postcompose with $\pi_\mu$ to recover $\phi$. 

    But by self duality of $X$, this gives a non-zero morphism in $Hom_{{C_r}_{A_r} /\gamma}(X^{\tensor 2(m -f + k) - 1}  \to A)$. Since $m - f + k \le n$, this contradicts our choice of $r$. 

    Thus far we have proven the consistency of the labeling of simple objects in the fermion condensation appearing in $X^{\tensor m}$ for $m \le n$ by strict partitions of size no greater than $n$. A nearly identical (but somewhat easier) argument shows that the simple objects appearing in ${C_r}_{A_r}$ appearing in $X^{\tensor m}$ for $n \le m$ are indexed by signed strict partitions of size no greater than $n$. 

    The branching rules for both ${C_r}_{A_r}$ and ${C_r}_{A_r}/\gamma$ are a direct consequence of the fact that the inclusion $CW_n(\sigma, \epsilon) \sbs CW_{n + 1}(\sigma, \epsilon)$ and $CW_n(\sigma, \epsilon)[0]\sbs CW_{n + 1}(\sigma, \epsilon)[0]$ is given by padding a Clifford BMW diagram with a single strand, which in turn corresponds to tensoring with $X$.  
\end{proof}

\begin{corollary}
    For $V_{(\lambda, \kappa)}$ a simple object indexed by a signed strict partition with $\kappa\in \set{(-1)^{\ell(\lambda) + 1}, 1}$ in an appropriate $C_A$ with $r$ large enough relative to $|\lambda|$, and with $g\in C_A$ the fermion object, we have
    \[
    g \tensor V_{(\lambda, \kappa)} \cong V_{(\lambda, (-1)^{\ell(\lambda) + 1}\kappa)}
    \]
\end{corollary}
\begin{proof}
    By \cite{CainAndHans}, Lemma~3.16, it follows that for $q$ a $4N$-th root of $1$, we have that the irreducible representations of $G_n(q)$ and $G_n(q)[0]$ are indexed by strict (signed) partitions when $n < N$. Since for each of the relevant categories, $q$ is a primitive $4N$-th root of unity, with $N \to \infty$ as $r\to \infty$.

    As a result, for $\ell(\lambda)$ odd, a minimal idempotent $p_\lambda \in End_{C_A}(X^{\tensor m})\cong CW_n(\sigma, \epsilon)$ with $p_\lambda CW_n(\sigma, \epsilon)\cong W_\lambda$, we have by the restriction rules for $CW_n[0] \sbs CW_n$ that $p_\lambda$ splits as $p_\lambda = p_{\lambda, 1} + p_{\lambda, -1}$ for $p_{\lambda, \pm 1}$ inhomogeneous minimal idempotents in $CW_n$. As a result, $p_\lambda$ is the even projection onto a type $Q$ simple $X_\lambda$ in $C_A/\gamma$. Hence $X_\lambda \tensor g = X_\lambda$.

    A similar argument works if $\ell(\lambda)$ is even, only showing that $X_\lambda \tensor g$ is not isomorphic to $X_\lambda$.
\end{proof}

As a result, just as with the Hecke Clifford category of \cite{CainAndSnyder}, tensoring by $g$ acts as a ``parity swap'' operator. 

With a little bit of work, one can show the following:
\begin{itemize}
    \item There exists a limiting Grothendiek ring for $C_A$ for all conformal embeddings of Lie types B C and D into Lie types B and D.
\end{itemize}

We will not prove this here, although it is a fairly straightforward consequence of the results of this paper. We have some partial progress for the explicit specializations at roots of unity, which we will save for future works. We hope to eventually classify all simple objects and all fusion rules for the quantum subgroups $C_A$ of types B C and D coming from conformal embeddings.

\end{document}